\documentclass[12pt,leqno,fleqn,epsfig]{article}
\usepackage{amssymb, epsfig, amsmath, amsthm}
\usepackage{mathrsfs}       
\usepackage{color}  
     
\newcommand{\R}{\mathbb{R}}

\newcommand{\Z}{\mathbb{Z}}

\newcommand{\N}{\mathbb{N}}

\newcommand{\osc}{\operatorname*{osc}}

\newcommand{\supp}{\operatorname*{supp}}

\newcommand{\const}{\operatorname*{const}}

\newcommand{\bb}{\begin{equation}}
\newcommand{\ee}{\end{equation}}
\newcommand{\bq}{\begin{eqnarray}}
\newcommand{\eq}{\end{eqnarray}}
\newcommand{\bqn}{\begin{eqnarray*}}
\newcommand{\eqn}{\end{eqnarray*}}

\newcommand{\var}{\varepsilon}

\newcommand{\intl}{\int\limits}

\newcommand{\Beweisende}{\rule{0.2cm}{0.2cm}}

\newcommand{\D}{\displaystyle}

\newcommand{\intmw}{{\int\hspace{-830000sp}-\!\!}}

\newcounter{secnum}

\newtheorem{thm}{Theorem}[section]

\newtheorem{cor}[thm]{Corollary}
\newtheorem{lem}[thm]{Lemma}

\theoremstyle{definition}

\newtheorem{rem}[thm]{Remark}

\title{Transport equation  in generalized Campanato spaces}  %Liouville problem for the steady 3D Navier-Stokes equations}
 
\author{Dongho Chae$^{(*)}$ and J\"{o}rg Wolf$^{(\dagger)}$\\
\ \\
 Department of Mathematics$^{(*), (\dagger)}$ \\
Chung-Ang University\\
Dongjak-gu Heukseok-ro 84\\
Seoul 06974, Republic of Korea\\
and \\
School of Mathematics$^{(*)}$ \\
Korea Institute for Advanced Study\\
Dongdaemun-gu   Hoegi-ro 85 \\
Seoul 02455, Republic of Korea\\
$^{(*)}$e-mail: dchae@cau.ac.kr \\
$^{(\dagger)}$e-mail: jwolf2603@cau.ac.kr}
\date{}
\begin{document}
\maketitle
\begin{abstract}
In this paper we study the  transport equation in $ \R^{n} \times  (0,T)$, $T >0$,
\[
\partial _t f + v\cdot \nabla f = g, \quad  f(\cdot ,0)= f_0 \quad   \text{in}\quad  \R^{n}
\]
in  generalized Campanato spaces $ \mathscr {L}^{ s}_{ q(p, N)}(\mathbb {R}^{n})$. The critical case is particularly interesting, and is 
applied to the local well-posedness problem in a space close to the Lipschitz space in our companion paper\cite{cw}.
More specifically, in the critical case  $s=q=N=1$
we have the embedding relations,   $ B^{1}_{\infty, 1} (\Bbb R^n)  \hookrightarrow \mathscr {L}^{ 1}_{ 1(p, 1)}(\mathbb {R}^{n}) \hookrightarrow C^{0, 1} (\Bbb R^n)$, 
where  $B^{1}_{\infty, 1} (\Bbb R^n)$ and $ C^{0, 1} (\Bbb R^n)$ are the Besov space and the Lipschitz space respectively. 
For  $f_0\in   \mathscr {L}^{ 1}_{ 1(p, 1)}(\mathbb {R}^{n})$, $v\in L^1(0,T;  \mathscr {L}^{ 1}_{ 1(p, 1)}(\mathbb {R}^{n}))),$  and  
$ g\in L^1(0,T;  \mathscr {L}^{ 1}_{ 1(p, 1)}(\mathbb {R}^{n})))$, 
we prove the existence and uniqueness 
of solutions to the transport equation in $ L^\infty(0,T;  \mathscr {L}^{ 1}_{ 1(p, 1)}(\mathbb {R}^{n}))$ such that
\[
\hspace*{-1cm}\|f\|_{ L^\infty(0,T;  \mathscr {L}^{ 1}_{ 1(p, 1)}(\mathbb {R}^{n})))} \le C 
\Big( \|v\|_{ L^1(0,T;  \mathscr {L}^{ 1}_{ 1(p, 1)}(\mathbb {R}^{n})))}, \|g\|_{ L^1(0,T;  \mathscr {L}^{ 1}_{ 1(p, 1)}(\mathbb {R}^{n})))}\Big). 
\]
Similar results in the other cases are also proved.
  \noindent{\bf
keywords:} transport equation, generalized Campanato space, well-posedness

\end{abstract}%\renewcommand{\thesection}{\mbox{\arabic{section}}.}
\newpage
\tableofcontents

%%% ----------------------------------------------------------------------
%       SECTION 1
%%% ----------------------------------------------------------------------
\section{Introduction}
\label{sec:-1}
\setcounter{secnum}{\value{section} \setcounter{equation}{0}
\renewcommand{\theequation}{\mbox{\arabic{secnum}.\arabic{equation}}}}

Let $ 0< T < +\infty$ and $Q = \R^{n}\times  \R_+ $ with $n\in \Bbb N, n\geq 2$. We consider the 
transport equation  
\begin{equation}
\begin{cases}
\partial _t f+ (v\cdot \nabla) f = g \quad  \text{ in}\quad Q,
\\[0.3cm]
v= v_0\quad  \text{ on}\quad  \R^{n}\times \{ 0\},
\end{cases}
\label{trans}
\end{equation}
where $ f=f(x_1, \ldots, x_n)$ is unknown, while $v=(v_1, \cdots, v_n)=v(x,t)$ represents a given drift velocity  
and $ g=g(x_1, \ldots, x_n)$ given function. 
 
Our aim in this paper is to obtain estimates of solutions to  \eqref{trans} in  generalized Campanato 
spaces.   
The proof relies on a key estimate in terms of local oscillation. 
As byproduct we get existence of solutions in Besov spases and Tribel-Lizorkin spases, which can be estimated by the data 
belonging to these spaces.  One of the  main motivations to study the transport equation in such generalized Campanato spaces  is to apply it to prove  local well-posedness  of the incompressible Euler equations in function space embedded in 
the Lipschitz space, which includes linearly growing functions at spatial infinity. For recent developments of the local well-posedness/ill-posedness of the Euler equations in  various critical function spaces embedded   in $C^{0,1} (\Bbb R^n)$ we refer \cite{bour1, bour2, lem, pak, vis1, vis2, bah, che, lio}). We would also like to refer \cite{dip} for  the study of transport equation with drift velocity in less regular space. For the  application of our new function spaces in the critical case  to the Euler equations  please see our  companion paper \cite{cw}.

\vspace{0.3cm}
Let us introduce the function spaces we will use throughout the paper. Let $ N\in \N \cup  \{0, -1 \}$.   By $ \mathcal{P}_N$  ($ \dot{\mathcal{P}} _N$ respectively) we denote the space of all polynomial (all homogenous polynomials respectively)
of degree less or equal $ N$.  We equip the space $ \mathcal{P}_N$ with the norm  $ \| P\|_{ (p)}= \| P\|_{ L^p(B(1))}$. Note that 
since $ \dim (\mathcal{P}_N)<+\infty$ all norms $ \| \cdot \|_{ (p)}, 1 \le p \le  \infty$,  are equivalent.  For notational convenience,  in case $ N=-1$ we  use the convention $\mathcal{P}_{ -1} = \{0\}$, which consists of  the trivial polynomial $ P \equiv 0$.

\vspace{0.5cm}  
 Let $ f\in L^p_{ loc}(\R^{n}), 1 \le p \le +\infty$. For  $ x_0\in  \R^{n}$ and $ 0< r<\infty$ we define 
the oscillation 
\begin{equation}\label{deff}
\osc_{ p, N} (f; x_0, r):= | B(r)|^{ - \frac{1}{p}} \inf_{ P\in \mathcal{P}_N}  \| f-P\|_{ L^p(B(x_0, r))}. 
\end{equation}
We note that from our convention above in the case $ N=-1$ we  have 
\[
\osc_{ p, -1} (f; x_0, r):= | B(r)|^{ - \frac{1}{p}}   \| f\|_{ L^p(B(x_0, r))}. 
\]

Then, we define for $ 1 \le q, p \le +\infty$ and $ s\in (-\infty, N+1]$  the spaces
$$
 \mathscr{L}^{ s}_{q (p, N)}(\R^{n})= \left\{ f\in L^p_{ loc}(\R^{n}) \,\Big|\,|f|_{  \mathscr{L}^{ s}_{q (p, N)}} 
  :=  \bigg\|\left( \sum_{j\in \Bbb Z} \Big(2^{ -sj}\osc_{ p, N} (f; \cdot , 2^{j})\Big)^q\right)^{\frac1q}\bigg\|_{ L^\infty} <+\infty\right\}.
$$
%where  $ |f|_{  \mathscr{L}^{ s}_{q (p, N)}}$ stands for the semi norm 
%\[
%  |f|_{  \mathscr{L}^{ s}_{q (p, N)}} 
%  = \supl_{ x_0\in \R^{n}} \left( \sum_{j\in \Bbb Z} 2^{ -sqj}|\osc_{ p, N} (f; x_0, 2^{j})|^q\right)^{\frac1q} .
%\]
Furthermore, by $ \mathscr{L}^{k, s}_{q (p, N)}(\R^{n}) $, $k\in \N$,  we denote the space of all $ f\in W^{k,\, p}_{ loc}(\R^{n})$ such that $ D^k f\in  \mathscr{L}^{ s}_{q (p, N)}(\R^{n})$.  The space $ \mathscr{L}^{k, s}_{q (p, N)}(\R^{n})$ will be equiped with the  norm   
\[
\| f\|_{ \mathscr{L}^{k, s}_{q (p, N)} }= | D^k f|_{ \mathscr{L}^{  s}_{q (p, N)} } + \| f\|_{ L^p(B(1))},\quad  f\in 
{ \mathscr{L}^{k, s}_{q (p, N)}}(\R^{n}). 
\]
According to the characterization theorem of the Triebel-Lizorkin spaces  in terms of  oscillation, we have
\begin{align*}
\hspace{-.5in}\begin{cases}
{\D f\in F^s_{ r, q} (\Bbb R^n) \quad  \Leftrightarrow   \quad  \|f\|_{ L^{ \min\{r,q\}}}+
\quad  \bigg\|\left( \sum_{j=-\infty} ^0\Big(2^{ -sj}\osc_{ p, N} (f; \cdot , 2^{j})\Big)^q\right)^{\frac1q}\bigg\|_{ L^r} <+\infty. }
\\[0.5cm]
0< r<+\infty, 0< q \le \infty, \quad s> \Big(\frac{1}{r} - \frac{1}{p}\Big)_{ +},\quad  s> \Big(\frac{1}{q} - \frac{1}{p}\Big)_{ +}
\end{cases}
\end{align*}
(cf. \cite[Theorem, Chap. 1.7.3]{tri}), and we could regard the spaces $  \mathscr{L}^{s}_{q (p, N)}(\R^{n} ) $ as an extension of the  limit case of 
$F_{ r,q}^s(\R^{n} ) $
as $ r\rightarrow +\infty$.

In fact in case $ q=+\infty$ and $ s >0$ we get the usual Campanato spaces with the isomorphism relation(cf. \cite{ca, gia})
\[
\mathscr{L}^{n+ps, p}_N(\R^{n}) \cong  \mathscr{L}^{s}_{\infty (p, N)}(\R^{n}).
\] 
 Furthermore, in the case $ N=0, s=0$ and $ q= \infty$ we get the space of bounded 
mean oscillation, i.e.,
\[
\mathscr{L}^{0}_{\infty (p, 0)}(\R^{n})\cong BMO. 
\]
In case $ N=-1$ and $ s\in (-\frac{n}{p}, 0) $ the above space coincides with the usual Morrey space $ {\cal M}^{ n+ ps}(\R^{n} )$. 

We note that the oscillation introduced in \eqref{deff} is attained  by a unique 
polynomial $ P_{ \ast}\in \mathcal{P}_N$.

\vspace{0.2cm}
According to Theorem\,\ref{thm5.4} (see Section\,3 below), for the spaces $ \mathscr{L}^{1}_{1 (p.  1)}( \R^{n})$ we have the following embedding properties
\begin{equation}
B^{ 1+ \frac{n}{r}}_{ r, 1} \hookrightarrow \mathscr{L}^{1}_{1 (p, 1)}( \R^{n}) \hookrightarrow  \mathscr{L}^{0,1}_{1 (p, 0)}( \R^{n}) 
 \hookrightarrow C^{0, 1}(\R^{n}). 
\label{1.8}
\end{equation}
Accordingly, 
\begin{equation}
\| \nabla u\|_{\infty}  \le c\| u\|_{\mathscr{L}^{1}_{1 (p, 1)}}.
\label{1.11}
\end{equation}

Furthermore,  for every $ f\in  \mathscr{L}^{k}_{1 (p, k)}( \R^{n}), k\in \{0,1\}$, there exists a unique $ \dot{P}^k_{ \infty}(f)\in \dot{\mathcal{P}} _1$,  such that for all $ x_0 \in \R^n$
\[
f \quad  \text{converge asymptotically to} \quad  \dot{P}^k_{ \infty}(f)\quad  \text{as }\quad  |x| \rightarrow +\infty.
\]
The  precise meaning of this asymptotic limit will be given in Section 3 below.

We are now in a position to present our first main result.
\begin{thm}[The case $ N=0$]
\label{thm1} 
Let $ 0< T< +\infty$. Let $s\in ( - \frac{n}{q}, 0), 1 <p< +\infty, 1 \le q \le +\infty$.  
Let  $v\in L^1(0,T; L^p_{ loc} (\R^{n}))$, with 
\begin{equation}
 \intl_{0}^T \|\nabla v(\tau )\|_{ \infty} d\tau  < +\infty. 
\label{cond1}
 \end{equation} 
Then for every  $ f_0\in \mathscr{L}^{ s}_{q (p, 0)} (\R^{n})$ and  $g\in L^1(0,T; \mathscr{L}^{s}_{q (p, 0)} (\R^{n}))$ there exists a unique solution $ f\in L^\infty(0, T; \mathscr{L}^{ s}_{q (p, 0)} (\R^{n}))$
to the transport equation  \eqref{trans}.  Furthermore, it holds for almost all $ t\in (0, T)$
\begin{align}
| f(t)|_{\mathscr{L}^{ s}_{q (p, 0)}}
\le   c\bigg\{| f_0|_{\mathscr{L}^{ s}_{q (p, 0)}}+ \intl_{0}^{T} | g(\tau )|_{\mathscr{L}^{ s}_{  q(p, 0)} } d\tau \bigg\}
 \exp \bigg(c\intl_{0}^{T} \| \nabla v(\tau )\|_{ \infty}  d\tau\bigg).
\label{estimate1}
\end{align}  
 
\end{thm}

In case $ N=1$ we get 

\begin{thm}[The case $ N=1$ and $ s=1$]
\label{thm2} 
Let $ 0< T< +\infty$ and $ 1 <p< +\infty, 1 \le q \le +\infty$.  
Let  $v\in L^1(0,T; \mathscr{L}^{ 1}_{q (p, 1)} (\R^{n}))$ with  \eqref{cond1} and 
\begin{equation}
 \intl_{0}^T \sup_{ x_0\in \R^{n}}\Big(\sum_{j=-\infty}^{0}  (-j)^{ q-1} 2^{ -jq} \osc_{ p,1}(v(\tau ); x_0, 2^j ) \Big)^{ \frac{1}{q}}d\tau  < +\infty
\label{cond2}
 \end{equation} 
Let $ f_0\in \mathscr{L}^{ 1}_{q (p, 1)} (\R^{n})$ and  $g\in L^1(0,T; \mathscr{L}^{1}_{q (p, 1)} (\R^{n}))$ 
satisfying the condition 
\begin{align}
 & \qquad \sup_{ x_0\in \R^{n}}  \osc_{p,0}(f_0; x_0, 1 ) +\intl_{0}^T  \sup_{ x_0\in \R^{n}}  \osc_{p,0}(g(\tau ); x_0, 1 ) d\tau   < +\infty.
 \label{cond3}
\end{align}
Then, there exists a unique solution $ f\in L^\infty(0, T; \mathscr{L}^{ 1}_{q (p, 1)} (\R^{n}))$
to the transport equation  \eqref{trans}.  Furthermore, it holds for  all $ t\in (0, T)$
\begin{align}
| f(t)|_{\mathscr{\tilde{L}}^{ 1}_{q (p, 1)}}
\le   c\bigg\{| f_0|_{\mathscr{\tilde{L}}^{ 1}_{q (p, 1)}}+ \intl_{0}^{T} 
| g(\tau )|_{\mathscr{\tilde{L}}^{ 1}_{  q(p, 1)} } d\tau \bigg\}
 \exp \bigg(c\intl_{0}^{T} C(\tau ) d\tau\bigg),
\label{estimate2}
\end{align}  
where we set
\begin{align*}
  C(\tau ) =  \| \nabla v(\tau )\|_{ \infty} + 
 \sup_{ x_0\in \R^{n}}\Big(\sum_{j\in \Z}  (j^-)^{ q-1} 2^{ -jq} \osc_{ p,1}(v(\tau ); x_0, 2^j) \Big)^{ \frac{1}{q}}, 
\end{align*}
 $ j ^- = -\min\{j,0\}$, and  $ | z|_{\mathscr{\tilde{L}}^{ 1}_{  q(p, 0)} }$ stands for the semi norm 
 \[
 | z|_{\mathscr{\tilde{L}}^{ 1}_{  q(p, 1)} }=  | z|_{\mathscr{L}^{ 1}_{  q(p, 1)} } +\sup_{ x_0\in \R^{n}} 
 |\nabla \dot{P}^1_{ x_0, 1}(z)|. 
 \]

\end{thm}

Our third  main result concerns the case $s>1$. 
\begin{thm}[The case $ N \ge 1$ and $ s>1$]
\label{thm3} 
Let $ 0< T< +\infty$,  $ N\in \N,  1 < s < +\infty, 1 <p< +\infty,$ and $ 1 \le q \le +\infty$.  
Let  $v\in L^1(0,T; \mathscr{L}^{ s}_{q (p, N)} (\R^{n}))$ with  \eqref{cond1} and 
Let $ f_0\in \mathscr{L}^{ s}_{q (p, N)} (\R^{n}) $ and  $g\in L^1(0,T; \mathscr{L}^{s}_{q (p, N)} (\R^{n}))$ 
satisfying the condition 
\begin{align}
 & \|\nabla f_0\|_{ \infty} +\intl_{0}^T \| \nabla g(\tau )\|_{ \infty}d\tau   < +\infty.
 \label{cond4}
\end{align}
Then, there exists a unique solution $ f\in L^\infty(0, T; \mathscr{L}^{ 1}_{q (p, 1)} (\R^{n})) $
to the transport equation  \eqref{trans} together with the estimate 
\begin{align}
| f(t)|_{\mathscr{L}^{ s}_{q (p, 0)}}
\le   c\bigg\{| f_0|_{\mathscr{\tilde{L}}^{ s}_{q (p, 0)}}+ \intl_{0}^{T} | g(\tau )|_{\mathscr{\tilde{L}}^{ s}_{  q(p, 0)} }\bigg\}
 \exp \bigg(c\intl_{0}^{T} \| v(\tau )\|_{ \mathscr{\tilde{L}}^{ s}_{  q(p, 0)} }  d\tau\bigg),
\label{estimate3}
\end{align}  
  where $ | z|_{\mathscr{\tilde{L}}^{ s}_{  q(p, 0)} }$ stands for the semi norm defined by
 \[
 | z|_{\mathscr{\tilde{L}}^{ s}_{  q(p, 0)} }=  | z|_{\mathscr{L}^{ s}_{  q(p, 0)} } +
 \|\nabla z\|_{ \infty}. 
 \]
\end{thm}

From Theorem\,\ref{thm2} we get the following corollary for the special case $ s=q=N=1$, which will be useful for our future application to the Euler equations in the critical spaces.

\begin{cor}
 \label{cor1.4}
Let $ 0< T< +\infty,  1 <p< +\infty$.  
Let  $v\in L^1(0,T; \mathscr{L}^{ 1}_{1 (p, 1)} (\R^{n}))$, 
$ f_0\in \mathscr{L}^{ 1}_{1 (p, 1)} (\R^{n})$ and  $g\in L^1(0,T; \mathscr{L}^{1}_{1 (p, 1)} (\R^{n}))$. 
Then there exists a unique solution $ f\in L^\infty(0, T; \mathscr{L}^{ 1}_{1 (p, 1)} (\R^{n}))$
to the transport equation  \eqref{trans}.  Furthermore, it holds for all $ t\in (0, T)$
\begin{align}
\| f(t)\|_{\mathscr{L}^{ 1}_{1 (p, 1)}}
\le   C\bigg\{1+ \intl_{0}^{T} | v(\tau )|_{\mathscr{L}^{ 1}_{ 1(p, 1)} } d\tau \bigg\}
\exp \bigg(c\intl_{0}^{T} \| \nabla v(\tau )\|_{ \infty}  d\tau\bigg).
\label{estimate4}
\end{align}  
where 
\[
C= c\bigg(\| f_0\|_{\mathscr{L}^{ 1}_{1 (p, 1)}}+\intl_{0}^{T} \| g(\tau )\|_{\mathscr{L}^{ 1}_{ 1(p, 1)} } d\tau \bigg),
\]
while $ c=\const>0$ depending on $ n$ and $p$. 
 \end{cor}

\begin{rem}
\label{rem1.7}
Using the well-known characterization  of $ B^1_{ \infty, 1}(\R^{n} )$ in terms of oscillation, 
we easily verify the embeddings
\begin{equation}
  B^1_{ \infty, 1}(\R^{n} ) \hookrightarrow  \mathscr{L}^{1}_{1 (p, 1)}(\R^{n})\cap L^\infty(\R^{n} )
  \hookrightarrow  \mathscr{L}^{1}_{1 (p, 1)}(\R^{n}). 
\label{emb}
 \end{equation}  
Indeed, referring   to  \cite[Theorem, Chap.1.7.3]{tri}), we see that 
\[
 v\in   B^1_{ \infty, 1}(\R^{n}) \quad  \Leftrightarrow \quad 
  \sum_{j=-\infty}^{0} 2^{-j} \|\osc_{p, 1}(v; \cdot , 2^j)\|_{L^\infty} + \|v\|_{L^\infty } < 
  +\infty. 
\]

This shows that for $ x\in \R^{n} $ it holds
\begin{align*}
 &\sum_{j\in \Z} 2^{-j} \osc_{p, 1}(v; x , 2^j)
\\
 &\quad  \le\sum_{j=-\infty}^{0} 2^{-j} \|\osc_{p, 1}(v; \cdot , 2^j)\|_{L^\infty } + 
   \sum_{j=1}^\infty 2^{-j} \osc_{p, 1}(v; x , 2^j)  +
  \|v\|_{L^\infty }. 
\end{align*}
On the other hand, it is readily seen that ${\D \osc_{p, 1}(v; x , 2^j)} \le 2\|v\|_{ L^\infty}$. Accordingly, 
the second sum on the right-hand side is bounded by $ \|v\|_{ L^\infty}$. This yields 
\begin{align*}
\|v\|_{\mathscr{L}^{1}_{1 (p, 1)} } & =   \sum_{j\in \Z} 2^{-j} \osc_{p, 1}(v; x , 2^j)  + \|v\|_{ L^2(B(1))}
\\ 
&\le \sum_{j\in \Z} 2^{-j} \osc_{p, 1}(v; x , 2^j) +c\|v\|_{ B^{ 1}_{ \infty,1}}
\\
&\le   c\sum_{j=-\infty}^{0} 2^{-j} \|\osc_{p, 1}(v; \cdot , 2^j)\|_{L^\infty} + c\|v\|_{L^\infty } \le c\|v\|_{ B^{ 1}_{ \infty,1}}.
\end{align*}
Secondly, according to \cite[p. 85]{tay} (see also \cite{bah}) we have the embedding 
\[
 B^1_{ \infty,1}(\R^{n} )  \hookrightarrow C^1(\R^{n} ) \cap L^\infty(\R^{n} ). 
\]

On the other hand, there exists a function $ f\in \mathscr{L}^{1}_{1 (p, 1)} (\R^{n} )$ 
which is not in $ C^1(\R^{n} )$ (see Appendix\,B). This clearly shows that $ \mathscr{L}^{1}_{1 (p, 1)} (\R^{n} )$ contains less regular functions then  
$  B^1_{ \infty,1}(\R^{n} ) $.

\vspace{0.2cm}
Thirdly, since $ \mathscr{L}^{1}_{1 (p, 1)} (\R^{n} )$ contains linearly  growing function at infinity, in particular polynomials of of degree 
less or equal one,  $ \mathscr{L}^{1}_{1 (p, 1)} (\R^{n} )$ is strictly bigger than  $ B^1_{ \infty,1}(\R^{n} )$ in terms of 
asymptotic  behaviors as infinity.   We also note that the use of our generalized Campanato spaces  to handle  the bounded domain problem is quite convenient as in the case of usual Campanato spaces.

\end{rem}

%%% ----------------------------------------------------------------------
%       SECTION 10
%%% ----------------------------------------------------------------------
\section{Preliminariy lemmas}
\label{sec:-10}
\setcounter{secnum}{\value{section} \setcounter{equation}{0}
\renewcommand{\theequation}{\mbox{\arabic{secnum}.\arabic{equation}}}}

Let $ X=\{ X_j\}_{ j\in \Z}$ be a sequence of non-negative real numbers. 
Given $ s \in \R$  and $ 0< q < +\infty$,  we denote 
$$ \{2^{ js}\} \cdot  X:= \{ 2^{ js} X_j\}_{j\in \Bbb Z}, \qquad  X^q:=\{ X_j^q\}_{j\in \Bbb Z}$$
respectively.
 We define $S_{\alpha, q} : X=\{ X_j\}_{j\in \Bbb Z}  \mapsto Y=\{ Y_j \}_{j\in \Bbb Z}$, where
\[
Y_j=(S_{ \alpha, q}(X))_j = 2^{ j\alpha }\Big(\sum_{i=j}^{\infty} (2^{ -i\alpha } X_i)^q\Big)^{ \frac{1}{q}},\quad  j\in \Z. 
\]
From the above definition, in case of $ \alpha =0$,  it follows that
\begin{equation}
\|S_{0, q}(X)\|_{\ell^\infty} = \|X\|_{ \ell^\infty }\leq \| X\|_{ \ell^q}\quad  \forall\,X \in \ell^q. 
\label{10.1}
\end{equation}

Clearly, for all $ \alpha , \beta  \in \R$ it holds 
\begin{equation}
2^{ \beta j }(S_{ \alpha, q}(X))_j =   S_{ \alpha + \beta ,q}(\{2^{\beta i} X_i\})_j, \quad j\in \Z.
%= S_{ \alpha + \beta ,q}(2^{\beta  }{\small \odot} X) . 
\label{10.2}
\end{equation}

Given $ X=\{ X_j\}_{j\in \Bbb Z} ,  Y=\{ Y_j \}_{j\in \Bbb Z}$,  we denote $X\leq Y$ if  $X_j\leq Y_j$ for all $j\in \Bbb Z$.
Throughout this paper, we frequently make use of the following lemma, which could be regarded as a generalization of the result in \cite{bou}.
 
\begin{lem}
\label{lem10.1}  For all $ \beta < \alpha $  and $ 0< p \le  q \le +\infty$ it holds 
\begin{equation}
 S_{ \beta , q} (S_{ \alpha, p}(X)) \le  \frac{1}{1- 2^{ -( \alpha -\beta) }}S_{\beta , q} (X). 
\label{10.3}
\end{equation}

\end{lem}

{\bf Proof}:  
We first observe
\begin{align}
\label{obs}
(S_{\beta, q}(S_{\alpha,p} X) )_j &= 2^{j\beta} \left\{ \sum_{i=j} ^\infty 2^{-i\beta q} (S_{\alpha ,p} X)_i ^q \right\}^{\frac1q}\cr
&= 2^{j\beta} \left\{ \sum_{i=j}^\infty 2^{-i\beta q} \left[ 2^{i\alpha } \left( \sum_{l=i}^\infty (2^{-\alpha l} X_l )^p \right)^{\frac1p} \right]^q \right\}^{\frac1q}\cr
&=2^{j\beta} \left\{ \sum_{i=j}^\infty 2^{i(\alpha-\beta) q}  
\left( \sum_{l=i}^\infty 2^{-(\alpha-\beta ) p l } 2^{ -\beta pl} X_l ^p\right)^{\frac{q}{p}}\right\}^{\frac1q}
\cr
&=(S_{0, q} (S_{ \alpha-\beta , p}(\{2^{ -\beta p i} X_i^p\}))_j .
\end{align}
{\it 1. The case $ p=1, \beta =0$}.  Let $ X $ be sequence with $ X_j =0$ except finite $ j\in \{ m, m+1, \ldots\}$.  
By the aid of  H\"older's inequality,  we get 
\begin{align*}
&(S_{ 0, q}(S_{ \alpha, 1}(X)))_j^q 
\\
&= \sum_{i=j}^{\infty} \Big(2^{ i \alpha }\sum_{l=i}^{\infty} 2^{ -\alpha l} X_l\Big)^q
=\sum_{i=j}^{\infty} 2^{ i q\alpha }\sum_{l=i}^{\infty} 2^{ -\alpha l} X_l
\Big(\sum_{l=i}^{\infty} 2^{ -\alpha l} X_l\Big)^{ q-1}
\\
& =\sum_{i=j}^{\infty} 2^{ i q\alpha }\sum_{l=0}^{\infty} 2^{ -\alpha (l+i)} X_{ l+i}
\Big(\sum_{l=i}^{\infty} 2^{ -\alpha l} X_l\Big)^{ q-1}
=\sum_{l=0}^{\infty} 2^{ -\alpha l}\sum_{i=j}^{\infty} X_{l+i} \Big(2^{i\alpha} \sum_{l=i}^{\infty} 2^{ -\alpha l} X_l\Big)^{ q-1}\\
& = \sum_{l=0}^{\infty} 2^{ -\alpha l}\sum_{i=j}^{\infty}   X_{ l+i}
S_{ \alpha , 1}(X)_i^{ q-1} \le \sum_{l=0}^{\infty} 2^{ -\alpha l}\left( \sum_{i=j}^{\infty} X_{l+i} ^q \right)^{\frac1q} \left( \sum_{i=j}^{\infty} (S_{\alpha,1} (X))_i ^q\right)^{\frac{q-1}{q}}
\\
&\le \frac{1}{1- 2^{ -\alpha }} (S_{0,  q}(X))_j (S_{0,  q}(S_{ \alpha, 1}(X)))_j^{ q-1},
\end{align*}
where we used the fact $ ( \sum_{i=j}^{\infty} X_{l+i} ^q)^{\frac1q} \leq  (\sum_{i=j}^{\infty} X_{i} ^q )^{\frac1q} = (S_{ 0,q} X)_j$ for all $l\ge 0$.
Dividing both sides by $(S_{0,  q}(S_{ \alpha, 1}(X)))_j^{ q-1}$, we get \eqref{10.3}.  

\vspace{0.2cm}
In the general case $ S_{0, q }(X)_j < +\infty$ we obtain   from \eqref{10.3} for the truncated sequence the property  
  $  S_{0,  q} (S_{\alpha , 1}(X))_j < +\infty$. This  shows  \eqref{10.3} for the general case.   
 
 \vspace{0.2cm}
{\it 2. The  case $ 0 < p \le q \le +\infty, \beta < \alpha $}. Recalling the definition of  $ S_{ \alpha , p}(X)$, we find 
\begin{equation}
S_{ \alpha , p}(X)_j=  \Big( 2^{ j \alpha p }\sum_{i=j}^{\infty} 2^{ -i \alpha p } X_i^p\Big)^{ \frac{1}{p}}
=  (S_{ \alpha p,1 }(\{X_i^p\}))_j^{ \frac{1}{p}},\quad j\in \Z.
\label{1.1cor}
\end{equation}
Using \eqref{1.1cor} with $ \alpha -\beta $ in place of $ \alpha $  together with   \eqref{1.1cor} with $ \beta =0$ 
and $ p=q$, we obtain the following two identities for $ j\in \Z$
\begin{align}
 (S_{ \alpha-\beta , p}(\{2^{ -\beta i} X_i\}))_j &= 
(S_{ (\alpha-\beta)p , 1}(\{2^{ -\beta p i}X_i^p\})^{ \frac{1}{p}}_j.
\label{1.1cor1}  
\\
\Big[ S_{0,  q} (
S_{ (\alpha-\beta)p , 1}(\{2^{ -\beta p i} X_i^p\}))\Big]_j ^{ \frac{1}{p}}
  &= \Big[S_{0,  1} \Big((S_{ (\alpha-\beta)p , 1}(\{2^{ -\beta p i} X_i^p\}))^{ \frac{q}{p}}
 \Big)\Big]^{ \frac{1}{q}}_j.
\label{1.1cor2}
\end{align}

Applying, $ S_{ 0,q}$ to both sides of and using first   \eqref{1.1cor1}, \eqref{1.1cor2}   together   \eqref{10.2}, and applying the inequality from the  first part of the proof, we arrive at   
\begin{align*}
\Big( S_{0, q} (S_{ \alpha-\beta , p}(\{2^{ -\beta i } X_i\})) \Big)_j&= 
\Big[ S_{0,  q} (
S_{ (\alpha-\beta)p , 1}(\{2^{ -\beta p i} X_i^p\}))\Big]_j ^{ \frac{1}{p}}
 \\
 &= \Big[S_{0,  1} \Big((S_{ (\alpha-\beta)p , 1}(\{2^{ -\beta p i} X_i^p\}))^{ \frac{q}{p}}
 \Big)\Big]^{ \frac{1}{q}}_j
  \\
 &= \Big(S_{0,  \frac{q}{p}} (S_{ (\alpha-\beta)p , 1}(\{2^{ -\beta p i} X_i^p\}))\Big)^{ \frac{1}{p}}_j
\\
& \le  \frac{1}{(1- 2^{-(\alpha-\beta)p })^{ \frac{1}{p}}} \Big(S_{0,  \frac{q}{p}} (\{2^{ -\beta p i} X_i^p\})\Big)_j^{ \frac{1}{p}}\\
&\le  \frac{1}{1- 2^{-(\alpha-\beta) }}  2^{-\beta j} (S_{\beta, q} (X))_j,
\end{align*}
where we used the fact $ (1- x^a)^{ \frac{1}{a}} \ge 1-x$ for all $ 0<x<1$ and $ a>1$. Combining this with  
\eqref{obs},  
we have \eqref{10.3}.  \hfill \Beweisende

%%% ----------------------------------------------------------------------
%       SECTION 2
%%% ----------------------------------------------------------------------
\section{Properties of the spaces  $ \mathscr{L}^{ s}_{ q(p, N)}(\R^{n})$}
\label{sec:-2}
\setcounter{secnum}{\value{section} \setcounter{equation}{0}
\renewcommand{\theequation}{\mbox{\arabic{secnum}.\arabic{equation}}}}

In this section our objective is to provide important properties of the space $  \mathscr{L}^{ k,s}_{  q(p, N)}(\R^{n})$ 
such as embedding properties, equivalent norms, interpolations properties and product estimates. 
First, let us recall the definition of the generalized mean 
for distributions $ f\in {\mathscr S}'$, where $ {\mathscr S}$ denotes the usual Schwarz class of rapidly decaying functions.  For $ f\in {\mathscr S}'$ 
and $ \varphi \in {\mathscr S}$ we dfine the convolution
\[
f \ast \varphi (x)=  \langle f, \varphi (x- \cdot )\rangle,\quad  x \in \R^{n},
\]
where $<\cdot, \cdot>$ denotes the dual pairing.
Below we use the notation $\N_0= \N\cup\{0\}$.
Then, $ f \ast \varphi \in  C^\infty(\R^{n} )$ and for every  multi index $ \alpha \in \N_0^n$ it holds 
\[
D^\alpha (f \ast \varphi) = f \ast (D^\alpha\varphi)= (D^\alpha f) \ast \varphi. 
\]

Given  $ x_0\in  \R^{n}, 0< r< +\infty$ and  $ f\in {\mathscr S}'$ we define the mean
\[
[f]^\alpha _{ x, r} = f \ast D^\alpha \varphi _r(x).
\]  
where $ \varphi _{r}(y) = r^{ -n}\varphi(r^{ -1}(y))$, and $ \varphi \in C^{\infty}_{ c}(B(1))$ stands for the standard mollifier, such that $ \intl_{ \R^{n} } \varphi dx =1$. Note that in case $ f\in L^1_{ loc}(\R^{n} )$ we get 
\[
[f]^0_{ x, r} =\intl_{ \R^{n}} f(x-y) \varphi _r(y)  dy  = 
\intl_{ B(x, r)} f(y) \varphi _{ x, r}(-y)  dy,
\]
where $ \varphi_{ x, r}= \varphi _r(\cdot +x)$. Furthermore,  from the above definition it follows that 
 \begin{equation}
[f]^{ \alpha }_{ x, r} = (D^\alpha f) \ast \varphi _{r}(x) = [D^\alpha f]^0_{ x, r}.  
 \label{5.1a}
 \end{equation}
For $ f\in L^1_{ loc}(\R^{n} )$ and  $ \alpha \in \N_0^n$ we immediately get 
\begin{equation}
 [f]^{ \alpha }_{ x, r} \le c r^{ -|\alpha |-n}  \|f\|_{ L^1(B(x, r))}\quad \forall x\in \R^{n}, r>0 . 
\label{5.1b}
 \end{equation}

\begin{lem}
\label{lem5.4}
Let $x_0\in \R^{n}, 0<r<+\infty$ and $ N\in \N_0$. For every $ f\in \mathscr{S}'$ there exists a unique polynomial 
$ P^N_{ x_0, r}(f) \in \mathcal{P}_N$ such that 
\begin{equation}
[f- P^N_{ x_0, r}(f)]_{ x_0, r}^\alpha =0\quad  \forall\,| \alpha | \le N. 
\label{5.1}
\end{equation}

\end{lem}

{\bf Proof}:  Set $ L =  \binom{n+N}{N}$.   Clearly, $ \dim \mathcal{P}_N = L$. 
We define the mapping $ T_N : \mathcal{P}_N \rightarrow  \R^{L} $, by 
\[
(T_N Q)_{ \alpha } = [Q]_{ x_0, r}^\alpha,\quad  | \alpha | \le N, \quad Q\in  \mathcal{P}_N. 
\]
In order to prove the assertion of the lemma it will be sufficient to show that $ T_N$ is injective, 
since by $ \mathcal{P}_N =L$ this implies,  $ T_N$ is also surjective. 
In fact, this can be proved by induction over $ N$. In case $ N=0$ we see this by the fact that 
\[
(T_0 1)_0 =  [1]_{ x_0, r}^0 = 1. 
\]
This $T_0  $ stands for the identity in $ \mathcal{P}_0 \cong \R$. 
 Assume $ T_{ N-1}$ is injective.  Let $ Q=  \sum_{| \alpha | \le N} a_\alpha x^\alpha \in  \mathcal{P}_N$ 
 such that  $ T_N(Q)=0$.
 Using \eqref{5.1a},  this implies for $ | \alpha |= N$
 \[
0= [Q]^{ \alpha }_{ x_0, r} = \Big[\sum_{| \beta | \le N} a_\beta D^\alpha  x^\beta \Big]^0_{ x_0, r}= [\alpha ! a_{ \alpha }]^0_{ x_0, r}= \alpha ! a_\alpha. 
\] 
Here, we used the formula $ D^\alpha  x^\beta  = \alpha ! \delta _{ \alpha \beta }$ for all $ |\beta | \le N$. 

Accordingly, $ Q \in \mathcal{P}_{ N-1}$, and it holds $ T_{ N-1} (Q)=T_{ N} (Q)=0$. By our assumption it follows $ Q=0$. This proves that $ T_N$ is injective and thus surjective. 
 \hfill \Beweisende 
 
 \begin{lem}
 \label{cor5.5}
  1. Let $ f \in \mathscr{S}'$.  Then for all $ | \beta  | \le N$ it holds
 \begin{equation}
 P^{ N- | \beta |}_{ x_0, r}(D^\beta f ) =  D^\beta  P^N_{ x_0, r}( f ).  
 \label{5.2}
 \end{equation}
 
2.  The mapping  $ P^N_{x_0, r }: L^p(B(x_0, r)) \rightarrow \mathcal{P}_N, 1 \le p \le +\infty$,  defines a projection, i.e.
\begin{align}
&\qquad  P^N_{ x_0, r}(Q) = Q\quad  \forall\,Q \in \mathcal{P}_N,
\label{5.3}
\\
&\| P^N_{ x_0, r} (f)\|_{ L^p(B(x_0, 4r))} \le c\| P^N_{ x_0, r} (f)\|_{ L^p(B(x_0, r))} \le  c\| P^N_{ 0, 1}\|_{ p}\| f\|_{ L^p(B(x_0, r))}.
\label{5.4}
\end{align}
where 
\begin{align}
 \| P^N_{ 0, 1}\|_p = \sup_{\substack{g\in L^p(B(1)) \\  g \neq 0}}\frac{\|P_{ 0,1}^N(g)\|_{ L^p(B(1))}}{\|g\|_{ L^p(B(1))}} 
 =\sup_{\substack{g\in L^p(B(x_0, r)) \\ g \neq 0}}
 \frac{\|P_{ x_0,r}^N(g)\|_{ L^p(B(x_0, r))}}{\|g\|_{ L^p(B(x_0, r))}}. 
 \label{5.4a}
\end{align}

3. For all  
 $ f\in W^{p,\, j}(B(x_0, r))$, $ 1 \le p < +\infty,  1 \le j \le N+1$,  it holds 
 \begin{equation}
 \| f- P^N_{ x_0, r}(f)\|_{ L^p(B(x_0, r))} \le c r^j\sum_{| \alpha |=j}\| D^\alpha  f- D^\alpha   P^N_{ x_0, r}(f)  \|_{ L^p(B(x_0, r))}. 
 \label{5.5}
 \end{equation}
    
 \end{lem}

 {\bf Proof}:  1. Let $ \gamma \in \N_0^n$ be a multi index with $ | \gamma | \le N- | \beta |$. 
 Obviously, $ |  \beta +\gamma  | \le N$. From the definition of $ P^N_{ x_0, r}$, observing \eqref{5.1}, and employing   \eqref{5.1a} we find 
 \begin{align*}
 [P^{ N-| \beta |}_{ x_0, r} (D^\beta f)]^{ \gamma }_{ x_0, r} &=  
   [D^\beta f]^{ \gamma }_{ x_0, r}
  \\ 
  & =D^\beta f \ast D^\gamma \varphi _r(x_0)
   =  f \ast D^{ \beta + \gamma } \varphi _r(x_0)
   \\
  & 
 =  [ f]^{\beta + \gamma }_{ x_0, r}  =  [ P^N_{ x_0, r}(f) ]^{\beta + \gamma }_{ x_0, r} = [ D^\beta  P^N_{ x_0, r}(f) ]^{ \gamma }_{ x_0, r}.
 \end{align*} 
 As we have seen in the proof of Lemma\,\ref{lem5.4}, the mapping  $ T_{ N- | \beta |}: \mathcal{P}_{ N- | \beta |} 
 \rightarrow \mathcal{P}_{ N- | \beta |}$ is injective. This yields \eqref{5.2}.   
 
 \vspace{0.3cm}
 2. We show that $ P^N_{ x_0, r}$ is a projection, i.e. $ P^N_{ x_0, r}(Q)=Q$ for all $ Q \in \mathcal{P}_N$. 
 Indeed, given $ Q \in \mathcal{P}_N$, by the definition of  $ P^N_{ x_0, r}$  \eqref{5.1} it follows that 
 \[
[Q - P^N_{ x_0, r}(Q)]^\alpha_{ x_0, r} =0\quad  \forall\,| \alpha| \le N. 
\]
Consequently, $ T_N(Q - P^N_{ x_0, r}(Q))=0$. Since $ T_N$ is injective we get $ P^N_{ x_0, r}(Q)=Q$. 
The inequality \eqref{5.4} can be verified by a standard scaling and translation argument. 
   
 \vspace{0.3cm}
 3. We prove \eqref{5.5} by induction over $ j$. For $ j=1$ \eqref{5.5} follows from the usual Poincar\'e inequality, since 
 $ [f- P^N_{ x_0, r}(f)]^0_{ x_0, r} =0$.   Assume \eqref{5.5} holds for $  j-1$.  Thus, 
 \begin{equation}
 \| f- P^N_{ x_0, r}(f)\|_{ L^p(B(x_0, r))} \le c r^{ j-1}\sum_{| \alpha |=j-1}\| D^\alpha  f- D^\alpha   P^N_{ x_0, r}(f)  \|_{ L^p(B(x_0, r))}. 
 \label{5.6}
 \end{equation}
  Thanks to \eqref{5.3}  for 
 all $ | \alpha |=j-1$ it holds,
 \[
D^\alpha   P^N_{ x_0, r}(f) = P^{ N-j+1}_{ x_0, r}(D^\alpha f). 
\]
 Hence, $ [D^\alpha f  - D^\alpha  P^{ N}_{ x_0, r}( f)]^0_{ x_0, r} =0$.   An application of the Poincar\'e inequality gives 
\begin{equation}
\| D^\alpha f  - D^\alpha  P^{ N}_{ x_0, r}( f)\|_{ L^p(B(x_0, r))} \le c r \| D D^\alpha  f- D D^\alpha   P^N_{ x_0, r}(f)  \|_{ L^p(B(x_0, r))}. 
\label{5.7}
\end{equation} 
 Combining \eqref{5.6} and \eqref{5.7}, we get \eqref{5.5}.   \hfill \Beweisende   
 
 \begin{rem}
 \label{rem1.6} From   \eqref{5.5} with $ j = N+1$ we get the generalized Poincar\'e inequality 
 \begin{equation}
 \begin{cases}
 \| f- P^N_{ x_0, r}(f)\|_{ L^p(B(x_0, r))} \le  c r^{ N+1}
 \| D^{ N+1}  f\|_{L^ p( B(x_0, r))}
  \\[0.3cm]
  \forall\,f \in W^{N+1,\, p}(B(x_0, r)).  
 \end{cases}
 \label{5.8}
 \end{equation}
 
 \end{rem}
 
 \begin{cor}
 \label{cor5.7}
 For all $ x_0 \in \R^{n}, 0< r < +\infty$, $ N\in \N_0$, and $ 1 \le p < +\infty$ it holds 
  \begin{equation}
 \| f - P^N_{ x_0, r}(f)\|_{ L^p(B(x_0, r))} \le c \inf_{ Q \in \mathcal{P}_N} \| f - Q\|_{ L^p(B(x_0, r))}
  =cr^{ \frac{n}{p}} \osc_{ p,N} (f; x_0, r). 
 \label{5.9}
 \end{equation}

 \end{cor}
 
 {\bf Proof}: Let $ Q \in \mathcal{P}_N$ be arbitrarily chosen. In view of \eqref{5.3} we find 
 \[
f - P^N_{ x_0, r}(f) = f-Q - P^N_{ x_0, r}(f-Q).
\] 
Hence, applying triangle inequality, along with \eqref{5.4} we get  
 \begin{align*}
 \| f - P^N_{ x_0, r}(f)\|_{ L^p(B(x_0, r))} &\le  \| f - Q\|_{ L^p(B(x_0, r))}+  \| P^N_{ x_0, r}(f-Q)\|_{ L^p(B(x_0, r))}
 \\
 & \le c \| f - Q\|_{ L^p(B(x_0, r))}. 
 \end{align*}
 This shows the validity of \eqref{5.9}.  \hfill \Beweisende  
 
\vspace{0.3cm}
In our discussion below and in the sequel of the paper it will be convenient to work with smooth functions.  
Using the standard mollifier we get the following estimate in  $ \mathscr{L}^{ k,s}_{q (p, N) }(\R^{n})$ for the mollification.

\begin{lem}
\label{lem5.3}
Let $ \varepsilon >0$. 
Given  $ f \in \mathscr{S}'$, we define the mollification  
\[
f_\var (x) =  [f]_{ x, \varepsilon }^{ 0}= f \ast \varphi _{ \varepsilon }(x), \quad  x\in \R^{n}.   
\]
1. For all $ f\in \mathscr{L}^{ k,s}_{q (p, N) }(\R^{n})$,  and all  $ \var >0$ it holds 
\begin{equation}
| f_\var |_{ \mathscr{L}^{ k,s}_{q (p, N) } } \le c| f|_{\mathscr{L}^{ k,s}_{q (p, N) } }. 
\label{5.10}
\end{equation} 
 
 2. Let $ f\in L^p_{ loc}(\R^{n})$ such that for all $ 0<\var<1 $,
\begin{equation}
| f_\var |_{ \mathscr{L}^{ k,s}_{q (p, N) } } \le c_0,
\label{5.11}
\end{equation}  
then $ f\in \mathscr{L}^{ k,s}_{q (p, N) }(\R^{n})$ and it holds $ | f |_{\mathscr{L}^{ k,s}_{q (p, N) }} \le c_0.$
\end{lem}

{\bf  Proof}: 1. We may restrict ourself to the case $ k=0$. Let $ x_0\in \R^{n}$ and $ j\in \Z$. Set $ 0< r<+\infty$.  By the definition
of $ P^N_{ x_0, r}(f)$ (cf. \eqref{5.1} ) together with \eqref{5.1a}  it follows that for all $ | \alpha | \le N$ and for almost all  $ y\in \R^{n}$, 
\begin{align*}
  f \ast D^\alpha \varphi _{ r}(x_0-y) 
= [f]^{ \alpha }_{ x_0-y, r}=
[P^N_{ x_0-y, r}(f)]^{ \alpha }_{ x_0-y, r} 
 =P^N_{ x_0-y, r}(f) \ast D^\alpha \varphi _{r}(x_0-y).
\end{align*}
Multiplying both sides by $ \varphi _{0, \var  } (y)$, integrate the result over $ \R^{n}$ and apply Fubini's theorem, we get for all $ | \alpha | \le N$
\begin{align*}
[f_\varepsilon ]^{ \alpha }_{ x_0, r}& = (f \ast \varphi _\varepsilon   \ast D^\alpha \varphi _{r}) (x_0) = 
(f \ast D^\alpha \varphi _{r} \ast \varphi _{ \var })(x_0)
\\
&=\intl_{ \R^{n}}  (f \ast D^\alpha \varphi _{  r})(x_0- y) \varphi _{ \var } (y) dy
\\
& =  \intl_{ \R^{n}}  P^N_{ x_0-y, r}(f) \ast D^\alpha \varphi _{r}(x_0-y) \varphi _{ \var } (y)  dy
\\
& =  \intl_{ \R^{n}} \intl_{ \R^{n}} P^N_{ x_0-y, r}(f)(x)  D^\alpha \varphi _{r}(x_0-y- x) \varphi _{ \var } (y) dx  dy
\\
& =  \intl_{ \R^{n}} \intl_{ \R^{n}} P^N_{ x_0-y, r}(f)(x-y)  D^\alpha \varphi _{r}(x_0- x) \varphi _{ \var } (y) dx  dy
\\
& =  \intl_{ \R^{n}} \intl_{ \R^{n}} P^N_{ x_0-y, r}(f)(x-y)  \varphi _{ \var } (y) dy  
D^\alpha \varphi _{r}(x_0-x) dx
\\
& =  \bigg[\intl_{ \R^{n}} P^N_{ x_0-y, r}(f)(x- y)  \varphi _\varepsilon   (y) dy\bigg]^{ \alpha }_{ x_0, r}. 
\end{align*}
This shows that 
\begin{align}
P_{ x_0, r}^N(f_\var )(x) &= \intl_{ \R^{n}} P^N_{ x_0-y, r}(f)(x- y)  \varphi_\varepsilon    (y) dy,\quad  x\in \R^{n},
\\
&=  \intl_{ \R^{n}} P^N_{ x_0-\varepsilon y, r}(f)(x- \varepsilon y)  \varphi   (y) dy,\quad  x\in \R^{n}. 
\label{5.12}
\end{align}

Accordingly, 
\begin{align*}
| f_\var (x) - P_{ x_0, 2^j}^N(f_\var )(x)|^p & \le  \bigg( \int_{ \R^{n}} | f(x- \var y) - P^N_{ x_0-\var y, 2^j}(f)(x- \var y)|  \varphi  (y) dy \bigg)^{p}. 
\end{align*}
Integration of both sides over $ B(x_0, 2^j)$ and multiplication with $ \frac{1}{| B(2^j)|}$,  using Jensen's inequality with respect 
to the probability measure $ \varphi  dy$, we find 
\begin{align*}
\osc_{ p,N} (f_\var ; x_0, 2^j)  
 & \le \Bigg(\intmw_{B(x_0, 2^j)} \bigg( \int_{ \R^{n}} | f(x- \var y) - P_{ x_0-\var y, 2^j}(f)(x- \var y)| \varphi  (y) dy \bigg)^{ p}dx\Bigg)^{  \frac{1}{p} }
\\
&= \int_{ \R^{n}} \bigg(\intmw_{B(x_0- \var y, 2^j)}| f(x) - P_{ x_0-\var y, 2^j}(f) (x)|^p  dx \bigg)^{  \frac{1}{p} } \varphi  (y) dy
\\
& \le c \intl_{ B(1)} \osc_{ p,N}(f; x_0- \var y; 2^j) \varphi (y)dy.
\end{align*}
Multiplying both sides by $ 2^{ -js}$ applying the $ \ell^q$ norm to both sides of the resultant  inequality, and using Minkowski's inequality, we are led to 
\begin{align*}
\Big(\sum_{j\in \Z}( 2^{ -js}\osc_{ p,N} (f_\var ; x_0, 2^j))^q   \Big)^{ \frac{1}{q}}
 &\le  c\intl_{ B(1)} \Big(\sum_{j\in \Z}( 2^{ -js} \osc_{ p,N}(f; x_0- \var y; 2^j))^q \Big)^{ \frac{1}{q}}\varphi (y)dy
\\
& \le c | f|_{\mathscr{L}^{ s}_{q (p, N)} }.
\end{align*}
Taking the supremum over all $ x_0 \in \R^{n}$ in the above inequality shows \eqref{5.10}.  

\vspace{0.3cm}
2. Let $ f\in L^p_{ loc}(\R^{n})$ satifying \eqref{5.11}.  This implies that $ f\in W^{k,\, p}_{ loc}(\R^{n})$. 
Let $ x_0 \in \R^{n}$ and $ l,m\in \Z$, $ l < m$. According to the absolutely 
continuity of the Lebesgue measure together with \eqref{5.11} it follows 
\[
\sum_{j=l}^{m}( 2^{ -js}\osc_{ p,N} (D^k f ; x_0, 2^j))^q   = \lim_{\var  \searrow 0}
 \sum_{j=l}^{m}( 2^{ -js}\osc_{ p,N} (D^kf_\var  ; x_0, 2^j))^q \le c_0 ^q. 
\]
 This shows that $  \{2^{ -sj}{\D \osc_{ p,N} (D^k f ; x_0, 2^j)}\}_{ j\in \Z}  \in \ell^q$,  and  its sum is bounded by $ c_0$. 
 Accordingly, $ f\in \mathscr{L}^{ k,s}_{q (p, N)}(\R^{n})$, and it holds $| f|_{ \mathscr{L}^{ k,s}_{q (p, N)}} \le c_0$. 
 \hfill \Beweisende 
\vspace{0.3cm}

We are are now in a position to prove  the following embedding properties.  
First, let us introduce  the definition of the  projection 
to the space of  homogenous polynomial $\dot{P}^N_{ x_0,r}:  \mathscr{S}' \rightarrow   \dot{\mathcal{P}} _N$ 
defined by means of  
\[
\dot{P}_{ x_0, r}^N(f) (x) = \sum_{| \alpha |=N}  \frac{1}{ \alpha !} [ f]^{ \alpha }_{ x_0, r} x^\alpha,\quad  x\in \R^{n}.
\]
Clearly, for all  $ f\in \mathscr{S}'$ it holds 
\begin{equation}
D^\alpha \dot{P}_{ x_0, r}^N(f) =\dot{P}_{ x_0, r}^{ N-| \alpha |}(D^{ \alpha } f)\quad  \forall\,| \alpha | \le k.
\label{5.12a}
\end{equation}

\begin{thm}
\label{thm5.4} 
1. For every $ N\in \N_0$ the following embedding holds true
\begin{equation}
\begin{cases}
 \mathscr{L}^{N}_{1 (p,  N)}( \R^{n}) \hookrightarrow  C^{N-1, 1}(\R^{n})\quad  &\text{ if}\quad N \ge 1
\\[0.3cm]
 \mathscr{L}^{0}_{1 (p, 0)}( \R^{n}) \hookrightarrow  L^\infty(\R^{n})\quad &\text{ if}\quad N = 0.
\end{cases}
\label{5.13}
\end{equation}

2. For every $ f\in  \mathscr{L}^{N}_{1 (p,  N)}( \R^{n})$ there exists a unique $ \dot{P} ^N_{ \infty}\in \dot{\mathcal{P}} _N$, such that 
for all $ x_0 \in \R^n$
\[
\lim_{r \to \infty}  \dot{P} _{ x_0, r}(f) \rightarrow \dot{P}^N_{ \infty}(f)\quad  \text{ in }\quad \mathcal{P}_N. 
\]
Furthermore, $\dot{P}^N_{ \infty}:  \mathscr{L}^{N}_{1 (p,  N)}( \R^{n}) \rightarrow \dot{\mathcal{P}} _N$ is a projection, with the property 
\begin{equation}
D^\alpha \dot{P}^N_{ \infty}(f) = \dot{P}^{ N-| \alpha |}_\infty (D^\alpha f) \quad  \forall\,| \alpha | \le N. 
\label{5.14}
\end{equation}
3. For all $ g, f\in \mathscr{L}^{1}_{1 (p,  1)}( \R^{n})$ it holds
 \begin{align}
 \dot{P}^1_{ \infty}(g \partial _k f) &= 
 \dot{P}^1_{ \infty}(g) \partial _k\dot{P}^1_{ \infty} (f) = \dot{P}^1_{ \infty}(g) \dot{P}^0_{ \infty} (\partial _k f), \quad 
 k=1, \ldots, n.
 \label{5.15}
\end{align}
In addition, for $ g\in C^{ 0,1}(\R^{n}; \R^{n}  )$, and for all 
$ f\in \mathscr{L}^{0}_{1 (p,  0)}( \R^{n})$  it holds 
\begin{align}
  \dot{P}^0_{ \infty}(g\partial _k  f) := \lim_{r \to \infty}  P^0_{ 0, r}(g\partial _k   f) =0,\quad k=1, \ldots, n, 
\label{5.15b}
\end{align}
where $ g\partial _k   f = \partial _k(gf)- \partial _kg f\in \mathscr{S}'$. 

\vspace{0.3cm}
4. For all $ v\in \mathscr{L}^{1}_{1 (p,  1)}( \R^{n}; \R^{n} )$ with $ \nabla \cdot  v=0$ almost everywhere in $ \R^{n}$  
and $ f\in \mathscr{L}^{1}_{1 (p,  1)}( \R^{n})$ it holds 
\begin{align}
 \dot{P}^0_{ \infty}(\nabla v \cdot \nabla f) &= 
 \dot{P}^0_{ \infty}(\nabla v)\cdot \dot{P}^0_{ \infty}(\nabla f). 
\label{5.15a}
\end{align}

\end{thm} 

{\bf Proof}: 1. Let $ \var >0$ be arbitrarily chosen. Let $ f\in \mathscr{L}^{N}_{1 (p, N) }(\R^{n})$.  Set $ f_\var = 
f \ast \varphi _{ \var  } $. 
By Lemma\,\ref{lem5.3} we get $ f_\var \in \mathscr{L}^{N}_{1 (p, N) }(\R^{n})$ and it holds 
\begin{equation}
| f_\var |_{ \mathscr{L}^{  N}_{1 (p, N) }} \le c| f|_{ \mathscr{L}^{  N}_{1 (p, N) }}. 
\label{5.17}
\end{equation}
Let $ x_0\in \R^{n}$ be fixed. Let $ j\in \Z$.   Clearly, $ f_\var \in C^\infty(\R^{n})$. Let $ \alpha \in \N_0^n$ be a multi index with 
$ | \alpha | =N$.  
Then 
\[
D^\alpha P^N_{ x_0, 2^j}(f_\var ) = P^0_{ x_0, 2^j}(D^{ \alpha }(f_\var) )=  [D^{ \alpha }(f_\var)]^0_{ x_0, 2^j}  = 
D^\alpha \dot{P} ^N_{ x_0, 2^j}(f_\var ).
\]
Let $ m\in \Z$. Since $D^\alpha f_\var $ is continuous we have 
\[
D^\alpha f_\var  (x) = \lim_{ j\to -\infty} [ D^\alpha f_\var ]^{ 0}_{ x, 2^{ j}}\quad \forall x\in \R^{n}.
\]

Using triangle inequality along with \eqref{5.3} and \eqref{5.10},  and using  \eqref{5.1b},  we get
\begin{align}
&| D^\alpha f_\var  (x) -  [D^\alpha f_\var ]^{ 0}_{ x, 2^m}|
= \Big |\sum_{j=-\infty}^{m}  [ D^\alpha f_\var ]^{ 0}_{ x, 2^{ j-1}}-  [ D^\alpha f_\var ]^{ 0}_{ x, 2^j}\Big |
\cr
&\quad \le \sum_{j=-\infty}^{m} \Big| [ D^\alpha f_\var ]^{ 0}_{ x, 2^{ j-1}}-  [ D^\alpha f_\var ]^{ 0}_{ x, 2^j}\Big|
 =\sum_{j=-\infty}^{m} \Big| [  f_\var ]^{ \alpha }_{ x, 2^{ j-1}}-  [f_\var ]^{ \alpha }_{ x, 2^j}\Big|
\cr
&\quad =\sum_{j=-\infty}^{m} \Big| [ f_\var -  P^N_{ x, 2^j} ( f_\var )]^{ \alpha }_{ x, 2^{ j-1}}-  [ f_\var - P^N_{ x, 2^j} (f_\var )]^{ \alpha }_{ x, 2^j}\Big|
\cr
&\quad =\sum_{j=-\infty}^{m} \Big| [ f_\var -  P^N_{ x, 2^j} ( f_\var )]^{ \alpha }_{ x, 2^{ j-1}}\Big|
\cr
&\quad \le c \sum_{j=-\infty}^{m}2^{ -jN }\osc_{ p, N}( f_\var ; x, 2^j ) 
\le c | f_\var |_{ \mathscr{L}^{ N}_{ 1 (p, N)}} \le
 c | f |_{ \mathscr{L}^{ N}_{ 1 (p, N)}}.  
 \label{5.18}
\end{align}
Thus, $\{  D^N f_\var \}$ is bounded in $ L^\infty(B(r) )$ for all $ 0<r<+\infty $.  By means of Banach-Alaoglu's theorem 
and Cantor's diagonalization principle we get a sequence $ \var _k \searrow 0$ as $ k \rightarrow +\infty$ and $ f\in W^{N,\, \infty} _{ loc}(\R^{n})$, such that for all $ 0< r< +\infty$
\[
 D^N f_{ \var _k} \rightarrow D^N f  \quad  \text{{\it weakly$-\ast$ in}}\quad L^\infty(B(r)) \quad  \text{{\it as}}\quad  k \rightarrow +\infty.
\]
Furthermore, from \eqref{5.18}  we get for almost all $ x\in \R^{n}$ and all $ m\in \Z$, 
\begin{equation}
| D^N f (x)|\le c | f |_{ \mathscr{L}^{ N}_{ 1 (p, N)}} + \sum_{| \alpha |=N}| [f ]^{ \alpha }_{ x, 2^m}|.   
\label{5.19}
\end{equation}  
Let $ x_0 \in \R^{n}$ be fixed. We now choose $ m\in \Z $ such that $ 2^{ m-1} \le | x_0| < 2^m$. Then noting 
$ B(x_0, 2^m) \subset B(2^{ m+1})$, employing  \eqref{5.3} and  \eqref{5.1b}, we get   
\begin{align*}
\Big|| [f ]^{ \alpha }_{ x_0, 2^m}| - | [f ]^{ \alpha }_{0, 2^m}|\Big|& \le 
 [f ]^{ \alpha }_{ x_0, 2^m} -  [f ]^{ \alpha }_{0, 2^m}
  = \Big| [f- P^N_{ 0, 2^{ m+1}} ]^{ \alpha }_{ x_0, 2^m} -  [f- P^N_{ 0, 2^{ m+1}} ]^{ \alpha }_{0, 2^m}\Big|
\\
&\le c 2^{ - mN} \osc_{ p, N}(f; 0, 2^{ m+1}) \le c | f |_{ \mathscr{L}^{ N}_{ 1 (p, N)}}.
\end{align*}
 Similarly, we get for all $ j\in \Z$
 \[
\Big|| [f ]^{ \alpha }_{ 0, 2^m}| - | [f ]^{ \alpha }_{0, 2^j}|\Big| \le c | f |_{ \mathscr{L}^{ N}_{ 1 (p, N)}}. 
\]
Thus, combining the two inequalities we have just obtained, using triangle inequality, we find for all  $ j\in \Z$,
\[
\sum_{| \alpha |=N}|[f ]^{ \alpha }_{ x, 2^m}| \le c | f |_{ \mathscr{L}^{ N}_{ 1 (p, N)}}+
 \sum_{| \alpha |=N}|[f ]^{ \alpha }_{0, 2^j}| . 
\] 
This together 
with \eqref{5.19}, we infer for all $ j\in \Z$
\begin{align}
\| D^N f\|_{ \infty}&\le c | f |_{ \mathscr{L}^{ N}_{ 1 (p, N)}} + c\sum_{| \alpha |=N} | [f ]^{ \alpha }_{ 0, 2^j}|    
\le c | f |_{ \mathscr{L}^{ N}_{ 1 (p, N)}} + c\| \dot{P}^N _{ 0, 2^j}(f)\|.  
\label{5.20}
\end{align}
 This completes the proof of \eqref{5.13}.

\vspace{0.5cm}  
2. Let $ x_0 \in \R^{n}$.  Let $ m, l\in \Z$, $ l< m$. Noting that  
$ \dot{P}^N_{ x_0, 2^j}(Q)=Q$ 
for all $ Q\in \dot{\mathcal{P}}_N$ and  $ \dot{P}^N_{ x_0, 2^j}(Q)=0$ for all $ Q\in \mathcal{P}_{ N-1}$,  
 we get the following identity for all $ j, k\in \Z$
\[
\dot{P}^N_{ x_0, 2^{ j}}(P^N_{ x_0, 2^k}(f)) =\dot{P}^N_{ x_0, 2^{ k}}(f).
\]
 Using triangle inequality together with  the above identity, \eqref{5.1b} and   \eqref{5.9} we estimate  
\begin{align*}
&\| \dot{P}^N_{ x_0, 2^l}(f)- \dot{P}^N_{ x_0, 2^m}(f) \| 
\\
&\le \sum_{j=l+1}^{m} \| \dot{P}^N_{ x_0, 2^{ j-1}}(f)- \dot{P}^N_{ x_0, 2^j}(f) \| 
\\
&\le c \sum_{j=l+1}^{m}2^{ -jN- j\frac{n}{p}} \| \dot{P}^N_{ x_0, 2^{ j-1}}(f- P^N_{ x_0, 2^j}(f))- \dot{P}^N_{ x_0, 2^j}(f- P^N_{ x_0, 2^j}(f)) \|_{ L^p(B(x_0, 2^j))}
\\
&\le c\sum_{j=l+1}^{m} 2^{ -jN} \osc_{ p, N}(f; x_0, 2^j). 
\end{align*}
Owing to  $ f\in \mathscr{L}^{ N}_{ 1 (p, N)}(\R^{n})$ the right-hand side of the above inequality tends to zero as $ m, l \rightarrow +\infty$. This shows that $ \{ \dot{P}^N_{ x_0, 2^m}(f)\}$ is a Cauchy sequence in $ \dot{P} _N$ and converges to a unique limit 
$ \dot{P} ^N_{ \infty, x_0}$. We claim that 
\begin{equation}
 \dot{P} ^N_{ \infty, x_0} =  \dot{P} ^N_{ \infty, 0}=:  \dot{P} ^N_{ \infty}(f). 
\label{5.21}
\end{equation}
In fact, for $ m \in \Z$  such that $ | x_0| \le 2^m $, we obtain
\begin{align*}
&\| \dot{P}^N_{ x_0, 2^m}(f)- \dot{P}^N_{ 0, 2^m}(f) \| 
\\
&\le c 2^{ -mN- m\frac{n}{p}} \| \dot{P}^N_{ x_0, 2^{ m}}
(f- P^N_{ 0, 2^{ m+1}}(f))- \dot{P}^N_{ 0, 2^m}(f- P^N_{ 0, 2^{ m+1}}(f)) \|_{ L^p(B(x_0, 2^m))}
\\
&\le c 2^{ -mN} \osc_{ p, N}(f; 0, 2^{ m+1}) \rightarrow 0\quad  \text{ as}\quad m \rightarrow +\infty. 
\end{align*}
Consequently, \eqref{5.21} must hold.   The identity \eqref{5.14} is an immediate consequence of \eqref{5.12a}.

\vspace{0.2cm}
3. Now, let $ g,f\in \mathscr{L}^{ 1}_{ 1 (p, 1)} (\R^{n})$. Let 
$ x_0\in \R^{n}$. Let $ \alpha \in \N_0^n$ with $ |\alpha |=1$. We first show that $\{ [g\partial _k f ]^{ \alpha }_{ x_0, 2^j}\}_{ j\in \N}, k\in \{1, \ldots,n\}$, is a Cauchy sequence.  Let $ j\in \N$ be fixed. We easily calculate, 
\begin{align*}
& [g\partial _k f]^{ \alpha }_{ x_0, 2^{ j-1}}- [g\partial _k f ]^{ \alpha }_{ x_0, 2^j} 
\\
&= \Big[g\partial _k f-P^1_{ x_0, 2^j}(g)  [\partial _k f]^0 _{ x_0, 2^j} \Big]^{ \alpha }_{ x_0, 2^{ j-1}}
 - 
\Big[g\partial _k f-P^1_{ x_0, 2^j}(g)  [\partial _k f]^0 _{ x_0, 2^j} \Big]^{\alpha }_{ x_0, 2^{ j}}. 
\end{align*}
Furthermore, applying integration by parts, we get,
\begin{align*}
&\Big[g\partial _k f-P^1_{ x_0, 2^j}(g)  [\partial _k f]^0 _{ x_0, 2^j} \Big]^{ \alpha }_{ x_0, 2^{ j-1}}
\\
&= \Big[(g- P^1_{ x_0, 2^j}(g)  ([\partial _k f]^0 _{ x_0, 2^j}) \Big]^{ \alpha }_{ x_0, 2^{ j-1}}
+ \Big[ g\cdot  (\partial _k f - [\partial _k f]^0 _{ x_0, 2^j}) \Big]^{ \alpha }_{ x_0, 2^{ j-1}} 
\\
&=- \intl_{ \R^{n}} 
  (g- P^1_{ x_0, 2^j}(g))   ([\partial _k f]^0 _{ x_0, 2^j}) D^\alpha \varphi _{x_0, 2^{ j-1}} dx  
  \\
&  \qquad \qquad -  \intl_{ \R^{n}} 
 g  (\partial _k f - [\partial_k  f]^0 _{ x_0, 2^j}) D^\alpha \varphi _{ x_0, 2^{ j-1}}dx  
\\
&= - \intl_{ \R^{n}} 
  (g- P^1_{ x_0, 2^j}(g))   ([\partial _k f]^0 _{ x_0, 2^j}) D^\alpha \varphi _{x_0, 2^{ j-1}} dx  
  \\
&  \qquad \qquad +  \intl_{ \R^{n}} 
  \partial _k g   ( f-P^1_{ x_0, 2^j}(f) - [f-P^1_{ x_0, 2^j}(f)]^1_{ x_0, 2^j})  D^\alpha  \varphi _{ x_0, 2^{ j-1}} dx
  \\
&  \qquad \qquad +  \intl_{ \R^{n}} 
   g   ( f-P^1_{ x_0, 2^j}(f) - [f-P^1_{ x_0, 2^j}(f)]^1_{ x_0, 2^j})  \partial _kD^\alpha  \varphi _{ x_0, 2^{ j-1}} dx.     
  \end{align*}
This together with  \eqref{5.9} yields 
\begin{align*}
&\Big[g\partial _k f-P^1_{ x_0, 2^j}(g)  [\partial _k f]^0 _{ x_0, 2^j} \Big]^{ \alpha }_{ x_0, 2^{ j-1}}
\\
&\qquad \le c \| \nabla f\|_{ \infty} 2^{ -j} \osc_{ p, 1} (v; x_0, 2^j) +  c \| \nabla v\|_{ \infty} 2^{ -j} \osc_{ p, 1} (f; x_0, 2^j). 
\end{align*}
By an analogous reasoning we find 
\begin{align*}
&\Big[g\partial _k f-P^1_{ x_0, 2^j}(g)  [\partial _k f]^0 _{ x_0, 2^j} \Big]^{ \alpha }_{ x_0, 2^{ j}}
\\
&\qquad \le c \| \nabla f\|_{ \infty} 2^{ -j} \osc_{ p, 1} (v; x_0, 2^j) +  c \| \nabla v\|_{ \infty} 2^{ -j} \osc_{ p, 1} (f; x_0, 2^j). 
\end{align*}
Let $l,  m \in \Z$ with $ l < m$ be arbitrarily chosen. Using triangle inequality together with the two estimates we have just obtained, we estimate 
\begin{align*}
&\Big|[g\partial _k f ]^{ \alpha }_{ x_0, 2^{ l}}- [g\partial _k f]^{ \alpha }_{ x_0, 2^m} \Big|
\\
& = \Big|\sum_{j=l+1}^{m}[g\partial _k f]^{ \alpha }_{ x_0, 2^{ j-1}} - [g\partial _k f]^{ \alpha }_{ x_0, 2^j} \Big|
\\
&\qquad \le c \| \nabla f\|_{ \infty} \sum_{j=l+1}^{m}2^{ -j} \osc_{ p, 1} (g; x_0, 2^j) +  c \| \nabla g\|_{ \infty} 
\sum_{j=l+1}^{m}2^{ -j} \osc_{ p, 1} (f; x_0, 2^j). 
\end{align*}
Since $ g, f \in \mathscr{L}^{ 1}_{ 1 (p, 1)} (\R^{n})$ the right-hand side converges to zero as $  l,m \rightarrow +\infty$. Thus, 
$\{  [g\partial _k f ]^{ \alpha }_{ x_0, 2^{ l}}\}$ is a Cauchy sequence,  and has a unique limit say $ a_{ x_0}$.  
Let $ j\in \N$ such that $ 2^j \ge | x_0|$. Thus, $ B(x_0, 2^j) \subset B(2^{ j+1})$. By the same reasoning as above we estimate 
\begin{align*}
&\Big|[g\partial _k f]^{ 0}_{ x_0, 2^{ j}}- [g\partial _k f ]^{ \alpha }_{ 0, 2^{ j+1}} \Big|
\\
& = c \| \nabla f\|_{ \infty} 2^{ -j} \osc_{ p, 1} (g; 0, 2^{ j+1}) +  c \| \nabla g\|_{ \infty} 2^{ -j} \osc_{ p, 1} (f; 0, 2^{ j+1}). 
\end{align*}
Since the right-hand side converges to zero as $ j \rightarrow +\infty$ we get $ a_{ x_0}= a_0$. Setting 
$ [g\partial _k f]^{ \alpha }_\infty= a_0$, we complete the proof of  \eqref{5.15}.

\vspace{0.3cm}
Next, we prove  \eqref{5.15b}. Let $ g\in C^{ 0,1}(\R^{n} )$ and $ f\in \mathscr{L}^{ 0}_{ 1 (p, 0)} (\R^{n})$. 
Applying integration by parts and product rule, we calculate  
\begin{align*}
[g\partial _k f]^0_{ x_0, r} 
&=-\intl_{ B(x_0,r)} \partial _k g(y) (f(y)-[f]^0_{ x_0, r}) \varphi_r(x_0-y) dy 
\\
&\quad +\intl_{ B(x_0,r)} g(y) (f(y)-[f]^0_{ x_0, r}) \partial _k\varphi_r(x_0-y) dy. 
\end{align*}
Applying H\"older's inequality, we easily get 
\begin{align*}
|[g\partial _k f]^0_{ x_0, r}| \le   c \|\nabla g\|_{ \infty} \osc_{ p,0} (f; x_0, r)
+
c r^{ -1}\| g\|_{ L ^\infty(B(x_0, r))} \osc_{ p,0} (f; x_0, r).
\end{align*}
Noting that $r^{ -1}\| g\|_{ L ^\infty(B(x_0, r))}  \le c |g(x_0)| + c \|\nabla g\|_{ \infty}  $, and using the fact that
$ \osc_{ p,0} (f; x_0, r) \rightarrow 0$ as $ r \rightarrow r+\infty$, we obtain  \eqref{5.15b}.

\vspace{0.2cm}
It remains to show the identity \eqref{5.15a}.  Let  $ v\in \mathscr{L}^{1}_{1 (p,  1)}( \R^{n}; \R^{n} )$ with $ \nabla \cdot  v=0$ and $ f\in \mathscr{L}^{1}_{1 (p,  1)}( \R^{n})$. Using  \eqref{5.14} together with $ \nabla \cdot  v=0$ and  \eqref{5.15},  we obtain
\begin{align*}
[\partial _k v \cdot \nabla f]^0_{ x_0, r} &= \partial _j P^1_{ x_0, r}((\partial _k v_j) f)= 
\partial _j \dot{P}^1_{ x_0, r}((\partial _k v_j) f)
\\
 &\rightarrow   \partial _j \dot{P}^1_{ \infty}(\partial _k v_j) \dot{P}^1_{ \infty}(f) = \dot{P}^0_\infty(\partial _kv) \cdot
 \dot{P}^0_\infty( \nabla f)\quad 
  \text{as}\quad  r\rightarrow +\infty.
\end{align*}
This shows that 
\[
 \dot{P}^0_\infty(\partial _k v \cdot \nabla f) = \lim_{r \to \infty}[\partial _k v \cdot \nabla f]^0_{ x_0, r}=
 \dot{P}^0_\infty(\partial _kv) \cdot
 \dot{P}^0_\infty( \nabla f).
\]

This completes the proof of the Lemma. 
  \hfill \Beweisende

\vspace{0.3cm}
 Next, we prove the following norm equivalence which is similar to the properties of the known Campanato space.  
 
 \begin{lem}
 \label{lem5.9}  Let  $1 \le p < +\infty, 1 \le q \le +\infty$, and $ N, N'\in \N_0, N< N',  s \in [- \frac{n}{p}, N+1)$. 
  If $ f\in \mathscr{L}^{k,s}_{q  (p, N')}(\R^{n})$, and satisfies   
 \begin{equation}
   \lim_{m \to \infty}\dot{P}^{ L}_{ 0, 2^m}(D^k f) = 0\quad \forall\,L=N+1, \ldots, N'.
 \label{5.23}
 \end{equation} 
then $ f\in \mathscr{L}^{k,s}_{q (p, N)}(\R^{n})$ and it holds, 
 \begin{equation}
  | f|_{  \mathscr{L}^{ k,s}_{q(p, N')}} \le | f|_{  \mathscr{L}^{ k,s}_{ q(p,N)}} \le c  | f|_{  \mathscr{L}^{ k,s}_{q(p, N')}}.
 \label{5.24}
 \end{equation}

 \end{lem} 
 
{\bf  Proof}:  We may restrict ourself to the case $ k=0$. First, let us prove that for all $  s \in [- \frac{n}{p}, N)$ 
and for all $ f\in \mathscr{L}^{ s}_{q (p, N)}(\R^{n}) $  such that 
 \begin{equation}
  \lim_{m \to \infty}\dot{P}^{ N}_{ 0, 2^m}(f) = 0. 
 \label{5.24a}
  \end{equation} 
 it follows that  $ f\in \mathscr{L}^{ s}_{q (p, N-1)}(\R^{n}) $, together with the estimate
\begin{equation}
 |f|_{\mathscr{L}^{ s }_{q(p, N-1)}} \le c  |f|_{\mathscr{L}^{ s }_{q(p, N)}}. 
\label{5.25}
\end{equation}
Let $ x_0\in \R^{n}, 0< r< +\infty$. Noting that $  P^N_{ x_0, 2r}(f) - \dot{P}^N_{ x_0, 2r}(f) \in \mathcal{P}_{ N-1}$, 
we see that 
\[
\dot{P}^N_{ x_0, r}( P^N_{ x_0, 2r}(f))  = \dot{P}^N_{ x_0, r}(P^N_{ x_0, 2r}(f) - \dot{P}^N_{ x_0, 2r}(f))+ 
\dot{P}^N_{ x_0, r}( \dot{P}^N_{ x_0, 2r}(f))=
\dot{P}^N_{ x_0, 2r }(f). 
\]

By a scaling argument and triangle inequality  we infer 
\begin{align*}
& r^{ -N- \frac{n}{p}}\| \dot{P}_{x_0, r}^N(f)\|_{ L^p(B(x_0, r))}- (2r)^{ -N- \frac{n}{p}}\| \dot{P}_{x_0,  2r}^N(f)\|_{ L^p(B(x_0,  2r))}
\\
&= \| \dot{P}_{x_0, r}^N(f)\| - \| \dot{P}_{x_0, 2r}^N(f)\|
 \le \|  \dot{P}_{x_0 r}^N(f)- \dot{P}_{x_0,  2r}^N(f) \|
\\
& = \| \dot{P}_{x_0,  r}^N( f- P^N_{ x_0, 2r}(f) \|
\\
& \le c (2r)^{ -N- \frac{n}{p}}\| f-  P^N_{ x_0, 2r}(f) \|_{ L^p(B(x_0,  2r))}
\\
& \le c r^{ -N} \osc_{ p,N}(f; x_0, 2r).
\end{align*}
Let $j , m \in \Z$, $ j < m$. Using the above estimate we deduce that  
\begin{align*}
& \Big|2^{ -jN- j\frac{n}{p}}\| \dot{P}_{x_0,  2^j }^N(f)\|_{ L^p(B(x_0,  2^j ))} -
2^{ -mN- m\frac{n}{p}}\| \dot{P}_{x_0, 2^N}^N(f)\|_{ L^p(B(x_0, 2^m))}\Big|
\\
&\le c\sum_{i=j}^{m-1} 2^{ - iN}\osc_{ p,N}(f; x_0, 2^{ i+1})
\\
&\le c 2^N\sum_{i=j}^{m-1} 2^{ - iN}\osc_{ p,N}(f; x_0, 2^{ i}).
\end{align*}
Observing \eqref{5.24a}, we see that  
\[
 \lim_{m \to \infty}\|2^{ -mN- m\frac{n}{p}}\| \dot{P}_{x_0, 2^N}^N(f)\|_{ L^p(B(x_0, 2^m))} =\lim_{m \to \infty}\| \dot{P}_{x_0, 2^m}^N(f)\| =0. 
\]
Thus, letting $ m \rightarrow +\infty$ in the above estimate, we arrive at 
\begin{align}
2^{ jN}\| \dot{P}_{x_0,  2^j }^N(f)\| &= 2^{- j\frac{n}{p}}\| \dot{P}_{x_0,  2^j }^N(f)\|_{ L^p(B(x_0,  2^j ))} 
\le c2^{ jN}\sum_{i=j}^{\infty} 2^{ - iN}\osc_{ p,N}(f; x_0, 2^{ i})
\cr
&= c \Big(S_{ N,1}(\osc_{ p,N}(f; x_0))\Big)_{ j},
\label{5.26}
\end{align}
where $ \osc_{ p,N}(f; x_0)$ stands for a  sequence defined as  
\[
\osc_{ p,N}(f; x_0)_i = \osc_{ p,N}(f; x_0, 2^ i),\quad  i\in \Z. 
\]

Using triangle inequality together with \eqref{5.26}, we obtain 
\begin{align}
& \osc_{ p,N-1}(f; x_0, 2^{ j})
\cr
& =   2 ^{ -j\frac{n }{p}}  \inf_{ P\in \mathcal{P}_{ N-1}}\| f - P\|_{ L^p(B(x_0, 2^j))}
\cr
& \le c 2 ^{ -j\frac{n }{p}} \| f - P_{ x_0, 2^j}^{ N}(f) + \dot{P}^{ N}_{ x_0, 2^j}(f)\|_{ L^p(B(x_0, 2^j))}
\cr
& \le c 2 ^{-j\frac{n }{p}}\| f - P_{ x_0, 2^j}^{ N}(f)\|_{ L^p(B(x_0, 2^j))}+c 2 ^{ -j\frac{n }{p}} \| \dot{P}_{ x_0, 2^j}^{ N}(f)\|_{ L^p(B(x_0, 2^j))} 
\cr
&\le c\osc_{ p,N}(f; x_0, 2^{ j}) +2^{j N} \| \dot{P}_{ x_0, 2^j}^{ N}(f)\|
\cr
&\le c\osc_{ p,N}(f; x_0, 2^{ j}) + c \Big(S_{ N,1}(\osc_{ p,N}(f; x_0))\Big)_{ j}. 
\label{5.27}
\end{align}
Noting that ${\D \osc_{ p,N}(f; x_0, 2^j) \le  S_{ N,1}(\osc_{ p,N}(f; x_0, 2^j))}$, we infer from \eqref{5.27}
\begin{equation}
\osc_{ p,N-1}(f; x_0)_j=\osc_{ p,N-1}(f; x_0, 2^j) \le c\Big(S_{ N,1}(\osc_{ p,N}(f; x_0))\Big)_{ j}, \quad j\in \Z. 
\label{5.28}
\end{equation}

Applying $ S_{s,  q}$ to both sides of \eqref{5.28}, and using Lemma\,\ref{lem10.1}, we get the inequality 
\[
|f|_{\mathscr{L}^{ s }_{q(p, N-1)}} =\sup_{ x_0\in \R^{n} } {\D  S_{ s, q}(\osc_{ p,N-1}(f; x_0) ) \le c 
\sup_{ x_0\in \R^{n} } S_{ s, q}(\osc_{ p,N}(f; x_0) )}
=|f|_{\mathscr{L}^{ s }_{q(p, N')}}, 
\]
which implies  \eqref{5.25}. 
We are now in a position to apply  \eqref{5.25} iteratively, replacing $ N$ by $ N+1$ to get 
\[
 |f|_{\mathscr{L}^{ s }_{q(p, N)}} \le c  |f|_{\mathscr{L}^{ s }_{q(p, N+1)}} 
 \le \ldots \le  c  |f|_{\mathscr{L}^{ s }_{q(p, N')}}. 
\]

This completes the proof of the lemma.  
\hfill \Beweisende 

\begin{rem}
\label{rem5.7} 
For all   $ f\in \mathscr{L}^{s}_{q (p, N)}(\R^{n}), 1 \le p < +\infty, 1 \le q \le +\infty, s\in [- \frac{n}{p}, N+1)$,  the condition  \eqref{5.23} is fulfilled, and therefore 
\eqref{5.24} holds for all $ f\in \mathscr{L}^{s}_{q (p, N)}(\R^{n})$ under the assumptions on $ p,q,s, N$ and $ N'$ 
of Lemma\, \ref{lem5.9}.  To verify this fact we observe for  $ f\in \mathscr{L}^{s}_{q (p, N)}(\R^{n})$ that   
\begin{equation}
\sup_{ m\in \Z} 2^{ -Nm} \osc_{ p, N}(f, 0, 2^m) \le |f|_{ \mathscr{L}^{s}_{q (p, N)}}. 
\label{5.29}
\end{equation} 
Then for $ L \in \N$, $ L>N$,  we estimate for multi index $ \alpha $ with $ |\alpha |=L$
\begin{align*}
  |D^\alpha \dot{P}_{0, 2^m}^L(f)| &= |D^{ \alpha }\dot{P}_{0, 2^m}^L((f- P^N_{ 0, 2^m})) | 
  \le c2^{ -Lm}\osc_{ p, N}(f, 0, 2^m) 
\\
&\le c 2^{m (N-L)} |f|_{ \mathscr{L}^{s}_{q (p, N)}}
\rightarrow 0\quad   \text{as}\quad m \rightarrow +\infty.   
\end{align*}
Hence, \eqref{5.23} is fulfilled.  
\end{rem}

\begin{rem}
 \label{rem5.7b} 
 In case $ q=\infty$,  since  $ \mathscr{L}^{s}_{\infty (p, N)}(\R^{n} ) $ coincides with the 
 usual Campanato space,  and Lemma\,\ref{lem5.9} is well known (cf. \cite[p. 75]{gia}).  
 
\end{rem}

A careful inspection of the proof of Lemma\,\ref{lem5.9} gives the following. 

\begin{cor}
 \label{cor3.10}
Let $ N, N'\in \N_0, N< N'$. Let $ f\in L^p_{ loc}(\R^{n} )$ satisfy  \eqref{5.23} with $ k=0$. Then, 
for all $ x_0\in \R^{n} $ and $ j\in \Z$ it holds,
\begin{equation}
\osc_{ p, N}(f; x_0, 2^j) \le c(S_{N+1 , 1}(   \osc_{ p, N'}(f; x_0) ))_j.  
\label{5.24a}
 \end{equation}  
\end{cor}

{\bf Proof}: Set $ k= N'-N$. 
Using  \eqref{5.28} with $ N'$ in place of $ N$, we find 
\begin{equation}
\osc_{ p,N'-1}(f; x_0, 2^j) \le c(S_{ N',1}(\osc_{ p,N'}(f; x_0)))_{ j}, \quad j\in \Z. 
\label{5.28a}
\end{equation} 
Iterating this inequality $ k$-times and applying Lemma\,\ref{lem10.1}, we arrive at 
\begin{align*}
  \osc_{ p,N}(f; x_0) &=  \osc_{ p,N'-k}(f; x_0)  \le cS_{ N+1,1} (S_{ N+2,1} \ldots S_{ N', 1}(\osc_{ p,N'}(f; x_0)))
  \\
  &\le cS_{ N+1,1}(   \osc_{ p, N'}(f; x_0)) .
  \end{align*}
Whence,  \eqref{5.24a}.  
\hfill \Beweisende

\vspace{0.3cm}
We also have the  following growth properties of functions in $\mathscr{L}^{s}_{q (p, N)}(\R^{n} ) $ as $ |x| 
\rightarrow +\infty$

\begin{lem}
 \label{growth}
 Let $ N\in \N_0$. Let  $f\in  \mathscr{L}^{s}_{q (p, N)}(\R^{n} ) , 1 \le q \le +\infty, 1 \le p < +\infty, 
 s \in  [N, N+1)$. 
 
 1. In case $ s\in (N, N+1)$ it holds
 \begin{equation}
  |f(x)| \le c (1+|x|^s) \|f\|_{ \mathscr{L}^{s}_{q (p, N)}}\quad \forall x\in \R^{n}. 
 \label{growth1}
  \end{equation} 
  
  2. In case $ s=N$ it holds
 \begin{equation}
  |f(x)| \le c (1+\log(1+|x|)^{ \frac{1}{q'}}|x|^{ N}) \|f\|_{ \mathscr{L}^{N}_{q (p, N)}}\quad \forall x\in \R^{n}. 
 \label{growth11}
  \end{equation} 
  Here $q'= \frac{q}{q-1}$,  $ c=\const>0$, depending on $ q,p,s, N$ and $ n$. 
 
\end{lem}

{\bf Proof}: {\it 1. The case  $ s\in (N, N+1)$.} Let $ x_0\in \R^{n} $. Let $ j\in \N_0$ such that 
$ 2^{ j} \le 1+|x_0| \le 2^{ j+1} $. Let $ \alpha $ be a multi index with $ |\alpha |=N$. 
Verifying that $ D^\alpha f(x_0)= \lim_{i \to -\infty}D^\alpha \dot{P}_{ x_0, 2^i} ^N(f) $, using triangle inequality 
we find 
\begin{align*}
|D^\alpha f(x_0)| &\le  \sum_{i=-\infty}^{j}  
|D^\alpha \dot{P}_{ x_0, 2^i} ^N(f) - D^\alpha \dot{P}_{ x_0, 2^{ i-1}} ^N(f)| + 
|D^\alpha \dot{P}_{ x_0, 2^j} ^N(f)|
\\
&\le  c\sum_{i=-\infty}^{j} 2^{ -iN} \osc_{ p, N}(f; x_0, 2^i)+ |D^\alpha \dot{P}_{ x_0, 2^j} ^N(f)|. 
\end{align*}

By the aid of H\"older's inequality we find 
\begin{align*}
\sum_{i=-\infty}^{j} 2^{ -iN} \osc_{ p, N}(f; x_0, 2^i) &= \sum_{i=-\infty}^{j} 2^{ -i (N-s)}2^{ -i s} \osc_{ p, N}(f; x_0, 2^i)  
\\
&\le c 2^{j(s-N)}  |f|_{ \mathscr{L}^{s}_{q (p, N)}} \le c (1+ |x_0|^{ s-N}) |f|_{ \mathscr{L}^{s}_{q (p, N)}}.   
\end{align*}
On the other hand, 
\begin{align*}
|D^\alpha \dot{P}_{ x_0, 2^j} ^N(f)| &= |D^\alpha \dot{P}_{ x_0, 2^j} ^N(f- P^N_{0, 2^{ j+1}}(f)|
+ |D^\alpha   (P^N_{0, 2^{ j+1}}(f) - P^N_{0, 1}(f))| + |D^\alpha P^N_{0, 1}(f)| 
\\
&\le 2^{ - jN- \frac{n}{p}} \|f- P^N_{0, 2^{ j+1}}(f)\|_{ L^p(x_0, 2^{ j+1})}   
+c \sum_{i=0}^{j} 2^{ -i (N-s)}2^{ -i s} \osc_{ p, N}(f; 0, 2^i)  
\\
& \qquad + c  \|f\|_{ L^p(B(1))}
\\
&\le  \osc_{ p, N}(f; 0, 2^{ j+1})  
+c \sum_{i=0}^{j} 2^{ -i (N-s)}2^{ -i s} \osc_{ p, N}(f; 0, 2^i)   +c \|f\|_{ L^p(B(1))}
\\
& \le c (1+|x_0|^{ s-N}) \|f\|_{ \mathscr{L}^{s}_{q (p, N)}}. 
\end{align*}
Accordingly, 
\begin{equation}
 \|D^N f(x)\| \le c(1+|x|^{ s-N})  \|f\|_{ \mathscr{L}^{s}_{q (p, N)}}. 
\label{growth2}
 \end{equation} 
This implies  \eqref{growth1}. 

\vspace{0.2cm}
{\it 2. The case  $ s=N$.} Let $ x_0\in \R^{n} $. As above we choose $ j\in \N_0$ such that $ 2^{ j} \le 1+|x_0| < 2^{ j+1}$.

In this case we first claim  
\begin{align}
\|D^N \dot{P}^N_{ x_0, 1}(f)\| \le (\log (1+|x_0|))^{ \frac{1}{q'}}  \|f\|_{ \mathscr{L}^{N}_{q (p, N)}}.
 \label{growth4}
\end{align}
Indeed, arguing as above using triangle inequality along with H\"older's inequality, we get 
\begin{align*}
\|D^N \dot{P}^N_{ x_0, 1}(f)\| &\le  \sum_{i=1}^{j}   \|D^N \dot{P}^N_{ x_0, 2^i}(f)\| - 
\|D^N \dot{P}^N_{ x_0, 2^{ i-1}}(f)\| + \|D^N \dot{P}^N_{ x_0, 2^{ j}}(f)\|
\\
&\le \sum_{i=1}^{j}   2^{ -Ni} \osc_{ p, N }(f; x_0, 2^i) + \|D^N \dot{P}^N_{ x_0, 2^{ j}}(f)\|
\\
&\le \sum_{i=1}^{j+1}   2^{ -Ni} \osc_{ p, N }(f; x_0, 2^i) + \|D^N \dot{P}^N_{0, 2^{ j+1}}(f)\|
\\
&\le c j^{ \frac{1}{q'}} |f|_{ \mathscr{L}^{N}_{q (p, N)}}  + \|D^N \dot{P}^N_{0, 2^{ j+1}}(f)\|.
\end{align*}
Similarly, 
\[
\|D^N \dot{P}^N_{0, 2^{ j+1}}(f)\| \le c j^{ \frac{1}{q'}} |f|_{ \mathscr{L}^{N}_{q (p, N)}}+ 
\|D^N \dot{P}^N_{0, 1}(f)\|.
\]
Combining the two inequalities we have just obtained, we get  \eqref{growth4}. 

Let $ i \in \Z$. Then by triangle inequality together with  \eqref{growth4} we find 
\begin{align*}
 &2^{ - \frac{n}{p} - iN}\|\dot{P}^N_{ x_0, 2^i}(f)\|_{ L^p(x_0, 2^i)} 
 \le c \|D^N \dot{P}^N_{ x_0, 2^i}(f)\|
\\
 &\le   c\sum_{l=i}^{1} \Big(\|D^N \dot{P}^N_{ x_0, 2^l}(f)\| -
  \|D^N \dot{P}^N_{ x_0, 2^{ l-1}}(f)\|\Big)
 +c \|D^N \dot{P}^N_{ x_0, 1}(f)\|
\\
 &\le   c\sum_{l=i}^{1} \Big(\|D^N \dot{P}^N_{ x_0, 2^l}(f)- D^N \dot{P}^N_{ x_0, 2^{ l-1}}(f)\| \Big)
 +c \|D^N \dot{P}^N_{ x_0, 1}(f)\|
\\
 &\le   c\sum_{l=i}^{1} 2^{ - Nl}\osc_{ p, N} (f; x_0, 2^l) +c \|D^N \dot{P}^N_{ x_0, 1}(f)\|
 \\
 &\le   c|i|^{ \frac{1}{q'}} |f|_{ \mathscr{L}^{N}_{q (p, N)}}+c  (\log (1+|x_0|))^{ \frac{1}{q'}}  \|f\|_{ \mathscr{L}^{N}_{q (p, N)}}.
 \end{align*}
 This shows that 
 \begin{align}
& 2^{ - i(N-1)}\osc_{ p,N-1} (f; x_0, 2^i ) 
\cr
&\le 2^{ - i(N-1)}\osc_{ p,N} (f; x_0, 2^i ) +  
2^{ - \frac{n}{p} }\|\dot{P}^N_{ x_0, 2^i}(f)\|_{ L^p(x_0, 2^i)} 
    \cr
& \le 2^{ - i(N-1)}\osc_{ p,N} (f; x_0, 2^i )  + c 2^{ i} \Big(|i|^{ \frac{1}{q'}} + (\log (1+|x_0|))^{ \frac{1}{q'}} \Big)
 \|f\|_{ \mathscr{L}^{N}_{q (p, N)}}.
  \label{5.25e} 
 \end{align}
Summing both sides over $ i =-\infty$ to $ i=1$ and applying H\"older's inequality, we get 
\begin{equation}
  \sum_{i=-\infty}^{1} 2^{ -i (N-1)}\osc_{ p,N-1} (f; x_0, 2^i ) \le c\Big(1 + (\log (1+|x_0|))^{ \frac{1}{q'}} \Big)
 \|f\|_{ \mathscr{L}^{N}_{q (p, N)}}.
 \label{5.25d}
 \end{equation} 
 
 Let $ \alpha $ be a multi index with $ |\alpha |=N-1$. Noting that $ D^\alpha f(x_0) = \lim_{ i\to -\infty} D^\alpha  \dot{ P}^{ N-1}_{ x_0, 2^i}(f)$, using triangle inequality together with  \eqref{5.25d}, we 
 infer  
  \begin{align*}
  | D^\alpha  f(x_0)| &\le |D^\alpha  \dot{ P}^{ N-1}_{ x_0, 2^i}(f) | + c 
   \sum_{i=-\infty}^{1}|D^\alpha  \dot{ P}^{ N-1}_{ x_0, 2^i}(f) - D^\alpha  \dot{ P}^{ N-1}_{ x_0, 2^{ i-1}}(f)|
  \\
 & \le    |D^\alpha  \dot{ P}^{ N-1}_{ x_0, 2^i}(f) | + c 
   \sum_{i=-\infty}^{1} 2^{ - (N-1)}\osc_{ p, N-1}(f; x_0, 2^i)
 \\  
 &\le  \|D^{ N-1}  \dot{ P}^{ N-1}_{ x_0, 2^i}(f) \| +c\Big(1 + (\log (1+|x_0|))^{ \frac{1}{q'}} \Big)
 \|f\|_{ \mathscr{L}^{N}_{q (p, N)}}. 
  \end{align*}
 
Arguing as above using triangle inequality, using  \eqref{5.25e}, we find 
\begin{align*}
\|D^{ N-1}  \dot{ P}^{ N-1}_{ x_0, 2^i}(f) \| &\le  c\sum_{i=0}^{j} 2^{ - (N-1)i}\osc_{ p, N-1}(f, x_0, 2^i)   
+ \|D^{ N-1}  \dot{ P}^{ N-1}_{ x_0, 2^j}(f) \|
\\
&\le  c\sum_{i=0}^{j} 2^{ - (N-1)i}\osc_{ p, N-1}(f, x_0, 2^i)    +\sum_{i=0}^{j+1} 2^{ - (N-1)i}\osc_{ p, N-1}(f, 0, 2^i)   
\\
& \qquad +    \|D^{ N-1}  \dot{ P}^{ N-1}_{0, 1}(f) \|
\\
&\le c 2^{ j} j^{ \frac{1}{q'}} \|f\|_{ \mathscr{L}^{N}_{q (p, N)}}
\le c(1+ \log(1+|x_0|)^{ \frac{1}{q'}} |x_0|)\|f\|_{ \mathscr{L}^{N}_{q (p, N)}}.  
\end{align*}
Combining the above inequalities we obtain 
\[
|D^{ N-1} f(x_0)| \le  (1+ \log(1+|x_0|)^{ \frac{1}{q'}} |x_0|)\|f\|_{ \mathscr{L}^{N}_{q (p, N)}}.    
\]
This yields  \eqref{growth11}.  \hfill \Beweisende

\vspace{0.3cm}
Using the Poincar\'e's inequality and Lemma\,\ref{lem5.9}, we get the following embedding.
 
\begin{lem}
\label{lem5.8}  Let $ N \in \N_0, k\in \N_0, 1< p< +\infty, 1 \le q \le +\infty, s\in [N, N+1)$. 

1. In case $ q=\infty$ and $ s\notin \N$ it holds 
\begin{equation}
 \mathscr{L}^{ k, s}_{ \infty (p, N)} (\R^{n})  \cong C^{ k +N, s- N}  (\R^{n}). 
\label{5.29a}
\end{equation}

2. In case $ q=\infty$ and $ s\in \N$ it holds 
\begin{equation}
 \mathscr{L}^{ k, s}_{ \infty (p, N)} (\R^{n})  \cong BMO_{ k+s}  (\R^{n}). 
\label{5.29b}
\end{equation}
where 
\[
 BMO_N  = \Big\{f\in L^1_{ loc}(\R^{n} ) \,\Big|\sup_{j\in \Z} 2^{ -Nj}\osc_{ 1, N}(f; x_0, 2 ^j) < +\infty\Big\}.
\]

3. In case $ 1 \le q<\infty$ it holds 
\begin{equation}
 \mathscr{L}^{ k, s}_{ q (p, N)} (\R^{n}) \hookrightarrow \mathscr{L}^{ k+s}_{ q (p, N+k)} (\R^{n})
\hookrightarrow \mathscr{L}^{ k, s}_{ \infty (p, N)} (\R^{n}).
\label{5.29c}
\end{equation}
\end{lem}

{\bf Proof}: 1. In case $ k=0$ the space $ \mathscr{L}^{s}_{ \infty (p, N)} (\R^{n})$ coincides with 
the Campanato space $ \mathscr{L}^{p, N n + p(s-N )}_{N} (\R^{n})$ which is isomorphic to $ C^{ N, s-N}(\R^{n})$, 
(cf. \cite[Chap. III 1.]{gia}, \cite{ca}). In case, $ k \ge 1$. 
For $ f\in \mathscr{L}^{k, s}_{ \infty (p, N)} (\R^{n})$ 
we get $ D^k f \in \mathscr{L}^{s}_{ \infty (p, N)} (\R^{n}) \cong C^{ N, s-N}(\R^{n})$, which shows \eqref{5.29a}.  

\vspace{0.3cm}
2.  In case $ k=0$, and $ s=N$ the space $ \mathscr{L}^{s}_{ \infty (p, N)} (\R^{n})$ coincides with 
the Campanato space $ \mathscr{L}^{p, N n}_{N} (\R^{n})$. According to \cite[Chap. III,1.]{gia} this space coincides 
with the space $ BMO_N$. In case $ k \ge 1$ we argue as above to verify \eqref{5.29b}.

\vspace{0.3cm}
3. Let $ \mathscr{L}^{ k, s}_{ q (p, N)} (\R^{n})$. Using Poincar\'e inequality \eqref{5.5} with $ j=k$, we find 
${\D \osc_{ p, N+k}(f; x_0, 2^j ) } $$\le c 2^{jk} {\D \osc_{ p, N}(D^k f; x_0, 2^j ) }$. 
Accordingly, 
\[
\| \{2^{ -(s+k)j} \osc_{ p, N+k}(f; x_0, 2^j )\}_{ j\in \Z} \|_{ \ell^q}\le c 
\|\{ 2^{ -sj}\osc_{ p, N}(D^k f; x_0, 2^j)\}_{ j\in \Z}\|_{ \ell^q}, 
\]
where 
\[
\osc_{ p, N}(f; x_0 ) = \{\osc_{ p, N}(f; x_0, 2^j )\}_{ j\in \Z}.
\]
Taking the supremum  over all $ x_0\in \R^{n}$ on both sides of the above estimate, we get the first embedding.  

\vspace{0.3cm}
It remains to show the second embedding. To see this we first notice that $ \mathscr{L}^{ k+s}_{ q (p, N+k)} (\R^{n})
\hookrightarrow \mathscr{L}^{ k+s}_{ \infty (p, N+k)} (\R^{n})$. Indeed, 
\[
2^{ -(s+k)j} \osc_{ p, N+k}(f; x_0, 2^j) \le 2^{-(s+k)j} \Big(S_{ k+s, q} (\osc_{ p, N+k}(f; x_0))\Big)_j \le | f|_{ \mathscr{L}^{ k+s}_{ q (p, N+k)}}. 
\]
Taking the supremum over all $ j\in \Z$ and  $ x_0\in \R^{n}$, we get the embedding
\[
\mathscr{L}^{ k+s}_{ q (p, N+k)} (\R^{n})
\hookrightarrow \mathscr{L}^{ k+s}_{ \infty (p, N+k)} (\R^{n}).
\]
On the other hand, in case $ s\in (N, N+1)$, from  \eqref{5.29a}
it follows 
$ \mathscr{L}^{ k+s}_{ \infty (p, N+k)} (\R^{n})\cong C^{ k+N, s-N}(\R^{n}) \cong 
\mathscr{L}^{ k, s}_{ \infty (p, N)} (\R^{n}) $. In case $ s=N$ using \eqref{5.29b}, we  also get  
$ \mathscr{L}^{ k+N}_{ \infty (p, N+k)} (\R^{n}) \cong 
\mathscr{L}^{ k, N}_{ \infty (p, N)} (\R^{n}) $. This shows desired embedding.
 \hfill \Beweisende

\vspace{0.5cm}  
Using Gagliardo-Nirenberg's inequalities,  we can get the interpolation properties. First let us recall the  
Gagliardo-Nirenberg inequalities.

\begin{lem}
\label{lem5.10}
Let $ j, N \in \N_0, 0 \le  j < k$.   Let $ 1 \le p, p_0, p_1 \le +\infty $, and $ \theta \in  \Big[\frac{j}{N}, 1\Big]$, satisfying  
\begin{equation}
\frac{1}{p}= \frac{j}{n} + \frac{1-\theta }{p_0} + \Big(\frac{1}{p_1} - \frac{k}{n}\Big)\theta.
\label{5.31}
\end{equation}
 Then, for all $ f \in  L^{p_0}(B(1)) \cap  W^{k,\, p_1}(B(1))$ it holds 
\begin{equation}
\| D^j f\|_{ L^p(B(1))} \le c \|  f\|_{ L^{p_0}(B(1))}^{ 1-\theta } \|  f\|_{ W^{k,\, p_1}(B(1))} ^\theta. 
\label{5.32}
\end{equation}
\end{lem}

Notice that, using the generalized Poincar\'e inequality, under the assumption of Lemma\,\ref{lem5.10}, for all $ f \in  L^{p_0}(B(1)) \cap  W^{k,\, p_1}(B(1))$, and $N \in \N_0, N \ge k-1$ the  following inequality holds 
 \begin{equation}
\| D^j (f - P^N_{ 0, 1}(f)) \|_{ L^p(B(1))} \le c \|  f- P^N_{ 0, 1}(f))\|_{ L^{p_0}(B(1))}^{ 1-\theta } \| D^k f- D^k P^N_{ 0, 1}(f)\|_{ L^{ p_1}(B(1))} ^\theta. 
\label{5.33}
\end{equation}
By a standard scaling and translation argument, we deduce from \eqref{5.33} that for all $ x_0 \in \R^{n}, 0< r< +\infty$, 
$N \in \N_0, N \ge k-1$, and for all  $ f \in  L^{p_0}(B(x_0, r)) \cap  W^{k,\, p_1}(B(x_0,r))$ the following inequality holds
\begin{align}
&\| D^j f - P^{ N-j}_{ x_0, r}( D^j f)\|_{ L^p(B(x_0, r))}  
\cr
&\quad  =\| D^j (f - P^N_{ x_0, r}(f)) \|_{ L^p(B(x_0, r))}
\cr
 &\quad \le c \|  f- P^N_{x_0, r}(f))\|_{ L^{p_0}(B(x_0,r))}^{ 1-\theta } 
 \| D^k f-P^{ N-k}_{x_0, r}(D^k f)\|_{ L^{ p_1}(B(x_0, r))} ^\theta. 
\label{5.34}
\end{align}

\begin{thm}
\label{thm5.10}
Let $ j, k, N \in \N_0, 0 \le  j < k \le N+1$.   Let $ 1 \le p, p_0, p_1<+\infty $, $ 1 \le  q, q_0, q_1 \le +\infty, 
-\infty< s, s _0, s _1 < N+1 $, 
and $ \theta \in  \Big[\frac{j}{N}, 1\Big]$, satisfying  
\begin{align}
\frac{1}{p} &= \frac{j}{n} + \frac{1-\theta }{p_0} + \Big(\frac{1}{p_1} - \frac{k}{n}\Big)\theta,
\label{5.35}
\\
\frac{1}{q} &= \frac{1-\theta }{q_0} +  \frac{\theta }{q_1},
\label{5.36}
\\
s +j  &= (1-\theta) s_0 + \theta( s_1+k). 
\label{5.37}
\end{align}
Then, for all $ \mathscr{L}^{  s_0}_{  q_0 (p_0, N) } (\R^{n})\cap \mathscr{L}^{ k, s_1}_{  q_1 (p_1, N) }(\R^{n})$ 
it holds  
\begin{equation}
\| f\|_{ \mathscr{L}^{ j, s}_{q (p, N-j) }} \le c
\| f\|_{ \mathscr{L}^{  s_0}_{  q_0 (p_0, N) }}^{ 1-\theta } \| f\|_{ \mathscr{L}^{ k, s_1}_{ q_1 (p_1, N-k) }}^{ \theta }. 
\label{5.38}
\end{equation}
\end{thm} 

{\bf Proof}:  Observing \eqref{5.35} and \eqref{5.36},  thanks to \eqref{5.34} we find 
\begin{align*}
& 2^{ -s l}  \osc_{ p, N-j}(D^j f; x_0, 2^{ l} ) 
\\
&\le c2^{- ls_0 (1-\theta )-  l s_1\theta    } \osc_{ p_0, N} (f; x_0, 2^{ l})^{ 1-\theta } \osc_{ p_1, N-k}(D^k f; x_0, 2^{ l} )]^{ \theta }  
\\
&= c  [2^{- ls_0    } \osc_{ p_0, N} (f; x_0, 2^{ l} )]^{ 1-\theta }
[2^{-  l s_1 }  \osc_{ p_1, N-k}(D^k f; x_0, 2^{ l} )]^{ \theta }.   
\end{align*} 
According to \eqref{5.37},  we may apply  $ \ell^q$ norm to both sides of the above inequality and use  H\"older's inequality. This 
gives  
\begin{align*}
& \Big(\sum_{l\in \Z}(2^{ -s l}  \osc_{ p, N-j}(D^j f; x_0, 2^{ l} ))^q \Big)^{ \frac{1}{q}}
\\
& \le  c \Big(\sum_{l\in \Z} (2^{- ls_0    } \osc_{ p_0, N} (f; x_0, 2^{ l} ))^{ q_0} \Big)^{ \frac{1-\theta }{q_0}}
\Big(\sum_{l\in \Z} (2^{-  l s_1 }  \osc_{ p_1, N-k}(D^k f; x_0, 2^{ l} ))^{ q_1} \Big)^{ \frac{\theta }{q_1}}.
\end{align*} 
Taking the supremum over all $ x_0\in \R^{n}$, we get the assertion \eqref{5.38}.  \hfill \Beweisende  

\begin{rem}
\label{rem5.11}
Consider the special case 
\begin{align}
\begin{cases}
 N=k, p=p_0=p_1,  \theta = \frac{j}{k}, s = s _0=s _1=0,   
\\[0.3cm]
 1 \le  q< +\infty, q_0=+\infty, q_1 = \frac{qk}{j}.   
\end{cases}
\label{5.39}
\end{align}
Then, \eqref{5.38} reads
\begin{equation}
\| f\|_{ \mathscr{L}^{ j, 0}_{q (p, k-j) }} \le c
\| f\|_{ \mathscr{L}^{  0}_{\infty (p, k) }}^{ 1- \frac{j}{k}} 
\| f\|_{ \mathscr{L}^{ k, 0}_{\frac{qk}{j} (p, 0) }}^{ \frac{j}{k} }
\le c \| f\|_{ BMO}^{ 1- \frac{j}{k}} \| f\|_{ \mathscr{L}^{ k, 0}_{q (p, 0) }}^{ \frac{j}{k} }. 
\label{5.40}
\end{equation}
Under the assumption that 
\begin{equation}
\lim_{m\to \infty} \dot{P}^L_{ 0, 2^m}(D^j u) =0\quad \forall\,L=1, \ldots, k-j,
 \label{5.41}
 \end{equation} 
we estimate the term on the left hand side by the aid of \eqref{5.24} with $ N=0$ and $ N'=k-j$. This yields
\begin{equation}
\| f\|_{ \mathscr{L}^{ j, 0}_{q (p, 0) }} \le c
 \| f\|_{ BMO}^{ 1- \frac{j}{k}} \| f\|_{ \mathscr{L}^{ k, 0}_{ q (p, 0) }}^{ \frac{j}{k} }. 
\label{5.42}
\end{equation} 
\end{rem}

We are now in a position to prove the following product estimate.

\begin{thm}
\label{thm5.12}
Let $ 1 < p< +\infty$. Let $ N\in \N_0$ and  $s\in (-\infty, N+1)$. Then for all  $ f, g \in 
\mathscr{L}^{ k, s}_{  q (p, N)}(\R^{n})\cap L^\infty(\R^{n})$, 
it holds 
\begin{equation}
\| fg\|_{  \mathscr{L}^{ k, s}_{  q (p, N)}} \le c \Big(\| f\|_{\infty}\| g\|_{  \mathscr{L}^{ k, s}_{ q (p, N)}}
+ \| g\|_{ \infty}\| f\|_{  \mathscr{L}^{ k, s}_{  q (p, N)}}\Big). 
\label{5.44}
\end{equation}
\end{thm}

{\bf Proof}:  Let $ \alpha , \beta \in N_0^n$ two multi index both are  not zero with $ | \alpha + \beta |=k$. Set $ | \alpha |= j$. 
Using triangle inequality, we see that 
\begin{align}
& \| D^{ \alpha } f D^\beta g - P^{ N+k-j}_{ x_0, r}(D^{ \alpha } f )P^{ N+j}_{ x_0, r}( D^\beta g) \|_{ L^p(B(x_0, r))}
\cr
&\quad \le c \| (D^{ \alpha } f - P^{ N+k-j}_{ x_0, r}(D^{ \alpha } f )) (D^\beta g-P^{ N+j}_{ x_0, r}( D^\beta g)) \|_{ L^p(B(x_0, r))}
\cr
&\qquad  + c \| (D^{ \alpha } f - P^{ N+k-j}_{ x_0, r}(D^{ \alpha } f ))  P^{ N+j}_{ x_0, r}(D^\beta g) \|_{ L^p(B(x_0, r))}
\cr
&\qquad + c \| P^{ N+k-j}_{ x_0, r}(D^{ \alpha } f ))  (D^\beta g-P^{ N+j}_{ x_0, r}( D^\beta g) )\|_{ L^p(B(x_0, r))}
\cr
&\quad = I+II+III.
\label{5.45}
\end{align}
Using H\"older's inequality together with Gaglirdo-Nirenberg's inequality  \eqref{5.34}, we estimate 
\begin{align*}
I &\le c\| D^{ \alpha } f - P^{ N+k-j}_{ x_0, r}(D^{ \alpha } f )\|_{ L^{ \frac{k}{j}p}(B(x_0, r))}
\|D^\beta g-P^{ N+j}_{ x_0, r}( D^\beta g) \|_{ L^{ \frac{k}{k-j}p}(B(x_0, r))}
\\
&= c\| D^{ \alpha} f - D^{ \alpha }P^{ N+k}_{ x_0, r} (f)\|_{ L^{ \frac{k}{j}p}(B(x_0, r))}
\|D^\beta g- D^\beta P^{ N+k}_{ x_0, r}(g) \|_{ L^{ \frac{k}{k-j}p}(B(x_0, r))}
\\
&\le c\| D^{ j} (f - P^{ N+k}_{ x_0, r} (f))\|_{ L^{ \frac{k}{j}p}(B(x_0, r))}
\|D^{ k-j} (g- P^{ N+k}_{ x_0, r}(g)) \|_{ L^{ \frac{k}{k-j}p}(B(x_0, r))}
\\
&\le c\|f - P^{ N+k}_{ x_0, r} (f) \|_{ L^{ \infty}(B(x_0, r))}^{ 1-\frac{j}{k}}
\| D^{ k} (f - P^{ N+k}_{ x_0, r} (f))\|_{ L^{p}(B(x_0, r))}^{ \frac{j}{k}}\times 
\\
&\qquad \times  c\|g - P^{ N+k}_{ x_0, r} (g) \|_{ L^{ \infty}(B(x_0, r))}^{ \frac{j}{k}}
\| D^{ k} (g - P^{ N+k}_{ x_0, r} (g))\|_{ L^{p}(B(x_0, r))}^{ 1-\frac{j}{k}}
\\
&\le c\|f\|_{ L^{ \infty}(B(x_0, r))}^{ 1-\frac{j}{k}} 
\| D^{ k} f - P^{ N}_{ x_0, r} (D ^k f))\|_{ L^{p}(B(x_0, r))}^{ \frac{j}{k}}\times 
\\
&\qquad \times \|g\|_{ L^{ \infty}(B(x_0, r))}^{ \frac{j}{k}} 
\| D^{ k} g - P^{ N}_{ x_0, r} (D ^k g))\|_{ L^{p}(B(x_0, r))}^{ 1-\frac{j}{k}}.
\end{align*}
Applying Young's inequality, we obtain 
\begin{align*}
I &\le c\|f\|_{ L^{ \infty}(B(x_0, r))} \| D^{ k} g - P^{ N}_{ x_0, r} (D ^k g)\|_{ L^{p}(B(x_0, r))} 
\\
&\qquad + 
c\|g\|_{ L^{ \infty}(B(x_0, r))} \| D^{k} f - P^{ N}_{ x_0, r} (D ^k f)\|_{ L^{p}(B(x_0, r))}. 
\end{align*}
In order to estimate $ II$ we make use of the inequality 
\[
\| P^{N+ j}_{ x_0, r} (D^\beta   g)\|_{ L^\infty(B(x_0, r))} \le c r^{-(k-j)} \| g\|_{ L^\infty(B(x_0, r))},  
\]
 which can be proved by a standard scaling argument.  Together with Poincar\'e's inequality we find 
 \begin{align*}
II &\le cr^{ k-j} \|D^{ k-j} (D^{ \alpha } f - P^{ N+k-j}_{ x_0, r}(D^{ \alpha } f ))   \|_{ L^p(B(x_0, r))}
r^{ -(k-j)}\| g\|_{ L^\infty(B(x_0, r))}
 \\
&\le  c\|g\|_{ L^{ \infty}(B(x_0, r))} \| D^{ k} f - P^{ N}_{ x_0, r} (D ^k f)\|_{ L^{p}(B(x_0, r))}. 
 \end{align*}
 By an analogous reasoning we get 
 \begin{align*}
III \le   c\|f\|_{ L^{ \infty}(B(x_0, r))} \| D^{ k} g - P^{ N}_{ x_0, r} (D ^k g)\|_{ L^{p}(B(x_0, r))}. 
 \end{align*}
 Inserting the estimates of $ I, II$ and $ III$ into the right-hand side of \eqref{5.45}, we arrive at 
 \begin{align}
& \| D^{ \alpha } f D^\beta g - P^{N+ k-j}_{ x_0, r}(D^{ \alpha } f )P^{ N+j}_{ x_0, r}( D^\beta g) \|_{ L^p(B(x_0, r))}
 \cr
 &\quad \le c\|f\|_{ L^{ \infty}(B(x_0, r))} \| D^{ k} g - P^{ N}_{ x_0, r} (D ^k g)\|_{ L^{p}(B(x_0, r))} 
\cr
&\qquad + 
c\|g\|_{ L^{ \infty}(B(x_0, r))} \| D^{ k} f - P^{ N}_{ x_0, r} (D ^k f)\|_{ L^{p}(B(x_0, r))}. 
 \label{5.46}
 \end{align}
Let $ \gamma \in \N_0$ be a multi index with $ | \gamma |=k$.   Using Leibniz formula, we compute 
 \[
D^\gamma  (fg) =  \sum_{ \alpha + \beta=\gamma } \binom{\gamma !}{\alpha ! \beta !} D^\alpha f D^\beta g. 
\]
 Thus,  employing Corollary\,\ref{cor5.7}, using triangle inequality together with \eqref{5.46}, we obtain  
 \begin{align*}
 &\|  D^\gamma  (fg) - P^{ 2N+k}_{ x_0, r}(D^\gamma (fg))\|_{ L^p(B(x_0, r))}
\\
&\le c \inf_{ Q\in {\cal P}_{2N+k }}\|  D^\gamma  (fg) -Q\|_{ L^p(B(x_0, r))}
\\
 & \le c\Big\|D^\gamma  (fg) -
 \sum_{ \alpha + \beta =\gamma } \binom{\gamma !}{\alpha ! \beta !}  
  P^{ N+k-j}_{ x_0, r}(D^{ \alpha } f )P^{ N+j}_{ x_0, r}( D^\beta g)\Big\|_{ L^p(B(x_0, r))}
 \\
&=  c\Big\|\sum_{ \alpha + \beta =\gamma } \binom{\gamma !}{\alpha ! \beta !} (D^\alpha f D^\beta g -
P^{ N+k-j}_{ x_0, r}(D^{ \alpha } f )P^{ N+j}_{ x_0, r}( D^\beta g))\Big\|_{ L^p(B(x_0, r))}
\\
&\le 
c\|f\|_{ \infty} \| D^{ k} g - P^{ N}_{ x_0, r} (D ^k g)\|_{ L^{p}(B(x_0, r))} 
\cr
&\qquad + 
c\|g\|_{\infty} \| D^{ k} f - P^{ N}_{ x_0, r} (D ^k f)\|_{ L^{p}(B(x_0, r))}.
 \end{align*} 
 This yields the product estimate 
 \begin{equation}
 \osc_{p, 2N+k}(D^{ k }(fg); x_0, r) \le c \|f\|_{ \infty} \osc_{p, N}( D^k g; x_0, r)  + 
 c\|g\|_{ \infty} \osc_{p, N}(D^k f; x_0, r).    
 \label{5.47}
 \end{equation}
Into \eqref{5.47} we insert $ r= 2^{ j}, j\in \Z$, and multiply this by $ 2^{ -sj}$. Then,  
applying the $ \ell^q$ norm to both sides of \eqref{5.47}, we are led to 
 \begin{equation}
 \| fg\|_{ \mathscr{L}^{ k,s}_{ q (p, 2N+k)}}  \le  c \Big(\|f\|_{ \infty} \| g\|_{ \mathscr{L}^{ k,s}_{ q (p, N)}} 
 +  \|g\|_{ \infty}\| f\|_{ \mathscr{L}^{ k,s}_{ q (p, N)}}  \Big). 
 \label{5.48}
  \end{equation}
 Verifying  \eqref{5.23} holds for $ N'=2N+k$,   
 we are in a position to apply Lemma\,\ref{lem5.9} with $ N'=2N+k$. This gives \eqref{5.44}.  \hfill \Beweisende    

%%% ----------------------------------------------------------------------
%       SECTION 2
%%% ----------------------------------------------------------------------
\section{Proof of the main theorems}
\label{sec:-2}
\setcounter{secnum}{\value{section} \setcounter{equation}{0}
\renewcommand{\theequation}{\mbox{\arabic{secnum}.\arabic{equation}}}}

We start with the following energy identity for solutions to the transport equation.  Let $ 1 < p < +\infty$, $ x_0\in \R$ and 
$ 0< r<+\infty$. We denote $ \varphi _{ x_0, r}= \varphi (r^{ -1}(x_0- \cdot ))$. 
We define the following minimal polynomial $ P^{ N,\ast}_{ x_0, r}(f)$, $ f\in L^p(B(x_0, r))$, 
 by 
\begin{equation}
\|( f- P^{ N,\ast}_{ x_0, r}(f)) \varphi_{ x_0, r}  \|_{ p} 
= \min_{ Q\in \mathcal{P}_N} \|( f- Q) \varphi_{ x_0, r}\|_{ p}.   
\label{2.1}
\end{equation}
The existence and uniqueness of such polynomial is shown in appendix of the paper. 

We recall the notation  $  \varphi _{ x_0, r}=  r^{ -n} \varphi (r^{ -1}(x_0 -\cdot ))$. We have the following.

\begin{lem}
\label{lem2.1}  Given $v\in  L^1(0, T; C^{ 0,1}(\R^{n}; \R^{n} ))$,  and 
$ g \in L^1(0, T; L^p_{ loc}(\R^{n} ))$,  
let $ f\in L^\infty(0, T; C^{ 0,1}(\R^{n}))
\cap C([0,T]; L^p_{ loc}(\R^{n}))$ be a weak  solution  to the transport equation 
\begin{equation}
\partial _t f + (v\cdot \nabla) f = g\quad  \text{ in}\quad  Q_T. 
\label{2.2}
\end{equation}
Let $ N \in \N_0$. Define,
\[
L=\begin{cases}
2N-1\quad   &\text{if}\quad N \ge 1
\\[0.3cm]
0\quad  &\text{if}\quad N =0.\end{cases}
\]
Then for all $ t\in [0,T]$ it holds
\begin{align}
& e(t)
= e(0)
 +  \intl_{0}^{t}  v\cdot \nabla \varphi _{ x_0, r}  
| f- P^{ L,\ast}_{ x_0, r}(f)|^{ p}\varphi _{ x_0, r}^{ p-1} e(\tau )^{ 1-p}
 dx  d\tau 
 \cr
 &\qquad +  \frac{1}{p}\intl_{0}^{t}  \intl_{B(x_0,r)} \nabla \cdot v  |f- P^{ L,\ast}_{ x_0, r}(f)|^{ p} \varphi_{ x_0, r}^{ p}  
 e(\tau )^{1- p} dxd\tau 
 \cr
 &\qquad +   \intl_{0}^{t}  \intl_{B(x_0,r)} 
 v\cdot \nabla P^{ L,\ast}_{ x_0, r}(f(\tau )) \cdot   
| f- P^{ L,\ast}_{ x_0, r}(f(\tau ))|^{ p-2}
(f- P^{ L,\ast}_{ x_0, r}(f))\varphi_{ x_0,r} ^p  
e(\tau )^{1- p} dxd\tau 
\cr
&\qquad  +  \intl_{0}^{t}
 \intl_{B(x_0, r)} (g-P_{ x_0,r}^{ N}(g)) | f- P^{ L,\ast}_{ x_0, r}(f)|^{ p-2}
  (f- P^{ L,\ast}_{ x_0, r}(f))\varphi_{ x_0, r} ^p 
e(\tau )^{1- p}dx d\tau
\cr
&= e(0) + I+II+III+IV, 
 \label{2.3}
  \end{align}
  where 
 \[
e(\tau )=  \|(f(\tau )- P^{ L,\ast}_{ x_0, r}(f(\tau ))) \varphi_{ x_0, r}  \|_{ p},\quad  \tau \in [0, T]. 
\]
In addition, the following inequality  holds  for all $ t\in [0,T]$
 \begin{align}
& \osc_{ p, L} \Big(f(t); x_0,  \frac{r}{2}\Big)
\le   c \osc_{ p, L} (f(0); x_0,  r) +c r^{ -1} \intl_{0}^{t} \|   v(\tau )\|_{ L^\infty(B(x_0, r))}  
\osc_{ p, N} (f(\tau ); x_0,  2r)d\tau  
 \cr
&\quad  +c \intl_{0}^{t}   \|  \nabla \cdot  v(\tau )\|_{ L^\infty(B(x_0, r))} \osc_{ p, N} (f(\tau ); x_0,  2r)d\tau  
\cr
& \quad + \delta_{ N0} c  \intl_{0}^{t}  \osc_{ p, N} (v(\tau ); x_0,  r) 
\| \nabla P^{ N}_{ x_0, r}(f(\tau ))\|_{ L^\infty(B(x_0, r))} d\tau 
\cr
&\quad +  c\intl_{0}^{t} \osc_{ p, N} (g(\tau ); x_0,  r)  d\tau,
  \label{2.4}
 \end{align}       
where $ \delta _{ N0}=0$ if $ N=0$ and $ 1$ otherwise. 

\end{lem}

{\bf Proof}: 
Let $ x_0 \in \R^{n}, 0< r< +\infty$ be fixed.  Let  $ \delta \ge 0$ we define  
\[
F_{ \delta }( z) = (\delta + | z|^2)^{ \frac{p-2}{2}} z,\quad  z\in \R^{n}.  
\]
Let $ N\in \N_0$.  Set $ L=0$ if $ N=0$ and $ L= 2N-1$ if $ L \ge 1$.  For $ \delta >0$ by  
$P^{L, \delta }_{ x_0, r}(f(\tau )) \in  \mathcal{P}_{ L}$, $ 0 \le \tau \le T$,  we denote the minimal polynomial, defined in the Appendix A. (cf. Lemma\,\ref{lemA.1}), such that 
\begin{equation}
\intl_{B(x_0, r)} F_\delta  (f(\tau )-P^{L, \delta }_{ x_0, r}(f(\tau ) )\cdot Q \varphi _{ x_0, r} ^p dx =0
\quad \forall\,\tau \in [0,T],\quad  \forall\,Q\in \mathcal{P}_{ L}.   
\label{2.5}
\end{equation}
Furthermore, for all $ \tau \in [0,T]$ it holds 
\begin{equation}
P^{L, \delta }_{ x_0, r}(f(\tau )) \rightarrow P^{ L,\ast}_{ x_0, r}(f(\tau ))\quad  \text{ in}\quad L^p(B(x_0,r))\quad  
\text{ as}\quad  \delta \searrow  0. 
\label{2.6}
\end{equation}

According to 
\eqref{A.2c} the function $s \mapsto  P^{ L, \delta } _{ x_0, r}(f(s))$  is differentiable for $ \delta >0$, and from \eqref{2.2}  we get 
\begin{align}
 & \partial _t (f- P^{ L, \delta } _{ x_0, r}(f)) + (v\cdot \nabla) (f-P^{ L, \delta } _{ x_0, r}(f))+ 
(v\cdot \nabla) P^{ L, \delta } _{ x_0, r}(f)
\cr
 & \qquad = g - \partial _t P^{ L, \delta } _{ x_0, r}(f)\quad  
\text{ in}\quad  Q_T. 
\label{2.7}
\end{align}
First let us verify that $ \partial _t P^{ L, \delta } _{ x_0, r}(f(\tau ))\in {\cal P}_L$ for all $ \tau \in [0,T]$. In fact, for any multi index 
$ \alpha \in \N_0$ with $ |\alpha |=L+1$, recalling $ P^{ L, \delta } _{ x_0, r}(f)\in {\cal P}_L$, we get $ D^\alpha \partial _t P^{ L, \delta } _{ x_0, r}(f)= \partial _t D^\alpha P^{ L, \delta } _{ x_0, r}(f) =0$. This shows the claim.

 We multiply \eqref{2.7} by $ F_\delta  (f(\tau )- P^{ L, \delta } _{ x_0, r}(f(\tau ))\varphi_{ x_0, r} ^p$,  
 integrate over $ B(x_0, r)$ and apply  integration by parts. This together with \eqref{2.5} yields 
 \begin{align*}
& \partial _t  \|(\delta +| f(\tau )- P^{ L, \delta } _{ x_0, r}(f(\tau) )|^2)^{ \frac{1}{2}} \varphi_{ x_0,r} \|_{ p} 
\|(\delta +| f(\tau )- P^{ L, \delta } _{ x_0, r}(f(\tau ))|^2)^{ \frac{1}{2}} \varphi_{ x_0,r} \|^{ p-1}_{ p} 
\\
&= \frac{1}{p} \partial_t  \|(\delta +| f(\tau )- P^{ L, \delta }_{ x_0, r}(f(\tau) )|^2)^{ \frac{1}{2}} \varphi_{ x_0,r} \|^p_{ p}
 \\
 &=  \intl_{B(x_0,r)} v(\tau )\cdot \nabla \varphi_{ x_0, r} ( \delta +| f(\tau )- P^{ L, \delta } _{ x_0, r}(f(\tau )) |^2) ^{ \frac{p}{2}} \varphi_{ x_0, r}^{ p-1} dx  
 \\
&  \qquad + \frac{1}{p}\intl_{B(x_0,r)}\nabla \cdot  v(\tau )  ( \delta +| f(\tau )- P^{ L, \delta } _{ x_0, r} (f(\tau ))|^2) ^{ \frac{p}{2}} \varphi_{ x_0, r}^{ p} dx  
 \\
&  \qquad +\intl_{B(x_0,r)}  v(\tau )\cdot \nabla P^{ L, \delta } _{ x_0, r} (f(\tau) ) 
F_\delta  (f(\tau )- P^{ L, \delta } _{ x_0, r}(f(\tau) ))\varphi_{ x_0, r}^p dx  
 \\
&\qquad  + 
 \intl_{B(x_0, r)} (g(\tau )- P^{ N}_{ x_0, r}(g(\tau ))) F_\delta  (f(\tau )- P^{ L, \delta } _{ x_0, r}(f(\tau) )) \varphi_{ x_0, r} ^p dx.  
 \end{align*}    
 In the last line we used identity  \eqref{2.5} for $ Q=P^{ N}_{ x_0, r}(g(\tau ))$.

Multiplying both sides of the above identity by $e_\delta (\tau )^{ 1-p}$, where $e_\delta (\tau ) := 
\|(\delta +| f(\tau )- P^{ N, \delta } _{ x_0, r}(f(\tau) )|^2)^{ \frac{1}{2}} \varphi_{ { x_0, r}} \|_{ p}$, integrating the 
result over $ (0,t), t\in [0, T]$,  with respect to $ \tau $, and applying integration by parts, we find 
\begin{align*}
&  e_\delta(t) =
e_\delta (0)
  +   \intl_{0}^{t}  \intl_{B(x_0,r)} v(\tau )\cdot \nabla \varphi_{ x_0, r} ( \delta +| f(\tau )- P^{ L, \delta } _{ x_0, r}(f(\tau )) |^2) ^{ \frac{p}{2}} \varphi_{ x_0, r}^{ p-1} 
  e_\delta (\tau ) ^{ 1-p}dx   d\tau 
 \\
 &\qquad +  \frac{1}{p} \intl_{0}^{t}  \intl_{B(x_0,r)} \nabla \cdot v(\tau ) ( \delta +| f(\tau )-  P^{ L, \delta } _{ x_0, r} (f(\tau )) |^2) ^{ \frac{p}{2}} \varphi_{ x_0, r}^{ p} e_\delta (\tau )^{ 1-p} dx  d\tau 
 \\
 &\qquad +   \intl_{0}^{t}  \intl_{B(x_0,r)} 
 v(\tau )\cdot \nabla P^{ L, \delta } _{ x_0, r} (f(\tau ))  
F_\delta  (f(\tau )- P^{ L, \delta } _{ x_0, r}(f(\tau )))\varphi_{ x_0, r} ^p e_\delta (\tau )^{ 1-p} dx d\tau 
 \\
&\qquad  +  \intl_{0}^{t}
 \intl_{B(x_0, r)} (g(\tau )- P_{ x_0, r}^N(g(\tau ))) F_\delta  (f(\tau )- P^{ L, \delta } _{ x_0, r}(f(\tau ))) 
 \varphi_{ x_0, r} ^p e_\delta (\tau )^{ 1-p} dx d\tau .
 \end{align*}     
In the above identity,  letting $ \delta \rightarrow 0$ and making use of \eqref{2.6},  we obtain \eqref{2.4}.

\vspace{0.2cm}
2. Using the triangle  inequality,  we estimate  
\begin{align*}
I &\le  c\intl_{0}^{t}  \| \nabla \varphi_{ x_0, r}\cdot   v(\tau )\|_{ \infty} 
\| f(\tau )- P^{ L, \ast} _{ x_0, r}(f(\tau ))  \|_{ L^p(B(x_0, r))}
e(\tau )^{ p-1} e(\tau )^{ 1-p}d\tau  
\\
& \le  c\intl_{0}^{t}  \| \nabla \varphi_{ x_0, r}\cdot   v(\tau )\|_{ \infty} 
\| f(\tau )- P^{ L, \ast} _{ x_0, r}(f(\tau ))\|_{ L^p(B(x_0, r))}d\tau  
\\
&\le c\intl_{0}^{t}  \| \nabla \varphi_{ x_0, r}\cdot   v(\tau )\|_{ \infty} 
\| (f(\tau ) -P^{ L, \ast} _{ x_0, 2r}(f(\tau ))) \varphi _{ x_0, 2r} \|_{p}d\tau 
\\
&\qquad + c\intl_{0}^{t}  \| \nabla \varphi_{ x_0, r}\cdot   v(\tau )\|_{ \infty} 
\| P^{ L, \ast} _{ x_0, r}(f(\tau ))- P^{ L, \ast} _{ x_0, 2r}(f(\tau ))  \|_{ L^p(B(x_0, r))} d\tau  = I_1+I_2. 
\end{align*}
Thanks to the minimizing property \eqref{2.1} we get 
\[
I_1 \le c r^{ -1}\intl_{0}^{t}  \| v(\tau )\|_{ L^\infty(B(x_0, r))} 
\| f(\tau ) -P_{ x_0, 2r}^{ L}(f(\tau ))\|_{L^p(B(x_0, 2r))}d\tau . 
\]
On the other hand,  for estimating $ I_2$, making use of \eqref{A.10}, we see that for all $ \tau \in [0,T]$, 
\begin{align*}
&P^{ L, \ast} _{ x_0, r}(f(\tau ))- P^{ L, \ast} _{ x_0, 2r}(f(\tau ))
= P^{ L, \ast} _{ x_0, r}(f(\tau )- P^{ L} _{ x_0, 2r}(f(\tau ))) - 
P^{ L, \ast} _{ x_0, 2r}(f(\tau )- P^{ L} _{ x_0, 2r}f(\tau )).
\end{align*}
This,  together with \eqref{A.2b} and \eqref{A.1a}, yields  
\begin{align*}
&\| P^{ L, \ast}_{ x_0, r}(f(\tau ))- P^{ L, \ast} _{ x_0, 2r}(f(\tau )) \|_{ L^p(B(x_0, r))} 
\\
& \le  \| P^{ L, \ast}_{ x_0, r}(f(\tau )- P^{ L} _{ x_0, 2r}f(\tau ))  \|_{ L^p(B(x_0, r))}  + 
\| P^{ N, \ast} _{ x_0, 2r}(f(\tau )- P^{ L} _{ x_0, 2r}(f(\tau ))) \|_{ L^p(B(x_0, r))} 
\\
& \le c\| f(\tau ) -P_{ x_0, 2r}^{ L}(f(\tau ))\|_{L^p(B(x_0, 2r))} .  
\end{align*}
Consequently, $ I_2$ enjoys the same estimate as $ I_1$, which gives 
\[
I \le c r^{ -1}\intl_{0}^{t}  \|v(\tau )\|_{ L^\infty(B(x_0, r))} 
\| f(\tau ) -P_{ x_0, 2r}^{ L}(f(\tau ))\|_{L^p(B(x_0, 2r))}d\tau . 
\]

Using \eqref{A.1a}, we immediately get 
\begin{align*}
II &\le c \intl_{0}^{t}  \| \nabla \cdot  v(\tau )\|_{ L^\infty(B(x_0,r))} 
\| (f(\tau ) -P  ^{  L, \ast}_{ x_0, r}(f(\tau ))) \varphi_{ x_0,r} \|_{L^p(B(x_0, r))}d\tau  
\\
&\le c \intl_{0}^{t}  \| \nabla \cdot v(\tau )\|_{ L^\infty(B(x_0,r))} 
\| f(\tau ) -P^{ L}_{ x_0, 2r}(f(\tau ))\|_{L^p(B(x_0, 2r))}d\tau . 
\end{align*}

We proceed with the estimation of $ III$.  Clearly, in case $ N=0$, since $ P^{ L, \ast} _{ x_0, r}(f(\tau )) =\const$ for all $ \tau \in [0,T]$,  the integral  $ III$ vanishes. Thus, it only remains the case   
$ N >0$. Let $ \tau \in [0, T]$ be fixed. Making use of \eqref{2.5} with $ \delta =0$,  we find 
\begin{align*}
&\intl_{B(x_0,r)} 
 v(\tau )\cdot \nabla P^{ L, \ast} _{ x_0, r}(f(\tau )) \cdot   
| f(\tau )- P^{ L, \ast} _{ x_0, r}(f(\tau ))|^{ p-2}(f(\tau )- P^{ L, \ast}_{ x_0, r}(f(\tau )))\varphi_{ x_0, r} ^pdx  
\\
&=\intl_{B(x_0,r)} 
 v(\tau )\cdot \nabla (P^{ L, \ast} _{ x_0, r}(f(\tau )) -P^{ N} _{ x_0, r}(f(\tau ))) \cdot   
F_0 \Big(f(\tau )- P^{ L, \ast}_{ x_0, r}(f(\tau ))\Big)\varphi_{ x_0, r}^pdx  
\\
&\quad  +\intl_{B(x_0,r)} 
v(\tau )\cdot \nabla P^{  N} _{ x_0, r}(f(\tau )) \cdot   
F_0 \Big(f(\tau )- P^{ L, \ast}_{ x_0, r}(f(\tau ))\Big)\varphi_{ x_0, r} ^pdx  
\\
& = 
\intl_{B(x_0,r)} 
 v(\tau )\cdot \nabla (P^{ L, \ast} _{ x_0, r}(f(\tau )) -P^{  N} _{ x_0, r}(f(\tau ))) \cdot   
F_0 \Big(f(\tau )- P^{ L, \ast}_{ x_0, r}(f(\tau ))\Big)\varphi_{ x_0, r}^pdx  
\\
&\quad   +\intl_{B(x_0,r)} 
  (v(\tau ) - P_{ x_0, r}^N(v(\tau )))\cdot \nabla P^{  N} _{ x_0, r}(f(\tau )) \cdot   
F_0 \Big(f(\tau )- P^{ L, \ast}_{ x_0, r}(f(\tau ))\Big)\varphi_{ x_0, r}^pdx  
\\
& =J_1+ J_2.   
\end{align*}
Using the fact that $ P^{ L, \ast}_{ x_0, r}(Q) = P_{ x_0, r}^{ N}(Q)=Q$ for all $ Q\in \mathcal{P}_{ N}$, we get with $ Q= 
P^{ N}_{ x_0, r}(f(\tau ))$ 
for all $ \tau \in (0,t)$
\[
\| \nabla (P^{ L, \ast} _{ x_0, r}(f(\tau )) - P^{ N} _{ x_0, r}(f(\tau ))\|_{ L^p(B(x_0, r))} \le c r^{ -1}
\| f(\tau )- P^{ N}_{ x_0, r}(f(\tau ))\|_{ L^p(B(x_0, r))}. 
\]
Then H\"older's inequality 
yields
\begin{align*}
J_1 &\le c r^{ -1}\| v(\tau )\|_{ L^\infty(B(x_0, r))} \| f(\tau )- P^{ N}_{ x_0, r}(f(\tau ))\|_{ L^p(B(x_0, r))} e(\tau )^{ p-1}.
\end{align*}
Similarly, 
\begin{align*}
J_2 &\le c \| v(\tau )- P_{ x_0, r}^N(v(\tau ))\|_{ L^p(B(x_0, r))} 
\| \nabla P^{ N}_{ x_0, r}(f(\tau ))\|_{ L^\infty(B(x_0, r))}
e(\tau )^{ p-1}.
\end{align*}

Inserting the estimates of $ J_1$ and $ J_2$ into the integral of $ III$, we obtain 
\begin{align*}
III & \le  c  r^{ -1}  \intl_{0}^{t}   \|  v(\tau )\|_{L^\infty(B(x_0, r))} \| f(\tau )- P^{ N}_{ x_0, 2r}(f(\tau ))\|_{ L^p(B(x_0, 2r))} d\tau 
\\
& \qquad +  c  \intl_{0}^{t}  \| v(\tau )- P_{ x_0, r}^N(v(\tau ))\|_{ L^p(B(x_0, r))} 
\| \nabla P^{ N}_{ x_0, r}(f(\tau ))\|_{ L^\infty(B(x_0, r))} d\tau. 
  \end{align*}
To estimate $ IV$, we use H\"older's inequality. This leads to 
\[
IV \le  \intl_{0}^{t} \| g(\tau )- P^N_{ x_0, r}(g(\tau ))\|_{ L^p(B(x_0, r))}   d\tau. 
\] 
Inserting the estimates of $ I,II, III$ and $ IV$ into the right-hand side of \eqref{2.3}, we find 
\begin{align}
e(t) &\le   e(0) + c r^{ -1} \intl_{0}^{t} \|   v(\tau )\|_{ L^\infty(B(x_0, r))}  
\| f(\tau ) -P_{ x_0, 2r}^{ N}(f(\tau ))\|_{L^p(B(x_0, 2r))}d\tau  
\cr
&\quad  +c \intl_{0}^{t}   \|  \nabla \cdot  v(\tau )\|_{ L^\infty(B(x_0, r))}\| f(\tau ) -P_{ x_0, 2r}^{ N}(f(\tau ))\|_{L^p(B(x_0, 2r))}d\tau  
\cr
& \quad + c  \intl_{0}^{t}  \| v(\tau )- P_{ x_0, r}^N(v(\tau ))\|_{ L^p(B(x_0, r))} 
\| \nabla P^{ N}_{ x_0, r}(f(\tau ))\|_{ L^\infty(B(x_0, r))}d\tau 
\cr
&\quad +  c\intl_{0}^{t} \| g(\tau )- P^{ N}_{ x_0, r}(g (\tau ))\|_{ L^p(B(x_0, r))}   d\tau. 
  \label{2.8}
 \end{align}       
Noting that 
\[
\| f(t)- P_{ x_0,  \frac{r}{2}}^{ L}(f(t)) \|_{L^p(B(x_0, \frac{r}{2})) }  \le c \| (f(t)- P^{L, \ast} _{ x_0, r}(f(t))) \varphi_{ x_0, r} \|_{ p} = c e(t), 
\]
and using \eqref{A.1a},  recalling that $ L=2N-1$, the inequality \eqref{2.4} follows from   \eqref{2.8}.  \hfill \Beweisende

 \begin{rem}
 \label{rem2.2} Given
 $v\in  L^1(0, T; C^{ 0,1}(\R^{n}; \R^{n} ))$,  and 
$ \pi  \in L^1(0, T; W^{1,\,2}_{ loc}(\R^{n}; \R^{n}  ))$, 
let $ f\in L^\infty(0, T; C^{ 0,1}(\R^{n}; \R^{n} ))$ with $ \nabla \cdot f=0$ be a weak solution to the system 
 \begin{equation}
\partial _t f + (v\cdot \nabla) f = -\nabla \pi \quad  \text{ in}\quad  Q_T. 
\label{2.2s}
\end{equation}

Then, repeating the proof of Lemma\,\ref{lem2.1} for the case $ p=2$ and $ N=1$ in the vector valued case, we find   
\begin{align}
& e(t)
= e(0)
 +  \intl_{0}^{t}  v\cdot \nabla \varphi _{ x_0, r}  
| f- P^{ 1,\ast}_{ x_0, r}(f)|^{ 2}\varphi _{ x_0, r} e(\tau )^{ -1}
 dx  d\tau 
 \cr
 &\qquad +  \frac{1}{2}\intl_{0}^{t}  \intl_{B(x_0,r)} \nabla \cdot v  |f- P^{ 1,\ast}_{ x_0, r}(f)|^2\varphi_{ x_0, r}^{ 2}  
 e(\tau )^{-1} dxd\tau 
 \cr
 &\qquad +   \intl_{0}^{t}  \intl_{B(x_0,r)} 
 v\cdot \nabla P^{ 1,\ast}_{ x_0, r}(f) \cdot   
(f- P^{ 1,\ast}_{ x_0, r}(f))\varphi_{ x_0,r} ^2  
e(\tau )^{-1} dxd\tau 
\cr
&\qquad  +  \intl_{0}^{t}
 \intl_{B(x_0, r)} (\nabla \pi -P_{ x_0,r}^{ 1}(\nabla \pi ))   (f- P^{ 1,\ast}_{ x_0, r}(f))\varphi_{ x_0, r} ^2
e(\tau )^{-1}dx d\tau
\cr
&= e(0) + I+II+III+IV, 
 \label{2.3a}
  \end{align}
  where 
 \[
e(\tau )=  \|(f(\tau )- P^{ 1,\ast}_{ x_0, r}(f(\tau ))) \varphi_{ x_0, r}  \|_{ 2},\quad  \tau \in [0, T]. 
\]
The integrals $ I, II$ and $ III$ can be estimated as in the proof of Lemma\,\ref{lem2.1}.  For the estimation 
of $ IV$ we proceed as follows. 

Assume that the mollifier $ \varphi \in C^\infty_c(B(1))$ is  radial symmetric. 
Let $ u\in L^1(B(x_0, r))$. It can be checked easily that the 
minimal polynomial $ P^{ 1, \ast}_{ x_0, r}(u)$ is given by 
\[
\hspace*{-1.5cm} P^{ 1, \ast}_{ x_0, r}(u)(x) = \frac{1}{ \intl_{ \R^{n} } \varphi _{ x_0, r}^2 dy} 
\intl_{ \R^{n} } u\varphi _{ x_0, r}^2 dy + \frac{n}{ \intl_{ \R^{n} } \varphi _{ x_0, r}^2 |x_0 - y|^2 dy} 
\intl_{ \R^{n} } u\varphi _{ x_0, r}^2 (y_i- x_{ 0, i}) dy (x_i- x_{ 0, i}).
\] 
In case $ u=(u_1, \ldots, u_n)$ with $ \nabla\cdot u= 0 $ almost everywhere in $ B(x_0, r)$, recalling that 
$ \varphi$ is radialsymmetric, by Gauss' theorem  we get 
\[
\nabla \cdot P^{ 1, \ast}_{ x_0, r}(u)(x)=  \frac{n}{ \intl_{ \R^{n} } \varphi _{ x_0, r}^2 |x_0 - y|^2 dy} \intl_{ B(x_0, r) } u\cdot (y- x_0) \varphi _{ x_0, r}^2 dy=0. 
\]   
Using integration by parts together with $ \nabla \cdot P^{ 1, \ast}_{ x_0, r}(f(\tau ))=0$,  and applying Sobolev-Poincar\'e inequality, we get
\begin{align*}
& \intl_{B(x_0, r)} (\nabla \pi(\tau )-P_{ x_0,r}^{ 1}(\nabla \pi(\tau ) )) 
  (f(\tau )- P^{ 1,\ast}_{ x_0, r}(f(\tau )))\varphi_{ x_0, r} ^2
e(\tau )^{-1}dx 
\\
&= - 2\intl_{B(x_0, r)} (\pi(\tau ) -P_{ x_0,r}^{ 2}( \pi(\tau ) )) 
  (f(\tau )- P^{1, \ast}_{ x_0, r}(f(\tau )))  \varphi_{ x_0, r}\cdot \nabla \varphi_{ x_0, r}
e(\tau )^{-1}dx  
\\
& \le cr^{ -1} \bigg( \intl_{B(x_0, r)} |\nabla \pi(\tau ) -P_{ x_0,r}^{ 1}(\nabla \pi(\tau ) ))|^{ \frac{2n}{n+2}}dx \bigg)^{ \frac{n+2}{2n}} \le c r^{ \frac{n}{2}} \osc_{ \frac{2n}{n+2}, 1} (\nabla \pi(\tau ) ; x_0, r). 
\end{align*}
This  yields
\[
IV \le c r^{ \frac{n}{2}} \intl_{0}^t \osc_{ \frac{2n}{n+2}, 1} (\nabla \pi(\tau ) ; x_0, r) d\tau. 
\]
Inserting the estimates of $ I, II, III$ and $ IV$ into the right-hand side of  \eqref{2.3a},  and arguing as in the proof of Lemma\,\ref{lem2.1}, we arrive at 
 \begin{align}
& \osc_{ 2, 1} \Big(f(t); x_0,  \frac{r}{2}\Big)
\le   c \osc_{2, 1} (f(0); x_0,  r) +c r^{ -1} \intl_{0}^{t} \|   v(\tau )\|_{ L^\infty(B(x_0, r))}  
\osc_{ 2, 1} (f(\tau ); x_0,  2r)d\tau  
 \cr
&\quad  +c \intl_{0}^{t}   \|  \nabla \cdot  v(\tau )\|_{ L^\infty(B(x_0, r))} \osc_{2, 1} (f(\tau ); x_0,  2r)d\tau  
\cr
& \quad +  c  \intl_{0}^{t}  \osc_{ 2, 1} (v(\tau ); x_0,  r) 
 |\nabla P^{ 1}_{ x_0, r}(f(\tau )| d\tau 
\cr
&\quad +  c\intl_{0}^{t} \osc_{ \frac{2n}{n+2}, 1} (\nabla \pi (\tau ); x_0,  r)  d\tau.
  \label{2.4a}
 \end{align}

\end{rem}

%%%%%%%%%%%
%%%%%%%%%%%

\vspace{0.5cm}
{\bf Proof of the main theorems }

\vspace{0.3cm}
 {\it 1. Existence and uniqueness in terms of particle trajectories.} Assume $ f_0\in \mathscr{L}^s_{ q (p,N)}(\R^{n} ), 
 g\in L^1(0,T; \mathscr{L}^s_{ q (p,N)}(\R^{n} ))$,  and $ \nabla v\in L^1(0,T; L^\infty(\R^{n} ))$.  Let $(x,t)\in Q_T$ be fixed. By $  X_{ t}(x, \cdot  )$ we denote the unique solution to the ODE
\begin{equation}
\frac{d}{d \tau } X_{  t}(x, \tau ) = v(X_{t}(x, \tau ), \tau ),\quad  \tau \in [0,T],\quad   X_{t}(x, t)=x,
\label{C.4}
\end{equation}
which is ensured by Carath\'eodory's theorem. 
We define the flow map $ \Phi _{ t, \tau }: \R^{n} \rightarrow  \R^{n} $ by means of 
\[
\Phi _{ t, \tau } (x ) = X_{t}(x, \tau ),\quad  x\in \R^{n}, \quad \tau , t\in [0,T]. 
\]
By the uniqueness of this flow  we get the inverse formula 
\[
\Phi ^{ -1}_{ t, \tau }(x) = \Phi _{ \tau , t} (x). 
\]
Furthermore, from \eqref{C.4} we deduce that 
\begin{equation}
\frac{d}{d \tau } \Phi _{t, \tau }(x ) = v(\Phi _{t, \tau }(x ), \tau ),\quad   \tau \in [0,T],\quad   \Phi _{t, t }(x )=x.
\label{C.5}
\end{equation}

Let $ (x, t)\in Q_T$. We set $ y = \Phi _{ t, 0}(x)$, which is equivalent to $x = \Phi _{0, t}(y) $. We define $ f$ by means of  
\begin{equation}
f(x, t) = f_0 (y) +  \intl_{0}^{t} g(\Phi _{ 0, s}(y), s)  ds. 
\label{C.6}
\end{equation}
Recalling that $ f(t)$ is Lipschitz  for almost all $ t\in (0,T)$, we see that $ f$ is differentiable with respect to time almost everywhere in $ (0,T)$. 
Recalling the inverse formula, it holds $ x= \Phi _{0, t}(y)$. Consequently, for $ y\in \R^{n}$ fixed we get from \eqref{C.6}
\begin{equation}
f(\Phi _{0, t}(y), t) = f_0 (y) +  \intl_{0}^{t} g(\Phi _{0, s}(y), s)  ds\quad \forall\,t\in (0,T). 
\label{C.7}
\end{equation} 
Differentiating \eqref{C.7} with respect to $ t$, and observing \eqref{C.5},  we obtain 
\begin{equation}
\partial _t f(\Phi _{0, t}(y), t)  + (v(\Phi _{0, t}(y), t)\cdot \nabla) f(\Phi _{0, t}(y), t) =  g(\Phi _{0, t}(y), t). 
\label{C.8}
\end{equation}
This shows that $ f$ solves \eqref{trans} in $ Q_T$.  In addition, verifying that 
$ \Phi _{ 0,0}(x)=x$, we get from \eqref{C.7}
\[
f(x,0) = f_0(x) \quad \forall\,x\in \R^{n}. 
\]   

This solution is also unique. In fact,  assume there is another solution $ \overline{f} $ solves \eqref{trans}.  Setting 
$ w= f-\overline{f} $, 
then $ w$ solves \eqref{trans} with homogenous data. In other words for every $ y\in \R^{n}$ the function $Y(t)= 
w(\Phi _{ 0, t}(y),t)$ solves the 
ODE
\[
\dot{Y} = 0, \quad Y(0)=0,
\] 
which implies $ Y \equiv 0$, and thus $ w(\Phi _{ 0, t}(y),t) =0$. With  $ y= \Phi _{ t, 0}(x)$ we get 
$ w(x, t)=0$ for all $ (x,t )\in Q_T$.

\vspace{0.3cm}
{\it 2. Growth of the solution as $ |x| \rightarrow +\infty$.}
Applying $ \nabla_x $ to both sides of \eqref{C.5},  and using the chain rule, we find that 
\begin{equation}
\frac{d}{d \tau } \nabla \Phi _{s, \tau }(x ) = \nabla v(\Phi _{s, \tau }(x ), \tau )\cdot \nabla \Phi _{s, \tau }(x ).
\label{C.9}
\end{equation}
Integration  with respect to $ \tau $  over $ (s, t)$  yields 
\[
\nabla \Phi _{s, t}(x ) = I +  \intl_{s}^{t}  \nabla v(\Phi _{s, \tau }(x ), \tau )\cdot \nabla \Phi _{s, \tau }(x ) d\tau,
\] 
 where $ I$ stands for the unit  matrix. Thus, for all $ s, t\in (0,T)$,
 \[
| \nabla \Phi _{ s,t}(x)| \le 1 + \intl_{s}^{t} \| \nabla v(\tau )\|_{ \infty}  | \nabla \Phi _{s, \tau }(x )| d\tau. 
\]
By means of Gronwall's lemma it follows that for all $ s,t\in (0,T)$
\begin{equation}
| \nabla \Phi _{ s,t}(x)| \le \exp  \bigg(\intl_{s}^{t} \| \nabla v(\tau )\|_{ \infty}    d\tau \bigg). 
\label{C.10}
\end{equation}

 From the definition \eqref{C.4} we deduce that 
 \begin{align*}
 \nabla f(x,t) &= \nabla f_0  (\Phi _{ t, 0}(x)) \cdot \nabla \Phi _{ t, 0}(x)   +  \intl_{0}^{t}  \nabla _x g(\Phi _{ 0, \tau }(\Phi _{ t, 0}(x)), \tau )  d\tau
 \\
 &= \nabla f_0  (\Phi _{ t, 0}(x)) \cdot \nabla \Phi _{ t, 0}(x)   
 \\
 &\qquad +  \intl_{0}^{t}  \nabla g(\Phi _{ 0, \tau }(\Phi _{ t, 0}(x)), \tau ) 
 \cdot \nabla \Phi _{ 0, \tau }(\Phi _{ t, 0}(x) ) \cdot \nabla \Phi _{ t, 0}(x)  d\tau.
 \end{align*} 
 Thus, in case $ \nabla f_0\in L^\infty(\R^{n} )$ and $ g\in L^1(0,T; L^\infty(\R^{n} ))$, 
  in view of \eqref{C.10} we get for all $ t\in (0,T)$ 
 \begin{equation}
 \|\nabla f(t) \| _{ \infty}\le  \bigg(\| \nabla f_0 \|_{ \infty} +   \intl_{0}^{T} \| \nabla g(\tau )\|_{ \infty}\bigg) \exp  \bigg(2\intl_{0}^{T} \| \nabla v(\tau )\|_{ \infty}    d\tau \bigg). 
 \label{C.11}
 \end{equation}

Using integration by parts,  from  \eqref{C.5} we get    for all $ s,t\in (0,T)$
\[
\Phi _{ s,t} (x) - x = \Phi _{ s,t}- \Phi _{s,s} (x)  = 
 \intl_{s}^{t} v(\Phi _{ s, \tau }(y), \tau )- v(0, \tau )  d\tau 
+\intl_{s}^{t} v(0, \tau )  d\tau.
\] 
This leads to the inequality 
\[
| \Phi _{ s, t}(x)| \le | x|+ \intl_{0}^{T} | v(0, \tau )  |d\tau + \intl_{s}^{t} 
\| \nabla v(\tau )\|_{ \infty} | \Phi _{ s, \tau }(y)| d\tau.
\]
By means of Gronwall's lemma we find for all $ s,t\in (0,T)$
\begin{equation}
| \Phi _{ s, t}(x)|  \le \bigg(| x|+ \intl_{0}^{T} | v(0, \tau )|  d\tau\bigg) 
\exp \bigg( \intl_{0}^{T} \| \nabla v(\tau )\|_{ \infty}  d\tau\bigg)
\le c(1+| x|).  
\label{C.12}
\end{equation}

 Let $ x \in \R^{n}$ and $ t\in (0,T)$.  In case $ N=0, s\in [0,1)$, using Lemma\,\ref{growth}, we get 
 \begin{align}
 |f_0(x)| &\le c (1+ |x|^s) \|f_{ 0}\|_{ \mathscr{L}^s_{ q (p, N)}},
  \label{C.12a1}
 \\
 |g(x, \tau )|  &\le  c (1+ |x|^s) \|g(\tau )\|_{ \mathscr{L}^s_{ q (p, N)}}. 
   \label{C.12a2}
 \end{align}
 In case $ N=1, s=1$ and $ 1 < q \le \infty$ we get by Lemma\,\ref{growth}
 \begin{align}
 |f_0(x)| &\le c (1+ \log(1+ |x|)^{ \frac{1}{q'}}|x|) \|f_{ 0}\|_{ \mathscr{L}^s_{ q (p, N)}},
   \label{C.12b1}
 \\
 |g(x, \tau )|  &\le  c (1+ \log(1+ |x|)^{ \frac{1}{q'}}|x|) \|g(\tau )\|_{ \mathscr{L}^s_{ q (p, N)}}
   \label{C.12b2}
 \end{align}
 with $q'=\frac{q}{q-1}$.
 In the remaining cases having $\nabla  f_0 \in L^\infty(\R^{n} )$ and $ \nabla g \in L^1(0,T; L^\infty(\R^{n} ))$, 
 we find, 
 \begin{align}
 |f_0(x)| &\le c (1+ |x|) (\|f_{ 0}\|_{ \mathscr{L}^s_{ q (p, N)}}+ \|\nabla f_0\|_{ \infty}),
   \label{C.12c1}
 \\
 |g(x, \tau )|  &\le  c (1+ |x|) (\|g(\tau )\|_{ \mathscr{L}^s_{ q (p, N)}}+ \|\nabla g(\tau )\|_{ \infty}). 
   \label{C.12c2}
 \end{align}

Setting  $ y= \Phi _{t, 0}(x)$,  
we get  from  \eqref{C.7}   
\begin{align*}
| f(x, t)| & \le  | f_0(y)| + \intl_{0}^{t} | g(\Phi _{ 0, s}(y), s) | ds 
\end{align*}
Employing  \eqref{C.12a1}-  \eqref{C.12c2} together with \eqref{C.12}, we see that 
for all $ (x,t)\in Q_T$
\begin{equation}
| f(x, t)| \le c\begin{cases}
  (1+  | x|^{ \min\{s, 1\}})\quad  &\text{if}\quad s \neq 1,
\\[0.3cm]
  (1+ \log(1+|x|)^{ \frac{1}{q'}} | x|)\quad  &\text{if}\quad s=1,
\end{cases}
 \label{C.13}
\end{equation}
where 
$ c$ stands for a constant depending on $ s,q, p, N, n$ and $f_0, g$ and  $v $.

\vspace{0.3cm}
{\it 3. Local energy estimation}.  Let $ x_0 \in \R^{n}$.  Let  $ \xi \in C^2([0,T]; \R^{n})$ be a solution to the ODE
\begin{equation}
\dot{\xi}(\tau ) = v(x_0+ \xi(\tau ), \tau )\quad \tau \in [0, T].  
\label{C.14}
\end{equation}
We set 
\begin{align*}
&F(x, \tau ) = f(x+ \xi(\tau ), \tau ),\quad V(x, \tau ) = v(x+ \xi(\tau ) , \tau )- \dot{\xi}(\tau ),
\\
&\qquad G(x, s) = g(x+ \xi(\tau ) , \tau ), \quad  (x, s)\in Q_{ T}. 
\end{align*}
It is readily seen that $ V$ solves the transport equation
\begin{equation}
\partial _t F + (V\cdot \nabla) F = G  \quad  \text{ in}\quad  Q_T. 
\label{C.15}
\end{equation}
In particular, from \eqref{C.14} we infer 
\begin{equation}
V(x_0, \tau )=0\quad  \forall\,\tau \in [0, T]. 
\label{C.16}
\end{equation}  
Set $ L=2N-1$ if $ N>0$ and $ L=0$ if $ N=0$. According to \eqref{2.4} of  Lemma\,\ref{lem2.1} with $ r= 2^{ j+1}, j\in \Z$, noting that  in view of  \eqref{C.16} it holds $2^{ -j} \|   V(\tau )\|_{ L^\infty(B(x_0, 2^{ j+1}))} \le c\|\nabla v(\tau )\|_{ \infty}$, we find 
 \begin{align}
& \osc_{ p, L} (F(t); x_0,  2^j)
\le   c \osc_{ p, L} (f_0(\cdot + \xi (0)); x_0,  2^{ j+1}) 
\cr
&\quad +c  \intl_{0}^{t} \| \nabla v(\tau ) \|_{\infty}  
\osc_{ p, N} (F(\tau ); x_0,  2^{ j+2})d\tau  
\cr
& \quad + \delta_{ N0} c  \intl_{0}^{t}  \osc_{ p, N} (V(\tau ); x_0,  2^{ j+1}) 
\| \nabla P^{ N}_{ x_0, r}(F(\tau ))\|_{ L^\infty(B(x_0, 2^{ j+1}))} d\tau 
\cr
&\quad +  c\intl_{0}^{t} \osc_{ p, N} (G(\tau ); x_0,  2^{ j+1})  d\tau,
  \label{C.17}
 \end{align}       
where $ \delta _{ N0}=0$ if $ N=0$ and $ 1$ otherwise.

\vspace{0.3cm}
{\bf Proof of  \eqref{estimate1} in  Theorem\,\ref{thm1}} 
Inequality  \eqref{C.17}   gives 
\begin{align}
& \osc_{ p, 0} (F(t); x_0,  2^j)
\cr
&\quad \le   c \osc_{ p, 0} (F (0); x_0,  2^{ j+1}) +c \intl_{0}^{t} \|   \nabla v(\tau )\|_{ \infty}  
\osc_{ p, 0} (F(\tau ); x_0,  2^{ j+2})d\tau  
\cr
 &\qquad +  c\intl_{0}^{t} \osc_{ p, 0} (G(\tau ); x_0,  2^{ j+1})  d\tau.
  \label{C.17a}
 \end{align}

Observing \eqref{C.13},  since $ s <1$,  we get  $ S_{ 1,1}( \osc_{p, 0}(f(\tau ); x_0)) < +\infty$.  
Thus, applying $S_{ 1,1} $ 
 to both sides of \eqref{C.17a}, we obtain
 \begin{align}
& S_{ 1,1}(\osc_{p, 0}(F(t); x_0))
\cr
&\le   c  S_{ 1,1}(\osc_{p, 0}(F (0); x_0))
 + c\intl_{0}^{t}  \| \nabla v(\tau )\|_{ \infty} 
S_{ 1,1}(\osc_{p, 0}(F(\tau ); x_0))d\tau 
\cr
&\qquad +  c\intl_{0}^{t} S_{ 1,1}(\osc_{p, 0}(G(\tau ); x_0) )  d\tau . 
  \label{C.18}
 \end{align}    
 Applying Gronwall's lemma, we deduce from \eqref{C.18}  
\begin{align}
&\osc_{p, 0}(F(t); x_0)
\cr
& \le  S_{ 1, 1} ( \osc_{p, 0}(F(t); x_0)) 
\cr
&\le   c  \bigg\{S_{ 1, 1}(\osc_{p, 0}(F (0); x_0) 
 +  \intl_{0}^{t}  S_{ 1, 1}(\osc_{p, 0}(G(\tau ); x_0) )  d\tau\bigg\} 
  \exp \bigg(c\intl_{0}^{t} \| \nabla v(\tau )\|_{ \infty} d\tau\bigg). 
  \label{C.21}
 \end{align}     
Let $ t\in [0,T]$. Clearly, the constant in \eqref{C.21} is independent of  the choice of the characteristic for $\xi $. Therefore, we may choose 
 $ \xi $ such that $ \xi (t)=0$, which implies $ F(t)=f(t)$. Hence, we may replace $ F(t)$ by $ f(t)$ 
 on the left-hand side of \eqref{C.21}. Afterwards,   with the help of Lemma\,\ref{lem10.1} we are in a position to  operate $ S_{s, q}$ to both sides of \eqref{C.21}, verifying $ F(0)= f_0(\cdot -\xi (0))$,  that  yields 
\begin{align}
&(S_{ s,q}(\osc_{p, 0}(f(t); x_0)))_j
\cr
&\le   c\bigg\{(S_{ s,q}(\osc_{p, 0}(f_0(\cdot  -\xi (0); x_0))))_j
\cr
 &\qquad +  \intl_{0}^{t}   (S_{ s,q}(\osc_{p, 0}(G(\tau ); x_0))))_jd\tau\bigg\}  \exp \bigg(c\intl_{0}^{t} \| \nabla v(\tau )\|_{ \infty} d\tau\bigg). 
    \label{C.22}
 \end{align}    
 Multiplying both sides by $ 2^{ -js}$, we get 
\begin{align}
&  \Big(\sum_{i=j}^{\infty} (2^{ -si }\osc_{p, 1}(f(t); x_0; 2^{ i}))^q\Big)^{ \frac{1}{q}}
\cr
&\le   c\bigg\{|f_0|_{ \mathscr{L} ^{s}_{q (p,0) }} +\intl_{0}^{t}   |G(\tau )|_{ \mathscr{L} ^{s}_{q (p,0) }}d\tau\bigg\}
  \exp \bigg(c\intl_{0}^{t} \| \nabla v(\tau )\|_{ \infty} d\tau\bigg). 
    \label{C.22a}
 \end{align}    
Passing $ j \rightarrow -\infty$ and taking the supremum   over  $ x_0\in \R^{n}$ in \eqref{C.22a}, we get
\eqref{estimate1}.     \hfill \Beweisende 
\vspace{0.2cm}
{\bf Proof of  \eqref{estimate2} in Theorem\,\ref{thm2}.}  Recalling that $ V(x_0, \tau )= 0$ for all $ \tau \in [0,T]$, we see that 
$  2^{ -j} \|   V(\tau )\|_{ L^\infty(B(x_0, 2^{ j+1}))} \le c\|   \nabla v(\tau )\|_{\infty}$ 
and $  2^{ -j} {\D \osc_{ p, 0} (V(\tau ); x_0,  2^{ j+1})} \le c\|   \nabla v(\tau )\|_{ \infty}$. Thus, \eqref{C.7} 
leads to 
\begin{align}
& \osc_{ p, 1} (F(t); x_0,  2^j)
\cr
&\quad \le   c \osc_{ p, 1} (F(0); x_0,  2^{ j+1}) +c \intl_{0}^{t} \|   \nabla v(\tau )\|_{ \infty}  
\osc_{ p, 1} (F(\tau ); x_0,  2^{ j+2})d\tau  
\cr
&\qquad  +c  \intl_{0}^{t}  \osc_{ p, 1} (V(\tau ); x_0,  2^{ j+1}) |\nabla \dot{P}^1_{ x_0, 2^{ j+1}} (F(\tau ))|d\tau 
\cr
 &\qquad +  c\intl_{0}^{t} \osc_{ p, 1} (G(\tau ); x_0,  2^{ j+1})  d\tau.
  \label{C.17a}
 \end{align}     
In case $ j \ge 0$,  using  triangle inequality,   we get 
\begin{align*}
|\nabla \dot{P}^1_{ x_0, 2^j}(F(\tau ))|&\le  c \sum_{i=0}^{j} 2^{ -i} \osc_{p, 1}(F(\tau ); x_0, 2^i) + 
|\nabla \dot{P}^1_{ x_0, 1}(F(\tau ))|
\\
&\le c 2^{ -j}(S_{ 3,1}  (\osc_{p, 1}(F(\tau ); x_0)))_{ j} + |\nabla \dot{P}^1_{ x_0, 1}(F(\tau ))|.
\end{align*}
 In case $ j < 0$, using triangle inequality along with H\"older's inequality, we find 
 \begin{align*}
|\nabla \dot{P}^1_{ x_0, 2^j}(F(\tau ))|&\le  c \sum_{i=0}^{j} 2^{ -i} \osc_{p, 1}(F(\tau ); x_0, 2^i) + 
|\nabla \dot{P}^1_{ x_0, 1}(F(\tau ))|
\\
&\le 
(-j)^{ \frac{1}{q'}}\Big(\sum_{i=j}^{0}  2^{ -iq}(\osc_{p, 1}(F(\tau ); x_0, 2^i))^q\Big)^{ \frac{1}{q}} + |\nabla \dot{P}^1_{ x_0, 1}(F(\tau ))|.
\end{align*} 
 Summing up the above estimates, we arrive at 
\begin{align}
& \osc_{ p, 1} (V(\tau ); x_0,  2^{ j+1}) |\nabla \dot{P}^1_{ x_0, 2^{ j+1}} (F(\tau ))|
\cr
&\le 2^{ -j} \osc_{ p, 1} (V(\tau ); x_0,  2^{ j+1})(S_{ 3,1}  (\osc_{p, 1}(F(\tau ); x_0)))_{ j} 
\cr
&\quad + 
  c( j^{ -})^{ \frac{1}{q'}}  \osc_{p, 1}(V(\tau ); x_0, 2^{ j+1}) \Big\{\Big(\sum_{i=-\infty}^{0}  2^{ -iq}(\osc_{p, 1}(F(\tau ); x_0, 2^i))^q\Big)^{ \frac{1}{q}}+ 
|\nabla \dot{P}^1_{ x_0, 1}(F(\tau ))|\Big\},
 \label{C.29a}
\end{align}
where $ j ^- = -\min\{j,0\}$. 
Applying the operator $ S_{ 2,1}$ to the both sides of the above inequality, and making use of Lemma\,\ref{10.1}, 
with $ p=q=1$, $ \alpha =3$ and $ \beta =2$,  we obtain 
\begin{align}
&S_{ 2,1} \Big(\Big\{\osc_{ p, 1} (V(\tau ); x_0,  2^{ i+1}) |\nabla \dot{P}^1_{ x_0, 2^{ i+1}} (F(\tau ))|\Big\}\Big)
\cr
&\le c|v(\tau )|_{ \mathscr{L}^1_{ q(p,1)}}S_{ 2,1}  (\osc_{p, 1}(F(\tau ); x_0))
\cr
&\quad + 
 cS_{ 2,1}\Big(\Big\{ ( i^{ -})^{ \frac{1}{q'}}  \osc_{p, 1}(V(\tau ); x_0, 2^i) \Big\}\Big)
\Big\{\Big(\sum_{i=-\infty}^{0}  2^{ -iq}(\osc_{p, 1}(F(\tau ); x_0, 2^i))^q\Big)^{ \frac{1}{q}}+ 
|\nabla \dot{P}^1_{ x_0, 1}(F(\tau ))|\Big\}. 
 \label{C.29b}
\end{align}
Observing  \eqref{C.13}, all sum in the above estimates are finite. Again appealing to  \eqref{C.12} we are in a position to apply $ S_{ 2,1}$ to both sides of  \eqref{C.17a} to get  
 \begin{align}
& \hspace*{-1cm}S_{ 2,1}(\osc_{p, 1}(F(t); x_0))
\cr
& \hspace*{-1cm}\le   c  S_{2,1}(\osc_{p, 1}(F(0); x_0))
 + c\intl_{0}^{t} ( \| \nabla v(\tau )\|_{ \infty} +|v(\tau )|_{ \mathscr{L}^1_{ q(p,1)}})
S_{2,1}(\osc_{p, 1}(F(\tau ); x_0))d\tau 
\cr
& \hspace*{-1cm}  +c\intl_{0}^{t} S_{ 2,1}\Big(\Big\{ ( i^{ -})^{ \frac{1}{q'}}  \osc_{p, 1}(V(\tau ); x_0, 2^i) \Big\}\Big)
 \Big\{\Big(\sum_{i=-\infty}^{0}  2^{ -iq}(\osc_{p, 1}(F(\tau ); x_0, 2^i))^q\Big)^{ \frac{1}{q}}+ 
|\nabla \dot{P}^1_{ x_0, 1}(F(\tau ))|\Big\}d\tau 
\cr
&\, +  c\intl_{0}^{t} S_{ 2,1}(\osc_{p, 1}(G(\tau ); x_0) )  d\tau . 
  \label{C.18}
 \end{align}    
 Applying Gronwall's lemma, we are led to 
 \begin{align}
 & \hspace*{-1cm}\osc_{p, 1}(F(t); x_0)
 \cr
&\hspace*{-1cm} \le  S_{ 2,1}(\osc_{p, 1}(F(t); x_0))
\cr
&\hspace*{-1cm}\le   \bigg\{c  S_{2,1}(\osc_{p, 1}(f_0(\cdot + \xi (0); x_0))
\cr
&\hspace*{-1cm}+c\intl_{0}^{t} S_{ 2,1}\Big(\Big\{ ( i^{ -})^{ \frac{1}{q'}}  \osc_{p, 1}(V(\tau ); x_0, 2^i) \Big\}\Big)
 \Big(\sum_{i=-\infty}^{0}  2^{ -iq}(\osc_{p, 1}(F(\tau ); x_0, 2^i))^q\Big)^{ \frac{1}{q}}+ 
c|\nabla \dot{P}^1_{ x_0, 1}(F(\tau ))|\Big)d\tau 
\cr
&\, +  c\intl_{0}^{t} S_{ 2,1}(\osc_{p, 1}(G(\tau ); x_0) )  d\tau \bigg\}
 \exp \intl_{0}^T 
(\| \nabla v(\tau )\|_{ \infty} +|v(\tau )|_{ \mathscr{L}^1_{ q(p,1)}}) d\tau. 
  \label{C.18b}
 \end{align}    
Observing  \eqref{cond2}, using Lemma\,\ref{lem10.1}, we may apply $ S_{ 1, q}$ to both sides of  \eqref{C.18b}. Accordingly,  
\[
\sup_{ t\in [0,T]} S_{ 1,q} (\osc_{ p,1}(F(t); x_0)) < +\infty.
\]
For given $ t\in [0,T]$ we 
may choose $ \xi $ such $ \xi (t)=0$. Thus, the same holds for $ f(t)$ in place of $ F(t)$.
Now, we are able to apply $ S_{ 1,q}$ to both sides of  \eqref{C.29a}, which yields
 \begin{align}
&S_{1,q} \Big(\Big\{\osc_{ p, 1} (V(\tau ); x_0,  2^{ i+1}) |\nabla \dot{P}^1_{ x_0, 2^{ i+1}} (F(\tau ))|\Big\}\Big)
\cr
&\le c|v(\tau )|_{ \mathscr{L}^1_{ q(p,1)}}S_{ 1,q}  (\osc_{p, 1}(F(\tau ); x_0))
\cr
&\quad + 
 cS_{ 1,q}\Big(\Big\{ ( i^{ -})^{ \frac{1}{q'}}  \osc_{p, 1}(V(\tau ); x_0, 2^i) \Big\}\Big)
 \Big\{\Big(\sum_{i=-\infty}^{\infty}  2^{ -iq}(\osc_{p, 1}(F(\tau ); x_0, 2^i))^q\Big)^{ \frac{1}{q}}+ 
|\nabla \dot{P}^1_{ x_0, 1}(F(\tau ))|\Big\}. 
 \label{C.29c}
\end{align}
Applying  $ S_{ 1,q}$ to both sides of  \eqref{C.17a} multiplying the result by $ 2^{ -j}$ and letting 
$ j \rightarrow -\infty$, we infer   
 \begin{align}
& \hspace*{-1cm}\Big(\sum_{i=-\infty}^{\infty}  2^{ -iq}(\osc_{p, 1}(F(t); x_0, 2^i))^q\Big)^{ \frac{1}{q}}
\cr
& \hspace*{-1cm}\le   c|f_0|_{ \mathscr{L}^1_{ q(p,1)}}
 + c\intl_{0}^{t}  \| \nabla v(\tau )\|_{ \infty} 
\Big(\sum_{i=-\infty}^{\infty}  2^{ -iq}(\osc_{p, 1}(F(t); x_0, 2^i))^q\Big)^{ \frac{1}{q}}d\tau 
\cr
&\hspace*{-1cm} \,  +c\intl_{0}^{t}  \Big(\sum_{i=-\infty}^{\infty} ( i^{ -})^{q-1} (2^{ -i} \osc_{p, 1}(V(\tau ); x_0, 2^i))^q\Big)^{ \frac{1}{q}}
 \Big(\sum_{i=-\infty}^{\infty}  2^{ -iq}(\osc_{p, 1}(F(\tau ); x_0, 2^i))^q\Big)^{ \frac{1}{q}}+ 
|\nabla \dot{P}^1_{ x_0, 1}(F(\tau ))|\Big)d\tau 
\cr
&\, +  c\intl_{0}^{t} |g(\tau )|_{ \mathscr{L}^1_{ q(p,1)}}d\tau . 
  \label{C.29d}
 \end{align}

Next, we require  to estimate $ |\nabla \dot{P}^1_{ x_0, 1}(F(\tau ))|$ by the initial data $ f_0$ and $ g$. 
We apply $ \dot{P}^1_{ x_0, 1}$ to both sides  \eqref{C.15}. This gives 
\begin{equation}
\partial _t \dot{P}^1_{ x_0, 1}(F) + \dot{P}^1_{ x_0, 1}(V\cdot \nabla F) = \dot{P}^1_{ x_0, 1}(G)
\quad  \text{in}\quad Q_T. 
\label{C.29e}
 \end{equation} 
  Noting that $ \dot{P}^1_{ x_0, 1}( P^1_{ x_0,1}(V)\cdot \nabla \dot{P}^1_{ x_0, 1}(F))= P^1_{ x_0,1}(V)\cdot \nabla \dot{P}^1_{ x_0, 1}(F)$, and applying $ \nabla $ to both sides of  \eqref{C.29e},  we infer 
  \begin{align}
   &\frac{d}{dt}\nabla \dot{P}^1_{ x_0, 1}(F) +(\nabla \dot{P}^1_{ x_0,1}(V))\cdot \nabla \dot{P}^1_{ x_0, 1}(F)   
   \cr
   &\qquad =
 \nabla \dot{P}^1_{ x_0, 1}\Big(P^1_{ x_0,1}(V)\cdot \nabla \dot{P}^1_{ x_0, 1}(F) - V\cdot \nabla F\Big)+ 
 \nabla \dot{P}^1_{ x_0, 1}(G)\quad  \text{in}\quad [0,T].  
    \label{C.30}
  \end{align}
   On the other hand, 
 \begin{align*}
 &\nabla \dot{P}^1_{ x_0, 1}\Big(P^1_{ x_0,1}(V)\cdot \nabla \dot{P}^1_{ x_0, 1}(F) - V\cdot \nabla F\Big)  
 \\
 &= \nabla \dot{P}^1_{ x_0, 1}\Big((P^1_{ x_0,1}(V)-V)\cdot \nabla \dot{P}^1_{ x_0, 1}(F) \Big) 
 +\nabla \dot{P}^1_{ x_0, 1}\Big( V \cdot \nabla  (\dot{P}^1_{ x_0, 1}(F) - F)\Big)   
  \\
 &= \nabla \dot{P}^1_{ x_0, 1}\Big((P^1_{ x_0,1}(V)-V)\cdot \nabla \dot{P}^1_{ x_0, 1}(F) \Big) 
 -\nabla \dot{P}^1_{ x_0, 1}\Big( \nabla \cdot V (P^1_{ x_0, 1}(F) - F)\Big).     
 \end{align*}
 Inserting this identity into the right-hand side of  \eqref{C.30}, multiplying the result by  $
 \frac{\nabla \dot{P}^1_{ x_0, 1}(F)}{|\nabla \dot{P}^1_{ x_0, 1}(F)|} $, we get the following differential inequality 
 \begin{align*}
&\frac{d}{dt}|\nabla \dot{P}^1_{ x_0, 1}(F) | 
 \le c \|\nabla v\|_{ \infty} |\nabla \dot{P}^1_{ x_0, 1}(F) |  + c \|\nabla v\|_{ \infty}
 \osc_{ p,1} (F, x_0; 1) + |\nabla \dot{P}^1_{ x_0, 1}(G)|. 
 \end{align*}
  Integrating this inequality over $ (0,t)$ and applying integration by parts, we obtain 
  \begin{align}
    |\nabla \dot{P}^1_{ x_0, 1}(F(t)) | &\le |\nabla \dot{P}^1_{ x_0, 1}(F(0)) |  +  
    c\intl_{0}^t \|\nabla v(\tau )\|_{ \infty} |\nabla \dot{P}^1_{ x_0, 1}(F(\tau )) | d\tau 
    \cr
    &\qquad +\intl_{0}^t \|\nabla v(\tau )\|_{ \infty}  \osc_{ p,1} (F(\tau ), x_0; 1) d\tau  +
        \intl_{0}^t |\nabla \dot{P}^1_{ x_0, 1}(G(\tau ))| d\tau 
     \cr
    & \le \|f_0\|_{ \mathscr{\tilde{L}}^1_{q (p,1) }}   +  
     c\intl_{0}^t \|\nabla v(\tau )\|_{ \infty} \osc_{ p,1} (F(\tau ), x_0; 1)  d\tau
     \cr
     & \qquad +   \intl_{0}^t \|g(\tau )\|_{ \mathscr{\tilde{L}}^1_{q (p,1) }} d\tau,
          \label{C.31}
  \end{align}
  where   $ | z|_{\mathscr{\tilde{L}}^{ 1}_{  q(p, 0)} }$ stands for the semi norm 
 \[
 | z|_{\mathscr{\tilde{L}}^{ 1}_{  q(p, 1)} }=  | z|_{\mathscr{L}^{ 1}_{  q(p, 1)} } +\sup_{ x_0\in \R^{n}}  
 |\nabla \dot{P}^1_{ x_0, 1}(z)|. 
 \]
Combining  \eqref{C.29d} and  \eqref{C.31}, we arrive at  
  \begin{align}
&\Big(\sum_{i=-\infty}^{\infty}  2^{ -iq}(\osc_{p, 1}(F(t); x_0, 2^i))^q\Big)^{ \frac{1}{q}}
+    |\nabla \dot{P}^1_{ x_0, 1}(F(t)) |
\cr
& \le   c|f_0|_{ \mathscr{\tilde{L}}^1_{ q(p,1)}}
 + c\intl_{0}^{t}  \| \nabla v(\tau )\|_{ \infty} 
\Big(\sum_{i=-\infty}^{\infty}  2^{ -iq}(\osc_{p, 1}(F(t); x_0, 2^i))^q\Big)^{ \frac{1}{q}}d\tau 
\cr
& \, \quad  +c\intl_{0}^{t}  \Big(\sum_{i=-\infty}^{\infty} ( i^{ -})^{q-1} (2^{ -i} \osc_{p, 1}(V(\tau ); x_0, 2^i))^q\Big)^{ \frac{1}{q}}
 \bigg\{\Big(\sum_{i=-\infty}^{\infty}  2^{ -iq}(\osc_{p, 1}(F(\tau ); x_0, 2^i))^q\Big)^{ \frac{1}{q}}
\cr
&\, \qquad \quad + 
|\nabla \dot{P}^1_{ x_0, 1}(F(\tau ))|\bigg\}d\tau +  c\intl_{0}^{t} |g(\tau )|_{ \mathscr{\tilde{L}}^1_{ q(p,1)}}d\tau.
  \label{C.32}
 \end{align}

  Applying Gronwall's lemma and for given $ t\in [0,T]$ choosing $ \xi $ such that $ \xi (t)=0$, and taking the supremum  over $ x_0\in \R^{n} $, we obtain 
   the desired estimate \eqref{estimate2}.   \hfill \Beweisende

\vspace{0.3cm}
{\bf Proof of \eqref{estimate3} in Theorem\,\ref{thm3}}. We first define
 \[
\chi (x_0, t) = \sup_{ j\in \Z}2^{ -j}\osc_{ p, 0} (F(t); x_0,  2^{ j}),\quad (x_0, t)\in \R^{n} \times  [0,T].  
\]
Clearly, thanks to  \eqref{C.13} $ \chi (x_0, t)$ is finite. 
Noting 
that $ \|\nabla P^N_{ x_0, 2^{ j+1}}(F(\tau ))\|_{ L^\infty(B(x_0, 2^{ j+1}))} \le c\chi (x_0, \tau )$, 
we get from  \eqref{C.17} with $ L = 2N-1$
\begin{align}
&\osc_{ p, 2N-1}(F(t); x_0, 2^j)
\le   c \osc_{ p, N} (F(0); x_0, 2^{ j+1} )
 \cr
&\quad  +c \intl_{0}^{t}   \|  \nabla v(\tau )\|_{\infty} 
\osc_{ p, N} (F(\tau ); x_0, 2^{ j+2})d\tau  
\cr
& \quad +  c 2^{ -j} \intl_{0}^{t}  \osc_{ p, N}(V(\tau ); x_0, 2^{ j+1})
 \chi (x_0, \tau )d\tau 
 +  c\intl_{0}^{t} \osc_{ p, N} (G(\tau ); x_0, 2^{ j+1})d\tau.
  \label{C.33}
 \end{align}  
 First let us estimate the term  $ \osc_{ p,0}(F(t); x_0, 2^{ j+1})$. In view of   \eqref{C.17a}  with $j+1$ 
 in place of $ j$, and recalling that $\nabla  f_0\in L^\infty(\R^{n} ), g\in L^1(0,T; L^\infty(\R^{n} ))$,  we see that 
 \begin{align}
& \osc_{ p, 0} (F(t); x_0,  2^{ j+1})
\cr
&\quad \le   c 2^{ j}\|\nabla f_0\|_{ \infty}+c \intl_{0}^{t} \|   \nabla v(\tau )\|_{ \infty}  
\osc_{ p, 0} (F(\tau ); x_0,  2^{ j+3})d\tau  
+  c\intl_{0}^{t} 2^j \|\nabla g(\tau )\|_{\infty }d\tau.
  \label{C.35}
 \end{align}       
 Multiplying both sides of  \eqref{C.35} by $ 2^{ -j}$ and taking the supremum over all $ j\in \Z$, using the triangle inequality, we obtain 
 \begin{align}
 \chi (x_0, t)  \le c  
   \|\nabla f_0\|_{ \infty} + 
c \intl_{0}^{t} \|   \nabla v(\tau )\|_{ \infty}   \chi (x_0, \tau )  d\tau  
+  c\intl_{0}^{t} \|\nabla g(\tau ) \|_{ \infty} d\tau,
\label{C.36b}
\end{align}
 Thanks to  \eqref{C.13} we have 
$ S_{ N+1, 1}(\osc_{ p,N}F(t); x_0) < +\infty$ for all $ t\in [0,T]$.  
Applying $ S_{ N+1,1}$ to both sides of  \eqref{C.33},  and using Corollary\,\ref{cor3.10} with $ N'=2N-1$, we get 
 \begin{align}
\osc_{ p, N}(F(t); x_0) 
& \le  S_{ N+1, 1}(\osc_{ p, 2N-1}(F(t); x_0) ) 
\cr
&\le   cS_{ N+1, 1}( \osc_{ p, N} (F(0); x_0))
   +c \intl_{0}^{t}   \|  \nabla v(\tau )\|_{\infty} 
S_{ N+1, 1}(\osc_{ p, N} (F(\tau ); x_0))d\tau  
\cr
& \quad +  c  \intl_{0}^{t} S_{ N+1, 1}( \osc_{ p, N}(V(\tau ); x_0))
\chi (x_0, \tau ) d\tau 
 \cr
 &\qquad +  c\intl_{0}^{t}  S_{ N+1, 1} (\osc_{ p, N}(G(\tau ); x_0))d\tau.
  \label{C.38}
 \end{align}      
Next, once more using  \eqref{C.13} we see that  $ S_{ s,q}(\osc_{ p, N}(F(t); x_0))< +\infty$, for all $ t\in [0,T]$. Thus,  we apply $ S_{s, q}$ to both sides of \eqref{C.38} and use  Lemma\,\ref{lem10.1}. This combined with  \eqref{C.36b} gives 
\begin{align}
& 2^{ -js} (S_{ s , q}(\osc_{ p, N}(F(t); x_0) ))_j + \chi (x_0, t)
\cr
&\le   c |f_0|_{ \mathscr{L}^s_{q(p,1)}} +  \chi (x_0, 0)
\cr
& \quad +  c  \intl_{0}^{t} (|v(\tau )|_{ \mathscr{L}^s_{q(p,1)}} +\|  \nabla v(\tau )\|_{\infty})
 \Big[ 2^{ - js} \Big(S_{s , q}\Big(\osc_{ p, N}(F(t); x_0) \Big)\Big)_j +\chi (x_0, \tau ) \Big]d\tau
 \cr
  &\qquad   +  c\intl_{0}^{t} (|g(\tau )|_{ \mathscr{L}^s_{q(p,1)}} + \|\nabla g(\tau )\|_{ \infty})d\tau.
  \label{C.38a}
 \end{align}      
 By virtue of Gronwall's lemma we deduce from \eqref{C.38a} 
\begin{align}
& 2^{ -js} (S_{ s , q}(\osc_{ p, N}(F(t); x_0) ))_j + \chi (x_0, t)
\cr
&\le  c \bigg\{ |f_0|_{ \mathscr{L}^s_{q(p,1)}} + \|\nabla f_0\|_{ \infty}
\cr
& \quad +   \intl_{0}^{T} (|g(\tau )|_{ \mathscr{L}^s_{q(p,1)}} + \|\nabla g(\tau )\|_{ \infty})d\tau\bigg\}
\exp \bigg(\intl_{0}^{T} (|v(\tau )|_{ \mathscr{L}^s_{q(p,1)}} +\|  \nabla v(\tau )\|_{\infty})
\bigg) d\tau. 
  \label{C.38b}
 \end{align}    
 Whence, \eqref{estimate3}.    \hfill \Beweisende 

\vspace{0.3cm}
{\bf Proof of \eqref{estimate4} in Corollary\,\ref{cor1.4}}.  In view of  Theorem\,\ref{thm5.4} we have 
$ \nabla f_0\in L^\infty(\R^{n} ), \nabla g\in L^1(0,T; L^\infty(\R^{n} ))$. More precisely,  
 \eqref{5.20} yields 
\begin{align*}
\|\nabla f_0\|_{ \infty} \le  c\|f_0\|_{\mathscr{L}^1_{1(p,1)} },\quad    
\intl_{0}^T\|\nabla g(\tau )\|_{ \infty}d\tau  \le c \|g\|_{L^1(0,T; \mathscr{L}^1_{1(p,1)}) }.    
\end{align*}
In particular, this shows that condition  \eqref{cond3} of Theorem\,\ref{thm2} is fulfilled. 
Furthermore, since $ v\in L^1(0,T; \mathscr{L}^1_{1(p,1)}(\R^{n} ))$, condition of Theorem\,\ref{thm2} \eqref{cond2} is also  satisfied. 
Now, we are in a position to apply of Theorem\,\ref{thm2}, which yields $ f\in L^\infty(0,T; \mathscr{L}^1_{1(p,1)}(\R^{n} ))$.  This allows to apply $ S_{ 1,1}$ to both sides of   \eqref{C.17}. This  together with Gronwall's Lemma and the inequality $ | \nabla \dot{P}^1_{ x_0, 2^{ j+1}}(F(\tau ))| \le c 2^{ -j}\|\nabla f(\tau )\|_\infty$ 
yields 
\begin{align}
&(S_{ 1, 1}(\osc_{p, 1}(F(t); x_0)))_j
\cr
&\le   c  \bigg\{(S_{ 1, 1}(\osc_{p, 1}(F (0); x_0) )_j
 +  \intl_{0}^{t} ( S_{ 1, 1}(\osc_{p, 1}(G(\tau ); x_0) ) )_j d\tau
   \cr
  &\qquad \qquad  +  \intl_{0}^{t} ( S_{ 1, 1}(\osc_{p, 1}(V(\tau ); x_0) ) )_j \|\nabla f(\tau )\|_{ \infty} d\tau\bigg\} 
   \exp \bigg(c\intl_{0}^{t} \| \nabla v(\tau )\|_{ \infty} d\tau\bigg). 
  \label{C.40}
 \end{align}     
Choosing $ \xi $ so that $ \xi (t)=0$, multiplying both sides by $ 2^{ -j}$ and letting $ j \rightarrow -\infty$ taking the supremum over $ x_0\in \R^{n} $, 
we deduce from  \eqref{C.40} 
 \begin{align}
&|f(t)|_{\mathscr{L}^1_{1(p,1)} }
\cr
&\le   c  \bigg\{\|f_0\|_{\mathscr{L}^1_{1(p,1)} }
 +  \intl_{0}^{t} \|g(\tau )\|_{\mathscr{L}^1_{1(p,1)} }d\tau
   \cr
  &\qquad \qquad  + \intl_{0}^{t} \|v(\tau )\|_{\mathscr{L}^1_{1(p,1)} }\|\nabla f(\tau )\|_{ \infty} d\tau\bigg\} 
   \exp \bigg(c\intl_{0}^{t} \| \nabla v(\tau )\|_{ \infty} d\tau\bigg). 
  \label{C.41}
 \end{align}     
 Combining  \eqref{C.41} and  \eqref{C.11} along with  \eqref{C.13} in order to estimate $ \|f(t)\|_{ L^p(B(1))}$, we get the desired estimate  \eqref{estimate4}.    \hfill \Beweisende  
 
 \vspace{0.3cm}
 
 Below we prove the uniqueness parts of Theorem \ref{thm1}, Theorem\ref{thm2}, Theorem\ref{thm3} and Corollary\ref{cor1.4}. In fact we prove the stronger version of it, namely the strong-weak uniqueness.
 
  \vspace{0.3cm} 
  
 {\bf  Strong-weak uniqueness}.  Let $ \overline{f}  \in L^2_{loc}(\R^{n} )$  be a weak  solutions to \eqref{trans}. 
 Then $ w=f-\overline{f} $ solves the transport equation with homogenous data 
 \begin{equation}
 \partial _t w + (v \cdot \nabla) w = 0\quad  \text{ in}\quad  Q_T,\quad  w=0 \quad  \text{ on}\quad  \R^{n}\times \{ 0\}
  \label{C.25}
 \end{equation}
 in a weak sense, i.e.  for all $ t\in (0,T)$, and for all  
 $\varphi \in L^\infty(0,t; W^{ 1,2}(\R^{n}))\cap W^{1,1}(0, t; L^2(\R^{n}))$ with $ \supp(\varphi ) \Subset \R^{n}\times [0,t]$, it holds 
 \begin{align}
 -  \intl_{0}^{t}  \intl_{ \R^{n}} w \partial _t \varphi + (v\cdot \nabla) \varphi w + \nabla \cdot v\varphi w  dx ds = -\intl_{ \R^{n}} w(t) \varphi (t)  dx.    
  \label{C.26}
 \end{align}
 Let $ \psi \in C^{\infty}_{c}(\R^{n})$ be a given function.  
 Using the method of characteristics,  for every $ \varepsilon >0$ we get a solution  $ \varphi^\varepsilon   \in L^\infty(0,t;  W^{ 1,2}(\R^{n}))\cap W^{ 1,1}(0, t; L^2(\R^{n}))$  of the the following dual problem 
 \begin{equation}
 \partial _t \varphi^\varepsilon  + v \cdot \nabla \varphi ^\varepsilon  + \nabla \cdot v_\varepsilon  \varphi^\varepsilon   =0 \quad  \text{ in}\quad  Q_t, 
 \quad  \varphi^\varepsilon  (t) = \psi\quad  \text{ in }\quad  \R^{n}.  
 \label{C.27}
 \end{equation}
 Noting that $ \|\nabla v_\varepsilon (\tau )\|_\infty \le  \|\nabla v (\tau )\|_\infty$, using Gronwall's lemma we see that $ \| \varphi ^\varepsilon  \|_{ 1} + \| \varphi^\varepsilon   \|_{ \infty} \le c$ with a  constant $ c>0$ independent of 
 $ \var>0$.   Since $  v(0, \cdot ),\|\nabla v(\cdot )\|_{ \infty} \in L^1(0,T)$  using  \eqref{C.12},  we get a number $ 0< R< +\infty$ 
 such that $ \supp (\varphi^\varepsilon   ) \subset B(R) \times [0, t]$. In  \eqref{C.26} putting $ \varphi = \varphi _\var $, and using \eqref{C.27},  we infer 
 \begin{align}
\intl_{ \R^{n}} w(t) \psi   dx  & =  \intl_{0}^{t}  \intl_{ \R^{n}} w \partial _t \varphi _\var + (v\cdot \nabla) \varphi_\var  w + \nabla \cdot v\varphi^\var  w  dx ds
  \cr
&  =  \intl_{0}^{t}  \intl_{ B(R)}  \nabla \cdot (v- v_\var )\varphi^\var  w  dx ds.   
 \label{C.28}
 \end{align}
 Noting that $ \nabla \cdot (v(s)- v_\var(s) ) \rightarrow 0$  in $ L^{ 2}(B(R))$ as $ \var  \searrow 0$ for almost all 
 $ s\in (0,t) $, by the aid of Vitali's convergence theorem (\cite[p. 180]{fol}) it follows that 
 \[
\intl_{0}^{t}  \intl_{ B(R)}  \nabla \cdot (v- v_\var )\varphi^\var  w  dx ds \rightarrow 0\quad  \text{ as}\quad  \var \searrow 0.
\]
 
Letting $ \var \searrow 0$  in \eqref{C.28}, we deduce that $ \intl_{ \R^{n}} w(t) \psi   dx=0$. Whence, $ w \equiv 0$.   This shows the uniqueness.  
  \hfill \Beweisende 
  
  \appendix
  %%% ----------------------------------------------------------------------
%       SECTION A
%%% ----------------------------------------------------------------------
\section{Minimal polynomials}
\label{sec:-A}
\setcounter{secnum}{\value{section} \setcounter{equation}{0}
\renewcommand{\theequation}{\mbox{A.\arabic{equation}}}}

Let $ <p< +\infty$.  Let $ x_0\in \R^{n}$ and $ 0< r<+\infty$ be fixed. Set $ \phi = \varphi(r^{ -1}(x_0-\cdot )) $, 
where $ \varphi  \in C^{\infty}_{c}(B(1))$, being radial symmetric, stands for the standard mollifier.  
For $ \delta  \ge 0$ we define the following functional $ 
J_\delta : L ^p(B(x_0, r)) \rightarrow \R$ by 
\[
J_\delta (f) = \intl_{B(x_0, )} (\delta + | f|^2 )^{ \frac{p}{2}} \phi ^p dx,\quad    f\in L ^p(B(x_0, r)). 
\]
Recall $ \mathcal{P}_N$, $ N\in \N_0$,  denotes the space of all polynomial of degree less or equal $ N$. Since 
$ J_\delta $ is strict  convex and lower semi continuous with $ J_\delta (f) \rightarrow +\infty $ as $ \| f\|_{ L ^p(B(x_0, r))}
\rightarrow +\infty$.  For each $ f\in L^p(B(x_0, r))$ there exists a unique $ P ^{ N, \delta  }_{ x_0, r}(f) \in \mathcal{P}_N$
with 
\begin{equation}
J_\delta (P ^{ N, \delta  }_{ x_0, r}(f)-f) = \min_{P\in \mathcal{P}_N }J_\delta (P-f)
\label{A.1a}
\end{equation}
Clearly, the mapping $J_{ \delta, f }: P \mapsto  J_\delta(P-f) $ is  differentiable  as a function from $ \mathcal{P}_N$ into $ \R$. Since the first variation 
must vanish at each minimizer, we get 
\begin{equation}
 \langle D J_{ \delta, f}(P ^{ N, \delta  }_{ x_0, r}(f), P)\rangle=0\quad  \forall\,P \in \mathcal{P}_N. 
\label{A.1b}
\end{equation}
This shows that 
\begin{equation}
 \intl_{B(x_0, r)} F_\delta ( P ^{ N, \delta  }_{ x_0, r}(f)- f)\cdot P   \phi ^pdx =0  \quad  \forall\,P \in \mathcal{P}_N. 
\label{A.1b}
\end{equation}
where 
\[
F_\delta ( u) =  (\delta + | u|^2)^{  \frac{p-2}{2}} u, \quad  u\in \R^{n}. 
\]
It is well known that $ F_\delta $ is monotone and continuously differentiable  for  each $ \delta >0$ . Furthermore, 
there exists a constant $ c>0$ independent of $ \delta $ such that for all $ u,v\in \R^{m}$,
\begin{align}
(F_\delta (u)- F_\delta (v))(u-v) &\ge c (p-1)  (\delta + | u|+ | u-v|)^{ p-2}| u-v|^{2},
\label{A.1}
\\
|F_\delta (u)- F_\delta (v) | &\le c p (\delta + | u|+ | u-v|)^{  \frac{p-2}{2}}| u-v|. 
\label{A.2}
\end{align}
We now define the mapping $ G_\delta: L^p(B(x_0, r))\times \mathcal{P}_N \rightarrow (\mathcal{P}_N)'$ by
\[
 \langle G_\delta (f, P), Q\rangle = \intl_{B(x_0, r)} F_\delta (f(x)-P)\cdot Q  \phi^2 (x) dx,\quad f\in L^p(B(x_0,r), \,\,    P, Q\in \mathcal{P}_N. 
\]
Clearly, \eqref{A.1b} is equivalent to 
\begin{equation}
G_\delta ( f, P ^{ N, \delta  }_{ x_0, r}(f))=0. 
\label{A.1c}
\end{equation}

We obtain the following properties of $ G_\delta $. 

\begin{lem}
\label{lemA.1}  
1. For every $ f\in L^p(B(x_0, r))$ the mapping $ G_\delta (f, \cdot  ): \mathcal{P}_N\rightarrow (\mathcal{P}_N)'$ is strictly monotone, bijective, 
and in case $ \delta >0$ stronly monotone and is a $ C^1$ diffeomorphism. 

2. In case $ \delta >0$, the mapping $ f \mapsto P ^{ N, \delta  }_{ x_0, r}(f): L^p(B(x_0, r)) \rightarrow \mathcal{P}_N $ is Fr\'echet differentiable, and its derivative 
is given by 
\begin{equation}
DP ^{ N, \delta  }_{ x_0, r}(f) =- [D_2 G_\delta (f, P ^{ N, \delta  }_{ x_0, r}(f))]^{ -1} \circ D_1
G_\delta (f, P ^{ N, \delta  }_{ x_0, r}(f)),\quad  f\in L^p(B(x_0, r)),
\label{A.2c}
\end{equation} 
where $ D_1 G_\delta (f, P)\in \mathscr{L}(L^p(B(x_0,r)), (\mathcal{P}_N)')$ stands for derivative with respect to the first variable, while 
$ D_2 G_\delta (f, P)\in \mathscr{L}( \mathcal{P}_N,  (\mathcal{P}_N)')$ stands for derivative with respect to the second variable. 

Furthermore  it holds for every $ f\in L^p(B(x_0, r))$
\begin{equation}
\| P ^{ N, \delta  }_{ x_0, r} (f)\|^p_{ L^p(B(x_0,  \frac{r}{2}))}  \le 2^p\intl_{B(x_0, r)} (\delta + | f| ^2)^{ \frac{p}{2}} \phi ^p dx.  
\label{A.2b}
\end{equation}

3. For all $ f\in L^p(B(x_0, r))$ it holds
\begin{equation}
P ^{ N, \delta  }_{ x_0, r}(f) \rightarrow P^{ N, \ast}_{ x_0, r}(f)  \quad  \text{{\it in}}\quad  \mathcal{P}_N\quad  \text{{\it as}}\quad  \delta  \searrow 0, 
\label{A.2e}
\end{equation}
where $ P^{ N, \ast}_{ x_0, r}(f) = P^{ N, 0}_{ x_0, r}(f)$. 

\end{lem}

{\bf Proof}: 1. Observing \eqref{A.1},  we get for  all $ f\in L^p(B(x_0, r))$, and  $ P, Q\in \mathcal{P}_N$ 
\begin{align}
& \langle (G_\delta (f, P)- G_\delta (f, Q)), (P-Q) \rangle
\cr
&\quad \ge c (p-1) {\D \intl_{B(x_0, r)}}     (\delta + | P-f|+ | P-Q|)^{ p-2} | P-Q|^{2}\phi ^2dx .
\label{A.3}
\end{align}
This immediately shows that $ G_\delta( f, \cdot ) $ is strictly monotone and in case $ \delta >0$ strongly monotone. 
Here we have used the fact that $ \| P\|_{ L^2(B(x_0, r))}$ defines an equivalent norm on $ \mathcal{P}_N$. 
Furthermore, if $ \delta >0$ we see that   $ G_\delta (f, \cdot ): \mathcal{P}_N $ into $ (\mathcal{P}_N)'$ 
is continuously differentiable and coercitive, i.e. 
\[
\frac{ \langle G_\delta (f, P), P\rangle }{\| P\|} \rightarrow 0\quad  \text{ as}\quad  \| P\| \rightarrow +\infty.  
\]
Applying the theory of monotone operators, we see that $ G_\delta (f, \cdot )$ is bijective, and is a
$ C^1$ diffeomophism.  
  
\vspace{0.2cm}
2. Let $ \delta >0$ and $ f\in L^p(B(x_0, r))$.  Let $ P ^{ N, \delta  }_{ x_0, r}(f)\in \mathcal{P}_N$ denote 
the minimizer of the functional $ J_\delta (\cdot -f)$ in $ \mathcal{P}_N$.  In view of \eqref{A.1c}  we have 

$ G_\delta (f, P ^{ N, \delta  }_{ x_0, r}(f))=0$. Since $ D_2 G_\delta $ is an isomorphism from $ \mathcal{P}_N$
 into $ (\mathcal{P}_N)'$,  by the implicit function theorem we infer  that the mapping  
 $ P ^{ N, \delta  }_{ x_0, r}: L^p(B(x_0, r)) \rightarrow \mathcal{P}_N$ is Fr\'echet, 
 differentiable, and it holds \eqref{A.2c}.  

\vspace{0.2cm}
{\it Proof of \eqref{A.2b}. } Since $ J_\delta$ is convex and recalling the minimizing property of  $ P ^{ N, \delta  }_{ x_0,r} (f)$, we get 
\begin{align*}
 J_\delta \Big( \frac{ P ^{ N, \delta  }_{ x_0,r} (f)}{2}\Big)&\le \frac{1}{2} (J_\delta (P ^{ N, \delta  }_{ x_0,r} (f)-f) + 
 J_\delta (f))  
\le   J_\delta (f). 
\end{align*}
This shows that 
\[
2^{ -p}\intl_{B(x_0, r)} |  P ^{ N, \delta  }_{ x_0,r} (f)|^p  dx \le   J_\delta (f). 
\]
Whence, \eqref{A.2b}    

\vspace{0.3cm}
3. Now, let $ \delta _k \searrow 0$ as $ k \rightarrow +\infty$. By \eqref{A.2b} we see that 
$ \{ P ^{ N,  \delta _k}_{ x_0, r}(f)\}$ 
is bounded. Thus, there exists a subsequence, and $P^{N, \ast}_{ x_0, r }(f)\in  \mathcal{P}_N$ such that $ P^{N,  \delta_{ k_j}}_{ x_0, r}(f) \rightarrow 
P^{N, \ast}_{ x_0, r }$    in $ \mathcal{P}_N$ as $ j \rightarrow +\infty$. 
 Since $ F_{ \delta ^{ k_j}} (f(x)- P^{N,  \delta_{ k_j}}_{ x_0, r}(f)) \rightarrow F_0(f(x)- P^{N, \ast}_{ x_0, r }(f))$ as $ j \rightarrow +\infty$ for all $ x\in B(x_0, r)$ by 
Lebesgue's theorem of dominated convergence  it follows $ 0= G_{\delta_{ k_j}} (P^{N,  \delta_{ k_j}}_{ x_0, r}(f)) \rightarrow 
G_0(f, P^{N, \ast}_{ x_0, r })$. Since $ G_0(f, \cdot )$ is strictly  monotone, the zero is unique, and thus 
$ P^{N, \ast}_{ x_0, r }(f) = P^{N,0}_{ x_0, r } (f)$. Thus,  convergence property \eqref{A.2e} 
is verified.  

\vspace{0.2cm}
Furthermore, in \eqref{A.2b} letting $ \delta \searrow 0$, we see that 
\begin{equation}
\| P^{N, \ast}_{ x_0, r } (f)\|_{ L^p(B(x_0, \frac{r}{2}))} \le 2 \| f \phi \|_{ L^p(B(x_0, r))}. 
\label{A.9}
\end{equation} 

This completes the proof of the lemma.  \hfill \Beweisende 

\begin{rem}
\label{remA.2}
The mapping $ P _{ x_0, r}^{ N, \delta }: L^p(B(x_0, r)) \rightarrow \mathcal{P}_N$ fulfills  the projection property 
\begin{equation}
P_{ x_0, r}^{ N, \delta }(Q) =Q\quad  \forall\,Q\in \mathcal{P}_N. 
\label{A.10}
\end{equation}
In fact, this follows immediately from \eqref{A.1a} by setting $ f=Q$ therein.  

%%% ----------------------------------------------------------------------
%       SECTION 
%%% ----------------------------------------------------------------------
\section{Example of a function in $  \mathscr{L}^{1}_{1 (p, 1)} (\R^n) \setminus  C^1(\R^n)$}
\label{sec:-B}
\setcounter{secnum}{\value{section} \setcounter{equation}{0}
\renewcommand{\theequation}{\mbox{B.\arabic{equation}}}}

The following example shows that $  \mathscr{L}^{1}_{1 (p, 1)} (\R^n)$ is not in $ C^1(\R^{n} )$. 
For simplicity we only consider the case $ n=1$ since general case $ n\in \N$ can be reduced to $ n=1$.   
We define  
\[
f(x)= \intl_{0}^x u(y) dy,\quad  x\in \R,
\]
where 
\[
u(x)=\begin{cases}
1- 2^{ 2m} |x-  2^{ -m}|\quad   &\text{if}\quad x\in I_m, \quad m\in \N, 
\\
0\quad   &\text{elsewhere},
\end{cases}
\]
and $ I_m = [2^{ -m}-2^{ -2m}, 2^{ -m}+2^{ -2m})$. 

\vspace{0.2cm}
{\bf Proof of $ f\in  \mathscr{L}^{1}_{1 (p, 1)} (\R)$}: Thanks to  \eqref{5.29c} it will be sufficient to show that 
$ u\in  \mathscr{L}^{0}_{1 (p, 0)} (\R) $.  In what follows we estimate $\osc_{ p, 0} (u; x, r) $ for $ x\in [0,1]$ and 
$ 0< r <+\infty$. 

We start with the  case  $ x=0$. For   $ 2^{- m-1} < r \le  2^{ -m}$ we get 
\[
\osc_{ p, 0} (u; 0, r) \le 2 \bigg( \frac{1}{2r}\intl_{-r}^{ r} |u(y)|^p dy\bigg)^{ \frac{1}{p}} 
\le 2 \bigg( \frac{1}{2r}\intl_{-r}^{ r}  \sum_{j= m}^{\infty} \chi _{ I_j} dy\bigg)^{ \frac{1}{p}} 
\le c r^{ - \frac{1}{p}} \Big(\sum_{j=m}^{\infty} 2^{ - 2j}\Big)^{ \frac{1}{p}}
\le c 2^{ - \frac{m}{p}}. 
\]
This yields, 
\begin{align}
   \sum_{j=-\infty}^{+\infty}\osc_{ p, 0} (u; 0, 2^{ j}) 
& =  \sum_{j=-\infty}^{-1}\osc_{ p, 0} (u; 0, 2^{ j})  + \sum_{j=0}^{\infty}\osc_{ p, 0} (u; 0, 2^{ j})  
\cr
&\le c \sum_{j=-\infty}^{-1} 2^{ \frac{j}{p}}  + c\sum_{j=0}^{\infty} 2^{- \frac{j}{p}} < +\infty. 
 \label{B.1}
\end{align}

Let $ x\in (0, 1]$.  Then there exists $ m\in \N$ such that $2^{ -m} < x \le  2^{ -m+1} $.  Let $ 0<r< +\infty$. 
We consider the following  three cases. 

1. First, in case $ 2^{-m-1}< r< +\infty$ by triangle inequality we get 
\[
\osc_{ p, 0} (u; x, r) \le c \osc_{ p, 0} (u; 0, 8r).  
\]
2. In  case $ 2^{ -2m}< r \le 2^{ -m-1}$, again  by triangle inequality we find 
\begin{align*}
\osc_{p, 0} (u; x, r) &\le 
2   \bigg(\frac{1}{2r}\intl_{-r}^{ r} |u|^p dy\bigg)^{ \frac{1}{p}} \le 2\bigg(\frac{1}{2r}\intl_{-r}^{ r} ( \chi _{I_{ m+1}  } +  \chi _{I_{ m}  }+ 
 \chi _{I_{ m-1}  }) dy\bigg)^{ \frac{1}{p}} \le c r^{- \frac{1}{p}} 2^{ -\frac{2m}{p}}. 
\end{align*}
3. In case $ 0 < r \le 2^{ -2m}$,  using Poincar\'e's inequality,  we obtain 
\begin{align*}
\osc_{ p, 0} (u; x, r) &\le c  r\bigg(\frac{1}{2r}\intl_{-r}^{ r}  |u'(y)|^pdy\bigg)^{ \frac{1}{p}} 
\\
&\le  cr^{1- \frac{1}{p}}\bigg(\intl_{-r}^{ r} (2^{2( m+1)} \chi _{I_{ m+1}  } + 2^{ 2m} \chi _{I_{ m}  }+ 
2^{ 2(m-1)} \chi _{I_{ m-1}  })^p dy\bigg)^{ \frac{1}{p}} \le c r^{ \frac{1}{p'}} 2^{  \frac{2m}{p'}},
\end{align*}
where $ p' = \frac{p}{p-1}$. 

Using the  the estimates above together with  \eqref{B.1}, we obtain 
\begin{align*}
  \sum_{j\in \Z} \osc_{ p, 0} (u; x, 2^j) &= \sum_{j=-m+1}^{\infty}  \osc_{ p, 0} (u; x, 2^j)   +
\sum_{j=-2m+1}^{-m}  \osc_{ p, 0} (u; x, 2^j)   + \sum_{j=-\infty}^{-2m}  \osc_{p, 0} (u; x, 2^j)  
\\
&\le  c \sum_{j=-m+1}^{\infty}  \osc_{ p, 0} (u; 0, 2^{ j+3})  + c 2^{-  \frac{2m}{p}}\sum_{j=-2m+1}^{-m} 2^{ -\frac{j}{p}}
\\
&\qquad + 2^{ \frac{2m}{p'} } \sum_{j=-\infty}^{-2m}  2^{ \frac{j}{p'}} \le c, 
\end{align*}
where the $ c$ stands for an absolute constant. Accordingly,
\begin{equation}
 \sup_{ x\in [0,1]} \sum_{j\in \Z} \osc_{ p, 0} (u; x, 2^j) <+\infty. 
\label{B.2}
 \end{equation}

In case $ x<0$ there exists $ m\in \Z$ such that  $ -2^{ m+1} < x \le - 2^{ m}$. Using triangle inequality 
together with  \eqref{B.1}, we easily see that 
\[
 \sum_{j\in \Z} \osc_{ p, 0} (u; x, 2^j) \le  \sum_{j= m}^\infty \osc_{ p, 0} (u; x, 2^j) \le   c   
  \sum_{j\in \Z} \osc_{ p, 0} (u; 0, 2^j) \le c.
\]
Similarly by the aid of  \eqref{B.2} we get   $ \sum_{j\in \Z} \osc_{ p, 0} (u; x, 2^j) \le c  
\sum_{j\in \Z} \osc_{ p, 0} (u; 1, 2^j) \le  c$ for all $ x \ge 1$. 
This shows that $ u\in \mathscr{L}^{0}_{1 (p, 0)} (\R^{n} )$, and thus $ f\in \mathscr{L}^{1}_{1 (p, 1)} (\R^{n} )$ but $ u \notin C^1(\R)$. 

 \hfill \Beweisende

\end{rem}
\hspace{0.5cm}
$$\mbox{\bf Acknowledgements}$$
Chae was partially supported by NRF grants 2016R1A2B3011647, while Wolf has been supported 
supported by NRF grants 2017R1E1A1A01074536.
The authors declare that they have no conflict of interest.


\begin{thebibliography}{99} 

\bibitem{bah}
H. Bahouri, J. -Y. Chemin  and R. Danchin, {\it Fourier Analysis and Nonlinear Partial Differential equations,}  Springer (2011).
%\bibitem{bar} C. Bardos, C., E. Titi, {\it Loss of smoothness and energy conserving rough weak solutions for the 3d Euler equations,} Discret. Contin. Dyn. Syst. Ser. S {\bf 3}(2), 185-197 (2010).
%\bibitem{bea} J.T. Beale, T., Kato, A.  Majda, {\it Remarks on the breakdown of smooth solutions for the 3D Euler equations,} Comm. Math. Phys. {\bf 94,} 61-66 (1984).
  \bibitem{bou}G. Bourdaud, {\it  $L^p$-estimates for certain non-regular pseudo-differential operators,}  Comm. Partial Diff. Eqs. {\bf 7},   1023-1033 (1982). 
 \bibitem{bour1} J. Bourgain and D. Li, {\it Strong illposedness of the incompressible Euler equations in integer $C^m$ spaces}, Geom. Funct. Anal. {\bf 25},  1-86  (2015).
  \bibitem{bour2} J. Bourgain and D. Li, {\it Strong illposedness of the incompressible Euler equations in  Borderline Sobolev spaces}, Invent. Math. {\bf 201}, 97-157 (2015). 
 \bibitem{ca} S. Campanato,  {\it Proprieti di una famiglia di spazi funzionali}, Ann. Sc.
Norm. Sup. Pisa {\bf 18}, 137-160 (1964).
\bibitem{cw}D. Chae and J. Wolf, {\it The Euler equations  in  a critical case of the generalized Campanato space,}  arxiv preprint (2019).
%Asymptot. Anal. {\bf 38}(3-4), 339-358 (2004).
\bibitem{che}J.-Y. Chemin, {\it Perfect incompressible fluids. Clarendon press}, Oxford (1998).
%\bibitem{ches}A. Cheskidov and R.  Shvydkoy, {\it  Ill-posedness of the basic equations of fluid dynamics in
%Besov spaces,} Proc. Amer. Math. Soc. {\bf 138} (3), 1059-1067 (2010).
%\bibitem{con}  P. Constantin, {\it On the Euler equations of incompressible fluids,} Bull. Amer. Math. Soc. 
%N.S. {\bf  44}(4), 603-621 (2007).
\bibitem{dip}R.J. DiPerna and P.L. Lions, {\it Ordinary differential equations, transport theory
and Sobolev spaces,}  Invent. Math., {\bf 98}, 511-547 (1989).
 \bibitem{fol} G.B. Folland,  Real analysis. Pure and Applied Mathematics (Second ed.), New York: John Wiley, (1999). 
 \bibitem{gia}M. Giaquinta, {\it Multiple integrals in the calculus of variations and nonlinear elliptic systems,} Ann. Math. Stud. No. 105, Princeton Univ. Press, Princeton (1983). 
% \bibitem{mis} G. Misioekand T. Yoneda, {\it  Continuity of the solution map of the Euler equations in H\"{o}lder spaces and weak norm inflation in Besov spaces,} Trans. Amer. Math. Soc. {\bf 370} (2018), no. 7, 4709-4730. 
%\bibitem{gia} M. Giaquinta, { Multiple integrals in the calculus of variations and nonlinear elliptic systems}, Princeton Univ. press, (1983).
%\bibitem{kat1}T. Kato,{\it Nonstationary flows of viscous and ideal fluids in $\Bbb R^3$}, J. Func. Anal. {\bf 9}, 296-305
%(1972).
%\bibitem{kat2}T. Kato and G.  Ponce, {\it  Commutator estimates and the Euler and Navier-Stokes equations},
%Comm. Pure Appl. Math. {\bf 41} (7), 891-907 (1988).
%\bibitem{koz}H. Kozono and Y. Taniuchi,  {\it Limiting case of the Sobolev inequality in BMO, with application to the Euler equations,} Comm. Math. Phys. {\bf 214},  no. 1, 191-200 (2000).
\bibitem{lem}P.G. Lemari\'{e}-Rieusset, {\it Euler Equations and Real Harmonic Analysis,}  Arch. Rational Mech. Anal. {\bf 204},  355-386 (2012).
\bibitem{lio}P. L. Lions, {\it Mathematical topics in fluid mechanics} Volume 1, Oxford University Press, New York, (1996).
%\bibitem{maj}A. Majda and A. Bertozzi,  {\it Vorticity and incompressible flow, } Cambridge Texts in Applied Mathematics, {\bf 27}, Cambridge University Press, Cambridge (2002).
\bibitem{pak}H.C. Pak and Y. J. Park, {\it Existence of solution for the Euler equations in a critical Besov space  $B^1_{\infty, 1} (\Bbb R^n)$,}
 Comm. P.D.E. {\bf 29},1149-1166  (2004).
 \bibitem{tay}M.E. Taylor, {\it Pseudodifferential operators, paradifferential operators, and layer potentials}, AMS, 
vol {\bf 81} (2000).
 \bibitem{tri}H.  Triebel, { Theory of function spaces}. Monographs in mathematics Vol {\bf 84},  Birkh\"auser, Basel 1992.
 \bibitem{vis1}M.  Vishik, {\it Hydrodynamics in Besov spaces,} Arch. Rat. Mech. Anal. {\bf 145}, 197-214 (1998).
\bibitem{vis2}M.  Vishik, {\it Incompressible flows of an ideal fluid with vorticity in borderline spaces of
Besov type,} Ann. Sci. \'{E}cole Norm. Sup. {\bf 32}(6), 769-812 (1999).
%\bibitem{yud}V. Yudovich,{\it  Nonstationary flow of an ideal incompressible liquid,} Zh. Vych. Mat. {\bf 3},
%1032-1066 (1963).
 \end{thebibliography}
 \end{document}